\documentclass[10pt]{amsart}

\pdfoutput=1

\usepackage[T1]{fontenc}
\usepackage{lmodern}

\usepackage[text={420pt,660pt},centering,marginparwidth=75pt,marginparsep=8pt]{geometry}
\usepackage{marginfix}
\makeatletter
\@mparswitchfalse
\makeatother

\usepackage{esint,amssymb} 
\usepackage{graphicx}
\usepackage{MnSymbol}
\usepackage{mathtools} 
\usepackage[colorlinks=true, pdfstartview=FitV, linkcolor=blue, citecolor=blue, urlcolor=blue,pagebackref=false]{hyperref}
\usepackage{orcidlink}
\usepackage{microtype}

\usepackage{bm}
\usepackage{dsfont}
\usepackage{mathrsfs}
\usepackage{xcolor}

\usepackage{caption}

\newtheorem{proposition}{Proposition}
\newtheorem{theorem}[proposition]{Theorem}
\newtheorem{lemma}[proposition]{Lemma}
\newtheorem{corollary}[proposition]{Corollary}

\theoremstyle{remark}

\theoremstyle{definition}

\numberwithin{equation}{section}
\numberwithin{proposition}{section}

\newcommand{\N}{\mathbb{N}}

\newcommand{\R}{\mathbb{R}}

\newcommand{\ep}{\varepsilon}

\renewcommand{\le}{\leqslant}
\renewcommand{\ge}{\geqslant}

\renewcommand{\geq}{\geqslant}

\renewcommand{\subset}{\subseteq}
\renewcommand{\bar}{\overline}

\newcommand{\Ll}{\left}
\newcommand{\Rr}{\right}
\renewcommand{\d}{\mathrm{d}}
\newcommand{\dr}{\partial}
\newcommand{\1}{\mathbf{1}}

\newcommand{\mcl}{\mathcal}
\newcommand{\msf}{\mathsf}
\newcommand{\mfk}{\mathfrak}
\newcommand{\msc}{\mathscr}
\newcommand{\al}{\alpha}
\newcommand{\be}{\beta}
\newcommand{\ga}{\gamma}
\newcommand{\de}{\delta}
\newcommand{\si}{\sigma}

\DeclareMathOperator{\tr}{tr}
\DeclareMathOperator{\supp}{supp}

\newcommand{\E}{\mathbb{E}}
\renewcommand{\P}{\mathbb{P}}
\DeclareMathOperator{\Cov}{Cov}
\DeclareMathOperator{\Var}{Var}
\DeclareMathOperator{\Law}{Law}
\newcommand{\norm}[1]{\lVert #1 \rVert}

\newcommand{\Is}{\mathrm{I}}   
\newcommand{\Sp}{\mathsf{S}}   
\newcommand{\Ga}{\mathsf{G}}   
\newcommand{\Id}{\mathrm{Id}}
\newcommand{\epsav}{\E_{\ep}}  
\newcommand{\cP}{\mcl P}       
\newcommand{\cQ}{\mcl Q}       
\newcommand{\cE}{\mcl E}       
\newcommand{\cR}{\mcl R}       
\newcommand{\cF}{\mcl F}       
\newcommand{\cG}{\mcl G}       
\newcommand{\cH}{\msf h}       
\newcommand{\cS}{\msc S}       

\newcommand{\sfC}{\msf C}
\newcommand{\sfD}{\msf D}

\newcommand{\Ball}{\mathbb{B}}
\newcommand{\Cube}{\mathbb{K}}
\newcommand{\Dom}{\mathbb{D}}
\newcommand{\Sph}{\mathbb{S}}
\renewcommand{\complement}{{\msf c}}

\begin{document}
\raggedbottom

\author[A.\ Ali]{Ayman Ali\,\orcidlink{0009-0009-0734-220X}}
\address[Ayman Ali]{Hebrew University of Jerusalem, Jerusalem, Israel}

\author[J.\ Boursier]{Jeanne Boursier}
\address[Jeanne Boursier]{Department of Mathematics, MIT, Cambridge MA, USA}

\author[H.-B.\ Chen]{Hong-Bin Chen\,\orcidlink{0000-0001-6412-0800}}
\address[Hong-Bin Chen]{NYU-ECNU Institute of Mathematical Sciences, NYU Shanghai, China}

\author[F.-H.\ Ho]{Fu-Hsuan Ho\,\orcidlink{0000-0003-3178-5472}}
\address[Fu-Hsuan Ho]{Weizmann Institute of Science, Rehovot, Israel}

\author[V.\ Issa]{Victor Issa\,\orcidlink{0009-0009-1304-046X}}
\address[Victor Issa]{LPSM, Sorbonne University, Paris, France}

\author[J.-C.\ Mourrat]{Jean-Christophe Mourrat\,\orcidlink{0000-0002-2980-725X}}
\address[Jean-Christophe Mourrat]{Department of Mathematics, ENS Lyon and CNRS, Lyon, France}

\date{}

\title{Free energy of the Hopfield model}

\begin{abstract}
We identify the limit of the free energy of the Hopfield model with independent symmetric Bernoulli patterns, for every temperature and loading ratio. A Hubbard--Stratonovich transformation reduces the problem to that of identifying the limit free energy of a bipartite spin glass with $\pm 1$ and spherical spins. One key difficulty is that the patterns, which become the coupling coefficients in the bipartite model, are Bernoulli instead of Gaussian, and the limit free energy is sensitive to the law of these coefficients. We show that the model is in fact equivalent to one involving Gaussian patterns and also subject to a finite number of Bernoulli random fields.  We then solve this model by adapting the techniques recently introduced to identify the limit free energies of spin glasses with either $\pm 1$ or spherical spins.
\end{abstract}

\maketitle


%
%
%
%
%
%
\section{Introduction}
\label{s.intro}

The Hopfield model \cite{hopfield1982neural, little1974existence, pastur1978theory} describes the behavior of a collection of $N$ binary variables $\si = (\si_1, \ldots, \si_N) \in \{-1,1\}^N$ representing neurons that interact with a given collection of patterns $x^1,\ldots, x^{M_N} \in \{-1,1\}^N$ meant to be memorized. To be precise, we give ourselves a sequence $(M_N)_{N \in \N}$ of nonnegative integers, and let $(x_i^\mu)_{i \le N, \mu \le M_N}$ be independent random variables, each uniform on $\{-1,1\}$; we call such random variables Bernoulli variables. We assume that
\begin{equation}
\label{e.loading}
\lim_{N \to \infty} \frac{M_N}{N} = \al \in (0,\infty),
\end{equation}
and call the parameter $\al$ the loading ratio. We write $X_N := (x_i^\mu)_{i \le N, \mu \le M_N}$ for the $N$-by-$M_N$ matrix of patterns, so that 
\begin{equation*}  
|X_N^\top \si|^2 = \sum_{\mu = 1}^{M_N} \Ll(\sum_{i = 1}^{N} x_i^\mu \si_i \Rr)^2 =  \sum_{i,j = 1}^N \sum_{\mu = 1}^{M_N} x^\mu_i x^\mu_j \si_i \si_j,
\end{equation*}
and we let $P_N$ denote the uniform probability measure on $\{-1,1\}^N$. For every inverse temperature $\be \ge 0$, we define the random free energy $F_N(\be)$ and its expectation $\bar F_N(\be)$ by
\begin{equation}
\label{e.def.free-energy}
F_N(\be) := \frac 1 N \log \int \exp \Ll( \frac{\be}{2N} |X_N^\top \si|^2 \Rr) \, \d P_N(\si)
\qquad \text{and} \qquad
\bar F_N(\be) := \E F_N(\be).
\end{equation}
The idea of the model is that the associated Gibbs measure represents the equilibrium state of the set of neurons, with the synaptic weight between neurons $i$ and $j$ being fixed by the patterns as $\sum_{\mu = 1}^{M_N} x^\mu_i x^\mu_j$, in agreement with Hebb's rule \cite{hebb1949organization}. 

Our main result is an explicit variational formula for the limit of the free energy of this model. We first describe the ingredients featured in this formula.
For $a \in [0,1)$, let $\cP([0,a])$ denote the set of Borel probability measures on $[0,a]$, and let
\begin{equation*}
\cP_{<1} := \bigcup_{a \in [0,1)} \cP([0,a])
\end{equation*}
be the set of probability measures on $[0,1)$ with support bounded away from $1$. For $\boldsymbol{\zeta} \in \cP_{<1}$, we define $D_{\boldsymbol{\zeta}} : [0,1) \to (0,\infty)$ and the function $\cE : \cP_{<1} \to \R$ by
\begin{equation}
\label{e.def.entropy}
D_{\boldsymbol{\zeta}}(x) := \int_x^1 \boldsymbol{\zeta}([0,s]) \, \d s 
\qquad \text{and} \qquad 
\cE(\boldsymbol{\zeta}) := \frac 1 2 \int_0^1 \Ll( \frac{1}{D_{\boldsymbol{\zeta}}(x)} - \frac{1}{1-x} \Rr) \, \d x.
\end{equation}
If $\supp \boldsymbol{\zeta} \subset [0,a]$ with $a < 1$, then $D_{\boldsymbol{\zeta}}(x) = 1 - x$ for $x \in [a,1)$, so the integrand in \eqref{e.def.entropy} vanishes on $[a,1)$. We will see in Lemma~\ref{l.spherical-duality} that $\cE$ is the convex conjugate of the spherical initial condition of \cite{chen2026spherical}, and that it is strictly convex.

For $a \ge 0$ and $\boldsymbol{\zeta} \in \cP_{<1}$, let $\phi_{a,\boldsymbol{\zeta}} : [0,1] \times \R \to \R$ be the solution of the Parisi equation
\begin{equation}
\label{e.scalar-PDE}
-\dr_s \phi_{a,\boldsymbol{\zeta}} = a \Ll( \dr_x^2 \phi_{a,\boldsymbol{\zeta}} + \boldsymbol{\zeta}([0,s]) (\dr_x \phi_{a,\boldsymbol{\zeta}})^2 \Rr), \qquad \phi_{a,\boldsymbol{\zeta}}(1,x) = \log \cosh x.
\end{equation}
When $\boldsymbol{\zeta}$ is finitely supported, the function $s \mapsto \boldsymbol{\zeta}([0,s])$ is a step function, and the solution of \eqref{e.scalar-PDE} is given by a finite recursion of Gaussian expectations; for a general $\boldsymbol{\zeta}$, the solution is the limit of these, see Section~\ref{s.scalar}. We define $\msc P_a : \cP_{<1} \times \R \to \R$ and $\cR : (0,\infty) \to \R$ by
\begin{equation}
\label{e.def.Parisi-radius}
\msc P_a(\boldsymbol{\zeta},h) := \phi_{a,\boldsymbol{\zeta}}(0,h)
\qquad \text{and} \qquad
\cR(r) := \log r - \frac{r^2-1}{2}.
\end{equation}
In particular, $\msc P_0(\boldsymbol{\zeta},h) = \log \cosh h$ for every $\boldsymbol{\zeta}$. The function $\msc P_a(\cdot, h)$ is the Parisi functional of the $\pm 1$ spins; the function $\cR$ is the exponential rate per pattern, as $N \to \infty$, of the density at $r$ of the normalized radius $|z|/\sqrt{M_N}$ of a standard Gaussian vector $z \in \R^{M_N}$, see \eqref{e.radial-density} and \eqref{e.radial-density-limit}.
\begin{theorem}[Free energy of the Hopfield model]
\label{t.hopfield}
Under \eqref{e.loading}, for every $\be > 0$, the limit
\begin{equation}
\label{e.def.FB}
F(\al,\be) := \lim_{N \to \infty} \bar F_N(\be)
\end{equation}
exists and is given by
\begin{equation}
\label{e.main-formula}
F(\al,\be) = \sup_{\substack{r > 0 \\ 0 \le h \le \be}} \ \inf_{\boldsymbol{\zeta} \in \cP_{<1}} \Ll\{ \al \cR(r) - \frac{h^2}{2\be} + \msc P_{\al \be r^2/2}(\boldsymbol{\zeta}, h) + \al \cE(\boldsymbol{\zeta}) \Rr\}.
\end{equation}
Moreover, $F_N(\be)$ converges to $F(\al,\be)$ almost surely and in $L^1$ as $N$ tends to infinity.
\end{theorem}

The next proposition collects the properties of the variational problem in \eqref{e.main-formula}. We stress that the uniqueness statement concerns the minimizing $\boldsymbol{\zeta}$ for a fixed pair $(r,h)$; we do not claim that the maximizing pair $(r,h)$ is unique.

\begin{proposition}[Structure of the variational problem]
\label{p.structure}
Let $\al, \be > 0$, $r > 0$ and $h \in \R$, and set $a := \al \be r^2/2$.
\begin{enumerate}
\item For every $b \in [0,1)$, the function $\boldsymbol{\zeta} \mapsto \msc P_a(\boldsymbol{\zeta},h) + \al \cE(\boldsymbol{\zeta})$ is continuous on $\cP([0,b])$ for the topology of weak convergence, and it is strictly convex on $\cP_{<1}$.
\item The infimum of this function over $\boldsymbol{\zeta} \in \cP_{<1}$ is attained at a unique $\boldsymbol{\zeta}$, and this minimizer belongs to $\cP([0, \be r^2/(1+\be r^2)])$.
\item The supremum in \eqref{e.main-formula} is attained, and its value is unchanged if the constraint $0 \le h \le \be$ is replaced by $h \in \R$.
\end{enumerate}
\end{proposition}
The starting point of our analysis is the observation (often called the Hubbard--Stratonovich transformation) that, denoting by $\boldsymbol{\gamma}_d$ the standard Gaussian measure on $\R^d$, we have
\begin{equation}  
\label{e.hubbard}
\bar F_N(\beta) = \frac 1 N \E \log \int \exp \Ll( \sqrt{\frac \beta N} \si^\top X_N z\Rr)  \, \d \boldsymbol{\gamma}_{M_N}(z) \, \d P_N(\si).
\end{equation}
This turns the Hopfield model into a bipartite spin glass. By decomposing the Gaussian variable into a radial and a spherical part, we can reduce the problem to that of identifying the limit free energy of a bipartite spin glass involving both $\pm 1$ and spherical spins. The identification of the limit free energy of general multi-species spin glasses has recently witnessed significant progress, with the cases of spherical and of $\pm 1$ spins being now fully resolved in \cite{chen2026spherical} and \cite{boursier2026free, chen2026ising} respectively. It may thus seem that, in order to fully resolve the problem, one only needs to combine these results and extend them to a mixed $\pm 1$/spherical setting.

This is however not the case. Indeed, the works mentioned above (and most works on the topic) assume that the coupling coefficients, in our case the entries of the matrix $X_N$, are Gaussian. In the context of models with only $\pm 1$ spins (or more generally when the reference measure is of the form $\boldsymbol{\pi}^{\otimes N}$ with $\boldsymbol{\pi}$ of compact support), this is not a significant restriction, as one can show a universality property of the limit free energy with respect to the law of the disorder: any choice of i.i.d.\ centered coefficients with unit variance yields the same limit free energy as the Gaussian case \cite{carmona2006universality, chatterjee2005simple}. It is however important to stress that this universality property holds only for terms that are of order two or more in the spin variables; for linear terms, if for instance we take i.i.d.\ random variables $(h_i)$ and consider
\begin{equation}  
\label{e.product.one-spin}
\frac 1 N \E \log \int \exp\Ll(\beta \sum_{i = 1}^N{h_i \sigma_i} \Rr)\, \d P_N(\sigma) = \E \log \cosh(\beta h_1),
\end{equation}
it is clear that this quantity depends on the full law of $h_1$ and not just on its first two moments. 

For spherical spin glasses, universality has been shown under suitable moment conditions on the disorder \cite{kim2026sharp, sawhney2024free}. The situation is more subtle in this case, and the main difficulty is that a handful of coordinates of a point on the sphere of radius $\sqrt{N}$ in $\R^N$ can take values of order $\sqrt{N}$, which suffices on its own to contribute meaningfully to the energy value. The key mechanism that allows universality results to hold despite this possibility is the rotational invariance of the model. And indeed, unlike in \eqref{e.product.one-spin}, for spherical models universality holds even for terms that are linear in the spin variables. As a simple illustration of this, notice that if we replace the measure $P_N$ by the uniform measure over the $N$-dimensional sphere of radius $\sqrt{N}$ in \eqref{e.product.one-spin}, then by rotational invariance, the integral only depends on the $\ell^2$ norm of the vector $h = (h_i)_{1 \le i \le N}$, whose asymptotic behavior depends only on the second moment of $h_1$. 

In the context of \eqref{e.hubbard} with mixed $\pm 1$ and spherical variables, the arguments for universality no longer apply because the spherical variable prevents a direct argument using boundedness of the coordinates, while the $\pm 1$ variables break rotational invariance. For instance, the vector $z$ can concentrate on a single coordinate $z_\mu$ to select the pattern $x^\mu$, which in turn acts on the $\sigma$ variables as an external field, yielding a contribution that is sensitive to the full law of the pattern coordinates. 

To circumvent this difficulty, we proceed as follows. To start with, we identify a sufficiently large constant $L$ so that, at the price of a vanishing error, we can restrict the integration variable~$z$ in \eqref{e.hubbard} to the ball $|z| \le L \sqrt{N}$. We then decompose the coordinates $\mu$ of $z$ according to whether $|z_\mu| \le \delta \sqrt{N}$ or not, for some small $\delta > 0$. Notice that the set $S$ of coordinates $\mu$ such that $|z_\mu| > \delta \sqrt{N}$ must be of bounded cardinality, due to the global restriction $|z| \le L \sqrt{N}$. Using that the number of possible choices for the set $S$ is only polynomial in $N$, the invariance under permutations of the coordinates, and a concentration argument, we can reduce the problem to that of studying the limit free energy of all models with an imposed fixed set $S$, and then take a supremum. For such models, we show using an interpolation argument that one can replace all the Bernoulli patterns indexed by $\mu \notin S$ by Gaussian patterns. The error in this replacement can be made arbitrarily small in the limit $N \to \infty$ by choosing $\delta > 0$ sufficiently small. This brings us to the study of a Gaussian model with a fixed finite number of Bernoulli patterns, which we now define precisely.

Let $G = (g_i^\mu)_{i \le N, \mu \le M}$ be a matrix with independent standard Gaussian entries. For $k \in \N \cup \{0\}$, let $(\eta_i^j)_{i \le N, j \le k}$ be Bernoulli variables, independent of each other and of $G$. For a fixed $\be > 0$, we define $\bar F^{(k)}_{N,M} : \R^k \to \R$ by
\begin{equation}
\label{e.def.hybrid}
\bar F^{(k)}_{N,M}(y) := \frac 1 N \E \log \int \exp \Ll( \frac{\be}{2N} |G^\top \si|^2 + \sqrt \be \sum_{i=1}^N \sum_{j=1}^k y_j \eta_i^j \si_i \Rr) \, \d P_N(\si).
\end{equation}
We understand that $\R^0$ consists of a single point and that an empty sum is zero, so that $\bar F^{(0)}_{N,M}$ is the expected free energy of a Hopfield model with Gaussian patterns. To lighten the notation, we do not display the dependence of $\bar F^{(k)}_{N,M}$ on $\be$. Similarly, its limit $F^{(k)}$ defined in \eqref{e.hybrid-formula} below depends on $\al$ and $\be$, but this is not reflected in the notation.
For every $y \in \R^k$, we define the random field
\begin{equation}
\label{e.def.field-law}
\cH_y := \sqrt \be \sum_{j=1}^k y_j \ep_j,
\end{equation}
where $(\ep_j)_{j \ge 1}$ are independent Bernoulli variables. In particular, for $k = 0$ we have $\cH_y = 0$.

\begin{theorem}[Gaussian patterns with finitely many Bernoulli fields]
\label{t.hybrid}
Under \eqref{e.loading}, for every $k \in \N \cup \{0\}$ and $\be > 0$, the functions $\bar F^{(k)}_{N,M_N}$ converge locally uniformly on $\R^k$ to
\begin{equation}
\label{e.hybrid-formula}
F^{(k)}(y) := \sup_{r > 0} \Ll\{ \al \cR(r) + \inf_{\boldsymbol{\zeta} \in \cP_{<1}} \Ll[ \E \, \msc P_{\al \be r^2/2}(\boldsymbol{\zeta},\cH_y) + \al \cE(\boldsymbol{\zeta}) \Rr] \Rr\}.
\end{equation}
For every $r > 0$ and $y \in \R^k$, the infimum over $\boldsymbol{\zeta}$ in \eqref{e.hybrid-formula} is attained at a unique $\boldsymbol{\zeta}$, which belongs to $\cP([0,\be r^2/(1+\be r^2)])$. The supremum over $r$ is attained, and for $y$ varying in a bounded subset of $\R^k$, the corresponding maximizing $r$ lies in a common compact subset of $(0,\infty)$.
\end{theorem}
For $k \in \N$, we define
\begin{equation}
\label{e.def.Ak}
A_k(\al,\be) := \sup_{y \in \R^k} \Ll\{ F^{(k)}(y) - \frac{|y|^2}{2} \Rr\}.
\end{equation}
The quadratic term is the cost coming from the Gaussian measure of the $k$ retained patterns associated with the set $S$ discussed earlier. Since $\cH_{(y,0)} = \cH_y$ for every $y \in \R^k$, we have $F^{(k+1)}((y,0)) = F^{(k)}(y)$, and the sequence $(A_k(\al,\be))_{k \in \N}$ is nondecreasing. The strategy outlined above yields the following identification of the limit free energy of the Hopfield model.

\begin{theorem}[Reduction to Gaussian patterns with finitely many Bernoulli fields]
Under \eqref{e.loading}, for every $\be > 0$, the limit $F(\al,\be)$ in \eqref{e.def.FB} exists and
\begin{equation}
\label{e.sup-Ak}
F(\al,\be) = \sup_{k \in \N} A_k(\al,\be).
\end{equation}
Moreover, the supremum in \eqref{e.sup-Ak} is attained at $k = 1$: for every $k \in \N$ and $y \in \R^k$, we have $F^{(k)}(y) \le F^{(1)}(|y|)$, and therefore $A_k(\al,\be) = A_1(\al,\be)$.
\label{t.sup-Ak}
\end{theorem}
In fact, one can refine the analysis we provide here to show that for every $y \in \R^k$ with at least two non-zero coordinates, we have $F^{(k)}(y) < F^{(1)}(|y|)$. We expect this to be related to the fact that a sample from the Gibbs measure has a non-vanishing overlap with at most one pattern. 

Let $F^\Ga(\al,\be;h)$ denote the limit of the expected free energy of the model with Gaussian patterns and the uniform external field $h$, that is, of the right-hand side of \eqref{e.def.hybrid} with the double sum $\sqrt \be \sum_{i=1}^N \sum_{j=1}^k y_j \eta_i^j \si_i$ replaced by $h \sum_{i = 1}^N \sigma_i$. The change of variables $\si_i \mapsto \eta^1_i \si_i$, $g_i^\mu \mapsto \eta_i^1 g_i^\mu$ preserves the law of $G$ and turns the field $\sqrt \be \, y_1 \eta_i^1 \si_i$ into $\sqrt \be \, y_1 \si_i$. Hence, using Theorem~\ref{t.hybrid} with $k = 1$ and the fact that  $\msc P_a(\boldsymbol{\zeta},\cdot)$ is even, we have
\begin{equation}
\label{e.gaussian-formula}
F^\Ga(\al,\be;h) = F^{(1)}(h/\sqrt \be) = \sup_{r > 0} \Ll\{ \al \cR(r) + \inf_{\boldsymbol{\zeta} \in \cP_{<1}} \Ll[ \msc P_{\al \be r^2/2}(\boldsymbol{\zeta},h) + \al \cE(\boldsymbol{\zeta}) \Rr] \Rr\}.
\end{equation}
We write $F^\Ga(\al,\be) := F^\Ga(\al,\be;0)$ for the limit free energy of the model with Gaussian patterns and no external field.

\begin{corollary}[Formula with one external field]
\label{c.one-pattern}
For every $\al, \be > 0$, we have
\begin{equation}
\label{e.one-pattern}
F(\al,\be) = A_1(\al,\be) = \sup_{0 \le h \le \be} \Ll\{ F^\Ga(\al,\be;h) - \frac{h^2}{2\be} \Rr\}.
\end{equation}
Moreover, for every $\al > 0$ and $\be \in (0,1]$, 
\begin{equation}
\label{e.beta-at-most-one}
F(\al,\be) = F^\Ga(\al,\be),
\end{equation}
and the maximum in \eqref{e.one-pattern} is attained only at $h = 0$.
\end{corollary}
We can verify that the limit free energies $F$ and $F^\Ga$ corresponding respectively to the Bernoulli and the Gaussian patterns differ when $\alpha$ is sufficiently small and $\beta$ is sufficiently large.
\begin{proposition}[Strict separation of the two models]
\label{p.separation}
For every $\al, \be > 0$, we have
\begin{equation}
\label{e.two-bounds}
F^\Ga(\al,\be) \le \frac \be 2 \Ll( \sqrt{2/\pi} + \sqrt \al \Rr)^2
\qquad \text{and} \qquad
F(\al,\be) \ge \frac \be 2 (1+\al) - \log 2.
\end{equation}
Consequently, if
\begin{equation}
\label{e.strict-range}
0 < \al < \frac{(\pi-2)^2}{8\pi} 
\qquad \text{and} \qquad 
\be > \frac{2 \log 2}{1 + \al - (\sqrt{2/\pi} + \sqrt \al)^2},
\end{equation}
then $F(\al,\be) > F^\Ga(\al,\be)$.
\end{proposition}
We expect that the regime where $F(\al,\be) > F^\Ga(\al,\be)$ (in other words, the regime where the choice $h = 0$ in \eqref{e.one-pattern} is suboptimal) corresponds to the situation where a sample from the Gibbs measure gets a non-vanishing overlap with one of the patterns. We do not believe the bound $(\pi-2)^2/(8\pi)$ in \eqref{e.strict-range} to be sharp.

A separate interesting phase transition is that of replica symmetry breaking. In the variational formula \eqref{e.main-formula}, this corresponds to asking whether the infimum is achieved at a probability measure~$\boldsymbol{\zeta}$ that is a Dirac mass. The fact that the model is replica symmetric when $(1 + \sqrt{\alpha}) \beta < 1$ was shown in \cite{talagrand1998rigorous}. For a model with an additional term $h \sum_{i = 1}^N x^1_i \sigma_i$ in the exponential in \eqref{e.def.free-energy}, we know from \cite[Theorem~10.7.1]{talagrand2011mean2} that for every $h \neq 0$ and $\beta > 1$, there exists $\alpha_0(h, \beta) > 0$ such that the model is replica symmetric for all $\alpha \le \alpha_0(h,\beta)$. In agreement with \cite{amit1987statistical}, we expect that for every $\al > 0$ sufficiently small, the model is replica symmetric for $\beta \in [0,(1+\sqrt{\al})^{-1}] \cup {(1+ c(\alpha) \sqrt{\alpha},  \beta_c(\alpha)]}$, whereas replica symmetry is broken for $\beta \in {((1+\sqrt{\al})^{-1}, 1+c(\al)\sqrt{\al})} \cup (\beta_c(\al),\infty)$, where $\lim_{\alpha \to 0} c(\alpha) \simeq 2.6$ and
\begin{equation*}  
\beta_c(\al) = \frac{3\sqrt{2\pi}}{4\sqrt\al}\exp\!\left(\frac1{2\al}\right)(1+o(1)) \qquad (\al \to 0).
\end{equation*}
For $\alpha > 0$ sufficiently small and $\beta = 1 + c(\alpha)\sqrt{\alpha}$, we expect that the supremum in \eqref{e.main-formula} admits at least two maximizers, one with $h > 0$ whose minimizing law is a Dirac mass, and one with $h = 0$ whose minimizing law is not a Dirac mass.

\subsection*{Proof ideas}
Some of the key ideas involved in the reduction of the original Hopfield model to a model with Gaussian patterns and a fixed number of additional Bernoulli random fields have already been touched upon, so we now discuss the strategy for analyzing the latter model. 

As already mentioned, the limit free energies of multi-species models with $\pm 1$ and with spherical spins have been identified in \cite{boursier2026free, chen2026ising} and in \cite{chen2026spherical} respectively. We adopt here the same approach as in~\cite{boursier2026free}, which was introduced in \cite{aronow2026quantitative} to obtain a rate of convergence of the free energy to its limit for the Sherrington--Kirkpatrick (SK) model with no external field. This method was also used in~\cite{ho2026convergence} to show the existence of the limit free energy of general mean-field spin glasses, and in \cite{chen2026formula} to give a representation of this limit in terms of a family of solutions of certain Hamilton--Jacobi equations.

The approach of \cite{aronow2026quantitative} requires the existence of transformations that leave the model invariant and act transitively on the space of configurations. For the SK model with no external field considered in \cite{aronow2026quantitative}, these transformations are, for each $\tau \in \{-1,1\}^N$, the maps $\sigma \mapsto (\tau_i \sigma_i)_{i \le N}$. However, for models with an external field, these transformations no longer leave the model invariant. So instead, the approach of \cite{boursier2026free} is to first work on a fixed ``slice'' of the model where, for each species involved, the number of $+1$ and $-1$ spins is imposed. Once the magnetization is thus constrained, the effect of the external field no longer depends on the configuration, and intra-species permutations act transitively on the model, thus allowing us to make use of the approach of \cite{aronow2026quantitative}. For more general models as in \cite{chen2026formula, ho2026convergence} where spins can take more than two values, this amounts to fixing the empirical distribution of the coordinates. The outcome is an identification of the limit free energy in terms of the solution to a Hamilton--Jacobi equation. The restriction to slices of fixed magnetization is also crucially used in \cite{boursier2026free} to show that the initial condition appearing in this equation is convex, as this allows us to write a variational formula for the solution known as the Hopf formula \cite{chen2025HJ}.

Coming back to the quantity in \eqref{e.def.hybrid}, we start by applying the Hubbard--Stratonovich transformation as in \eqref{e.hubbard} and by decomposing the Gaussian integration into a radial and a spherical part. In order to apply the analysis strategy discussed in the previous paragraphs, we freeze the coefficients $(\eta_i^j)_{i \le N, j \le k}$, and decompose the coordinates in $\{1,\ldots, N\}$ into $2^k$ types, with the coordinate $i$ belonging to type $\ep \in \{-1,1\}^k$ when $(\eta^j_i)_{j \le k} = \ep$. In this way, we obtain a model with $2^k + 1$ types (the additional type being the spherical variable), where all the coordinates of a given spin type are exposed to a common external field. This puts us in a position to transfer the arguments of~\cite{boursier2026free}, including the convexity of the relevant initial condition, to our mixed $\pm 1$/spherical setting. Once an explicit variational formula is obtained, we can collapse back the number of types to simplify the formula. For spin-glass models with convex covariance function in a random external field with i.i.d.\ coordinates, it is well known that the Parisi functional is obtained from that of the model with a deterministic field $h$ by averaging over the law of $h$, see for instance \cite[Theorem~14.5.1]{talagrand2011mean2}. The average over the random field~$\cH_y$ in \eqref{e.hybrid-formula} plays exactly this role here, although our proof of this fact is more indirect. Since the model is bipartite, one might expect its limit free energy to involve two probability measures, one for the $\pm 1$ spins and one for the spherical spins, whereas only one appears in \eqref{e.hybrid-formula}. This is because the nonlinearity in the Hamilton--Jacobi equation, namely the covariance function of the bipartite model, is bilinear in the two species. This structure fits well with the convex-analysis arguments underlying the Hopf formula, and allows us to optimize explicitly over one of these two probability measures, for each fixed value of the other.

\subsection*{Related works} 
The Hopfield model was introduced in \cite{hopfield1982neural} as a model of associative memory, with synaptic weights prescribed by Hebb's rule \cite{hebb1949organization}; a closely related model with synchronous dynamics had been considered earlier in \cite{little1974existence}. The model was also proposed independently in \cite{pastur1977rigorously, pastur1978theory} as a simple model of spin glass (for a fixed number of patterns). One important question is whether each pattern lies close to a local minimum of the energy function, so that a relaxation dynamics started in the vicinity of a pattern recovers it. For a fixed number of patterns, this question was answered positively in \cite{amit1985spin}. In the regime of extensive loading, the replica method was used in \cite{amit1985storing, amit1987statistical} to derive a replica-symmetric phase diagram. Corrections to these predictions due to replica symmetry breaking were computed in \cite{crisanti1986saturation, steffan1994replica}.

On the mathematical side, the regime $M_N/N \to 0$ is well understood, see \cite{bovier1998review} and the references therein. At positive loading ratio, the replica-symmetric equations of \cite{amit1987statistical} were derived in \cite{pastur1994replica} under the assumption of concentration of the overlap. In the regime of small $\al$ and large $\be$, the concentration of the Gibbs measure near the patterns was analyzed in \cite{bovier1997retrieval, bovier1995gibbs}. The validity of the replica-symmetric solution was then proved in \cite{talagrand1998rigorous, talagrand2000exponential} in a range of parameters that includes the high-temperature region $\be(1+\sqrt\al) < 1$, as well as a region of small $\al$ and large $\be$ for the model with the additional term $h \sum x^1_i \sigma_i$, $h \neq 0$, discussed earlier; see also the chapters on the Hopfield model in \cite{talagrand2011mean1, talagrand2011mean2}. The concentration of the random free energy around its expectation was already known for all $\al$ and $\be$ \cite{bovier1995gibbs, shcherbina1993free}, but to the best of our knowledge, the convergence of the expectation itself was only known in the regimes described above, or under a concentration assumption on the order parameters \cite{agliari2026infinite}.

The Hubbard--Stratonovich representation \eqref{e.hubbard} turns the Hopfield model into a bipartite spin glass \cite{amit1985spin, barra2012equivalence}. Bipartite models are perhaps the most basic examples of models with a non-convex covariance function; for such models, the multi-species Parisi formula of \cite{barra2015multispecies, panchenko2015multispecies} is invalid in general \cite[Section~6]{mourrat2021nonconvex}. The last few years have witnessed a major effort to make progress on this point. For models with a single type of spins, an explicit link between the limit free energy and Hamilton--Jacobi equations was identified in \cite{mourrat2022parisi}, inspired by related observations in \cite{barra2013mean, barra2010breaking, guerra2001sum}. For vector and multi-species spin glasses, this point of view first made it possible to obtain an upper bound on the limit free energy in terms of the solution to the PDE \cite{mourrat2021nonconvex, mourrat2023upper}, as well as a critical-point representation of the limit free energy \cite{chen2026species, chen2025vector, ho2026convergence}. More recently, a full identification of the limit free energy was achieved for several model classes \cite{boursier2026free, chen2026formula, chen2026ising, chen2026spherical}, as was already discussed. 

For the Hopfield model with Gaussian patterns, the annealed region $\be(1+\sqrt\al) < 1$ was identified in \cite{barra2012glassy, barra2008ergodic}, and the replica-symmetric approximation was studied in \cite{barra2010replica, barra2012glassy}; see also~\cite{barra2018phase} for models with general laws of the patterns and of the hidden variables, and \cite{agliari2017neural} for models with a mix of Gaussian and Bernoulli patterns. As this paper was being completed, a variational formula for the limit free energy of the model with Gaussian patterns was announced in \cite{camilli2026solution}, with proofs deferred to a forthcoming paper. Their formula corresponds to that in \eqref{e.gaussian-formula} at $h = 0$ (without the reduction from two to one probability measure). Even more recently, the existence of the limit free energy for the Hopfield model with Bernoulli patterns was shown in \cite{barra2026existence}.

\subsection*{Organization of the paper}
Section~\ref{s.scalar} recalls the scalar initial conditions of \cite{boursier2026free, chen2026spherical}, computes the convex conjugate of the spherical one, and proves the concavity of the Parisi functional in the squared external field. Section~\ref{s.mixed} identifies the limit free energy of slices of $\pm 1$ spins coupled to a sphere, by transferring the arguments of \cite{boursier2026free}. Section~\ref{s.reduction} reduces the Hopf formula to an optimization over one overlap law and releases the magnetization constraints. Section~\ref{s.gaussian} proves Theorem~\ref{t.hybrid}. Section~\ref{s.comparison} contains the comparison between Bernoulli and Gaussian disorder under cutoffs, and Section~\ref{s.finite-pattern} combines it with the retention of finitely many patterns to prove Theorem~\ref{t.sup-Ak}, Proposition~\ref{p.structure}, as well as \eqref{e.def.FB}--\eqref{e.main-formula} and \eqref{e.one-pattern}. Section~\ref{s.consequences} proves the almost sure and $L^1$ convergence in Theorem~\ref{t.hopfield}, the identity \eqref{e.beta-at-most-one} together with the last assertion of Corollary~\ref{c.one-pattern}, and Proposition~\ref{p.separation}.

\subsection*{Notation} Throughout the paper, $C$ and $c$ denote positive finite constants whose values may change from one occurrence to the next, with $C$ large and $c$ small; indices such as in $C_\be$ or $C_{a,b,\be}$ indicate the parameters on which these constants are allowed to depend. For a matrix $A$, we write $\norm{A}_{\mathrm{op}}$ for its operator norm. For $d \in \N$, we write $\Sph_d := \{z \in \R^d : |z|^2 = d\}$ for the sphere of radius $\sqrt d$ in $\R^d$, and $\boldsymbol{\mu}_d$ for the uniform probability measure on $\Sph_d$.

\subsection*{Statement of AI use} AA obtained a proof of existence of the limit free energy, with assistance from GPT 6 Astra. He reached out to FHH to share his finding, who in turn reached out to the rest of us. HBC and JCM already had a detailed draft identifying the limit free energy of the Hopfield model with Gaussian patterns, and realized that identifying the limit free energy for Bernoulli patterns was also possible. We worked together to write the final manuscript, with assistance from GPT 6 Astra and Fable 5.1. We have manually verified all statements and proofs and polished the paper, and we take full responsibility for its contents. 

\subsection*{Acknowledgments} HBC acknowledges funding from the NYU Shanghai Start-Up Fund and support from the NYU--ECNU Institute of Mathematical Sciences at NYU Shanghai. FHH acknowledges support from Eliran Subag's ISF grant (No. 2055/21) and ERC grant (No. 101165541). VI acknowledges support from the FSMP postdoctoral program. JCM acknowledges support from the MSCA grant SLOHD (101203974).

%
%
%
%
%
%
\section{Scalar functions}
\label{s.scalar}

From now on, we call $\pm 1$ spins Ising spins. The limit free energy of the bipartite model studied in Section~\ref{s.mixed} is expressed, through a Hopf formula, in terms of two scalar functions of a path: the initial condition $g_m$ of Ising spins with prescribed magnetization $m$, introduced in \cite{boursier2026free}, and the spherical initial condition $\psi^\Sp$ of \cite{chen2026spherical}. In this section, we recall the properties of these functions that will be needed, and we compute the convex conjugate of $\psi^\Sp$. We also prove the properties of~$\msc P_a$, as a function of the external field, that will allow us to reduce several retained patterns to a single one. 

\subsection{Ising spins with prescribed magnetization}
\label{s.Ising-data}

Let $\cQ$ denote the set of right-continuous nondecreasing maps from $[0,1)$ to $\R_+ := [0,\infty)$, and for $r \in [1,\infty]$, let $\cQ_r := \cQ \cap L^r([0,1))$. We call the elements of $\cQ$ paths. For $\rho \in \cQ_\infty$, we set $\rho(1) := \lim_{u \uparrow 1} \rho(u)$, and we define the law of the path by
\begin{equation}
\boldsymbol{\nu}_\rho := \Law(\rho(U)), \qquad \text{where } U \text{ is uniform on } [0,1].
\label{e.def.path-law}
\end{equation}
Conversely, for $\boldsymbol{\zeta} \in \cP_{<1}$, we write $Q_{\boldsymbol{\zeta}} \in \cQ_\infty$ for the quantile function of $\boldsymbol{\zeta}$, which is the right-continuous nondecreasing map from $[0,1)$ to $[0,1)$ such that $\boldsymbol{\nu}_{Q_{\boldsymbol{\zeta}}} = \boldsymbol{\zeta}$. For $\rho \in \cQ_\infty$ and $T \ge \rho(1)$, let $\phi_\rho : [0,T] \times \R \to \R$ be the solution of the Parisi equation
\begin{equation}
\label{e.initial-PDE}
-\dr_v \phi_\rho = \dr_x^2 \phi_\rho + \boldsymbol{\nu}_\rho([0,v]) (\dr_x \phi_\rho)^2, \qquad \phi_\rho(T,x) = \log \cosh x.
\end{equation}
When $\boldsymbol{\nu}_\rho$ is finitely supported, the coefficient $v \mapsto \boldsymbol{\nu}_\rho([0,v])$ is a step function, and \eqref{e.initial-PDE} is solved by a finite recursion. Indeed, on an interval $[v_0, v_0 + \ell]$ on which this coefficient is constant equal to $a$, the Cole--Hopf transformation gives $\phi_\rho(v_0, \cdot) = \mcl T_{a,\ell} \phi_\rho(v_0 + \ell, \cdot)$, where
\begin{equation}
\label{e.Gaussian-recursion}
\mcl T_{a,\ell} f(x) :=
\begin{cases}
a^{-1} \log \E \exp \bigl( a f(x + \sqrt{2\ell} Z) \bigr) & \text{if } a > 0, \\
\E f(x + \sqrt{2 \ell} Z) & \text{if } a = 0,
\end{cases}
\end{equation}
and $Z$ is a standard Gaussian random variable. For a general bounded path, the solution is obtained by approximation with step paths, see for instance \cite[Section~6]{boursier2026free} or \cite{auffinger2015unique, jagannath2016dynamic}. For every bounded path $\rho$, $h \in \R$ and $m \in [-1,1]$, we define
\begin{equation}
\label{e.def.gm}
\psi^\Is(\rho;h) := T - \phi_\rho(0,h)
\qquad \text{and} \qquad
g_m(\rho) := \sup_{h \in \R} \{ \psi^\Is(\rho;h) + hm \}.
\end{equation}
These quantities do not depend on the choice of $T \ge \rho(1)$: on $[\rho(1),T]$ the coefficient in \eqref{e.initial-PDE} is one, and $\mcl T_{1,\ell} \log \cosh = \ell + \log \cosh$. 
The superscript $\Is$ stands for Ising: in the notation of \cite{boursier2026free}, $\psi^\Is(\rho;h)$ is the function $\psi^\circ(\boldsymbol{\nu}_\rho;h)$, the initial condition of the Hamilton--Jacobi equation for Ising spins in the external field $h$, and $g_m(\rho)$ is the initial condition for Ising spins with prescribed magnetization $m$. The function $g_m$ is the limit, as the number of spins tends to infinity, of the enriched free energy at time zero of a system of Ising spins constrained to have magnetization $m$, by \cite[Proposition~7.5]{boursier2026free}; see also Lemma~\ref{l.mixed-initial-limit} below. The extensions of $\psi^\Is$ and $g_m$ to $\cQ_1$ are provided by the Lipschitz estimates in Lemma~\ref{l.Ising-data}. Let $I : [-1,1] \to [0, \log 2]$ be the entropy
\begin{equation}
\label{e.def.I}
I(m) := \frac{1+m}{2} \log(1+m) + \frac{1-m}{2} \log(1-m),
\end{equation}
with the convention $0 \log 0 := 0$. This is the convex conjugate of $\log \cosh$.

\begin{lemma}[The Ising initial condition]
\label{l.Ising-data}
The functions $g_m$ and $\psi^\Is$ have the following properties.
\begin{enumerate}
\item For every $m \in [-1,1]$, the function $g_m$ is convex and $1$-Lipschitz on $\cQ_1$ for the $L^1$ norm, and $g_m(0) = I(m)$. For every $\rho \in \cQ_\infty$, there exists $z \in \cQ_\infty$ with values in $[m^2,1]$ such that, for every $\rho' \in \cQ_\infty$,
\begin{equation}
\label{e.gm-supporting}
g_m(\rho') \ge g_m(\rho) + \int_0^1 z(u) (\rho'(u) - \rho(u)) \, \d u.
\end{equation}
\item For every $\rho \in \cQ_1$, the function $m \mapsto g_m(\rho) - I(m)$ is even and convex on $[-1,1]$, and for every $m, m' \in [-1,1]$,
\begin{equation}
\label{e.gm-continuity}
|g_m(\rho) - g_{m'}(\rho) - I(m) + I(m')| \le 2 |\rho|_{L^1} |m - m'|.
\end{equation}
\item For every $\rho \in \cQ_1$, $h \in \R$, $a \ge 0$ and $\boldsymbol{\zeta} \in \cP_{<1}$,
\begin{equation}
\label{e.field-duality}
\psi^\Is(\rho;h) = \min_{m \in [-1,1]} \{ g_m(\rho) - hm \}
\qquad \text{and} \qquad
\msc P_a(\boldsymbol{\zeta},h) = a - \psi^\Is(a Q_{\boldsymbol{\zeta}};h).
\end{equation}
\item For every $m \in [-1,1]$, $a \ge 0$ and $b \in [0,1)$, the function $\boldsymbol{\zeta} \mapsto g_m(a Q_{\boldsymbol{\zeta}})$ is continuous on $\cP([0,b])$ for the topology of weak convergence. Moreover, for every $\boldsymbol{\zeta}, \boldsymbol{\zeta}' \in \cP([0,b])$ and $\theta \in [0,1]$, we have
\begin{equation}
\label{e.gm-law-concavity}
g_m(a Q_{\theta \boldsymbol{\zeta} + (1-\theta) \boldsymbol{\zeta}'}) \ge \theta g_m(a Q_{\boldsymbol{\zeta}}) + (1-\theta) g_m(a Q_{\boldsymbol{\zeta}'}).
\end{equation}
\end{enumerate}
\end{lemma}
\begin{proof}
The convexity, the Lipschitz estimate, and \eqref{e.gm-supporting} in part (1) follow from \cite[Proposition~6.10]{boursier2026free}, where a law is identified with its quantile path; the extension to $\cQ_1$ follows from the Lipschitz estimate. Part (2) and the first identity in \eqref{e.field-duality} follow from \cite[Lemma~6.7]{boursier2026free}. Since $\psi^\Is(0;h) = -\log \cosh h$ (take $T = 0$ in \eqref{e.initial-PDE} and \eqref{e.def.gm}), the definition \eqref{e.def.gm} gives $g_m(0) = \sup_{h \in \R} \{hm - \log \cosh h\} = I(m)$.

We now prove the second identity in \eqref{e.field-duality}. The case $a = 0$ follows from $\msc P_0(\boldsymbol{\zeta},h) = \log \cosh h = -\psi^\Is(0;h)$. For $a > 0$, the quantile path $a Q_{\boldsymbol{\zeta}}$ is bounded by $a$, and we may take $T = a$ in \eqref{e.initial-PDE}. The law of $a Q_{\boldsymbol{\zeta}}(U)$ is the image of $\boldsymbol{\zeta}$ under $x \mapsto ax$, so that $\boldsymbol{\nu}_{aQ_{\boldsymbol{\zeta}}}([0,as]) = \boldsymbol{\zeta}([0,s])$. The change of variables $v = as$ in \eqref{e.initial-PDE} therefore shows that $(s,x) \mapsto \phi_{aQ_{\boldsymbol{\zeta}}}(as,x)$ solves \eqref{e.scalar-PDE}, that is, $\phi_{aQ_{\boldsymbol{\zeta}}}(as,x) = \phi_{a,\boldsymbol{\zeta}}(s,x)$. Evaluating at $s = 0$ gives $\psi^\Is(aQ_{\boldsymbol{\zeta}};h) = a - \msc P_a(\boldsymbol{\zeta},h)$.

The inequality \eqref{e.gm-law-concavity} is \cite[Corollary~6.8]{boursier2026free}, applied to the image of $\boldsymbol{\zeta}$ under $x \mapsto ax$. Continuity follows from the $1$-Lipschitz estimate of part (1) and the identity 
\begin{equation*}
|a Q_{\boldsymbol{\zeta}} - a Q_{\boldsymbol{\zeta}'}|_{L^1} = a \, W_1(\boldsymbol{\zeta}, \boldsymbol{\zeta}'),
\end{equation*}
where $W_1$ is the Wasserstein distance of order one, since on $\cP([0,b])$, weak convergence is equivalent to convergence for $W_1$.
\end{proof}

\subsection{The spherical initial condition and its conjugate}
\label{s.spherical-data}

For $\rho \in \cQ_\infty$, we define $K_\rho : [0,\rho(1)] \to \R_+$ by
\begin{equation}
\label{e.def.K}
K_\rho(s) := \int_0^1 \bigl( \rho(1) - \max\{\rho(u), s\} \bigr) \, \d u.
\end{equation}
The function $K_\rho$ is nonincreasing, $K_\rho(0) = \rho(1) - \int_0^1 \rho(u) \, \d u$, and $K_\rho(\rho(1)) = 0$. The spherical initial condition $\psi^\Sp : \cQ_1 \to \R$ is the continuous extension of the function defined, for $\rho \in \cQ_\infty$, by
\begin{equation}
\label{e.def.psiS}
\psi^\Sp(\rho) := - \inf_{b > 2 K_\rho(0)} \Ll\{ \int_0^{\rho(1)} \frac{\d s}{b - 2 K_\rho(s)} + \frac{b - 1 - \log b}{2} - \rho(1) \Rr\}.
\end{equation}
This is the function $\psi^\circ_0$ of \cite[(2.24)--(2.25)]{chen2026spherical}, with external field zero. In particular, $\psi^\Sp(0) = 0$. We define the convex conjugate $\psi^{\Sp,*} : \cQ_\infty \to \R \cup \{+\infty\}$ of $\psi^\Sp$ by
\begin{equation}
\label{e.def.psiS-conj}
\psi^{\Sp,*}(p) := \sup_{\rho \in \cQ_\infty} \Ll\{ \int_0^1 \rho(u) p(u) \, \d u - \psi^\Sp(\rho) \Rr\}.
\end{equation}
The next lemma collects the regularity of $\psi^\Sp$ and the characterization of the minimizer in \eqref{e.def.psiS} from \cite{chen2026spherical}. The bound on the minimizer will determine the compact set of overlap laws to which the variational problems can be restricted.

\begin{lemma}[Derivative of the spherical initial condition]
\label{l.spherical-derivative}
The function $\psi^\Sp$ is convex and $1$-Lipschitz on $\cQ_1$ for the $L^1$ norm. For every $\rho \in \cQ_\infty$, the infimum in \eqref{e.def.psiS} is attained at a unique $b$, which is characterized by
\begin{equation}
\label{e.multiplier}
1 = \frac 1 b + 2 \int_0^{\rho(1)} \frac{\d s}{(b - 2 K_\rho(s))^2}.
\end{equation}
For this $b$, we define $p_\rho \in \cQ_\infty$ by
\begin{equation}
\label{e.spherical-derivative}
p_\rho(u) := 2 \int_0^{\rho(u)} \frac{\d s}{(b - 2 K_\rho(s))^2}
\qquad \text{and} \qquad
p_\rho(1) = 1 - \frac 1 b \le \frac{2 \rho(1)}{1 + 2 \rho(1)}.
\end{equation}
For every $\rho' \in \cQ_\infty$, we have
\begin{equation}
\label{e.spherical-subgradient}
\psi^\Sp(\rho') \ge \psi^\Sp(\rho) + \int_0^1 p_\rho(u) (\rho'(u) - \rho(u)) \, \d u.
\end{equation}
Consequently, for every $\rho \in \cQ_\infty$,
\begin{equation}
\label{e.spherical-biconjugacy}
\psi^\Sp(\rho) = \max_{\substack{p \in \cQ_\infty \\ 0 \le p \le 1}} \Ll\{ \int_0^1 \rho(u) p(u) \, \d u - \psi^{\Sp,*}(p) \Rr\}.
\end{equation}
\end{lemma}
\begin{proof}
The convexity and $1$-Lipschitz continuity of $\psi^\Sp$ on $\cQ_1$ follow from \cite[Proposition~2.8 and Lemma~2.9]{chen2026spherical}, with external field zero. We fix $\rho \in \cQ_\infty$. Differentiating the function minimized in \eqref{e.def.psiS} with respect to $b > 2K_\rho(0)$, we obtain
\begin{equation*}
\frac{\d}{\d b} \Ll( \int_0^{\rho(1)} \frac{\d s}{b - 2 K_\rho(s)} + \frac{b - 1 - \log b}{2} - \rho(1) \Rr) = \frac 1 2 \Ll( 1 - \frac 1 b - 2 \int_0^{\rho(1)} \frac{\d s}{(b - 2K_\rho(s))^2} \Rr).
\end{equation*}
The expression subtracted from $1$ on the right-hand side is continuous and strictly decreasing in $b$, and it tends to $0$ as $b \to \infty$. It also tends to $+\infty$ as $b \downarrow 2K_\rho(0)$. Indeed, for $s \in [0,\rho(1)]$, \eqref{e.def.K} gives
\begin{equation*}
K_\rho(0) - K_\rho(s) = \int_0^1 \bigl( \max\{\rho(u),s\} - \rho(u) \bigr) \, \d u = \int_0^1 (s - \rho(u))_+ \, \d u \le s.
\end{equation*}
If $\rho(1) > 0$, this bound yields
\begin{equation*}
2 \int_0^{\rho(1)} \frac{\d s}{(b - 2K_\rho(s))^2} \ge \frac{1}{b - 2K_\rho(0)} - \frac{1}{b - 2K_\rho(0) + 2\rho(1)}.
\end{equation*}
If $\rho(1) = 0$, then $K_\rho(0) = 0$ and $1/b$ tends to $+\infty$ as $b \downarrow 0$. Therefore the function minimized in \eqref{e.def.psiS} is strictly convex and has a unique minimizer, characterized by \eqref{e.multiplier}.

To bound this minimizer, we use $K_\rho(s) \le \rho(1) - s$ for $s \in [0,\rho(1)]$, which follows from \eqref{e.def.K}. This gives
\begin{equation*}
\frac{1}{1 + 2\rho(1)} + 2 \int_0^{\rho(1)} \frac{\d s}{(1 + 2\rho(1) - 2K_\rho(s))^2} \le \frac{1}{1 + 2\rho(1)} + 2 \int_0^{\rho(1)} \frac{\d s}{(1 + 2s)^2} = 1.
\end{equation*}
Since the expression on the right-hand side of \eqref{e.multiplier} is strictly decreasing in $b$, its solution satisfies
\begin{equation}
\label{e.spherical-minimizer-bound}
b \le 1 + 2\rho(1).
\end{equation}
The formula for $p_\rho$ in \eqref{e.spherical-derivative} is \cite[(2.38)]{chen2026spherical}. Evaluating it at $u = 1$ and using \eqref{e.multiplier} and \eqref{e.spherical-minimizer-bound}, we obtain
\begin{equation*}
p_\rho(1) = 1 - \frac 1 b \le \frac{2\rho(1)}{1 + 2\rho(1)},
\end{equation*}
which is the second relation in \eqref{e.spherical-derivative}. The derivative formula and the convexity of $\psi^\Sp$ imply \eqref{e.spherical-subgradient}.

By \eqref{e.spherical-subgradient} and \eqref{e.def.psiS-conj}, we have
\begin{equation}
\label{e.spherical-fenchel}
\psi^{\Sp,*}(p_\rho) = \int_0^1 \rho(u) p_\rho(u) \, \d u - \psi^\Sp(\rho).
\end{equation}
For every $p \in \cQ_\infty$, the definition \eqref{e.def.psiS-conj} also gives $\psi^\Sp(\rho) \ge \int_0^1 \rho p - \psi^{\Sp,*}(p)$. Since $0 \le p_\rho \le 1$ by \eqref{e.spherical-derivative}, equality holds at $p = p_\rho$ by \eqref{e.spherical-fenchel}. This proves \eqref{e.spherical-biconjugacy}.
\end{proof}

We recall the definition of $\cE$ from \eqref{e.def.entropy}. For $\boldsymbol{\zeta} \in \cP([0,a])$ with $a < 1$, we have $D_{\boldsymbol{\zeta}}(x) = 1 - x$ for $x \in [a,1)$, and therefore
\begin{equation}
\label{e.entropy-compact}
\cE(\boldsymbol{\zeta}) = \frac 1 2 \Ll( \int_0^a \frac{\d x}{D_{\boldsymbol{\zeta}}(x)} + \log(1-a) \Rr).
\end{equation}
In particular, $\cE(\boldsymbol{\delta}_0) = 0$. The expression \eqref{e.entropy-compact} agrees with the term independent of the covariance in the spherical Crisanti--Sommers functional; see \cite[(2.2)]{subag2026geometry}. The next lemma identifies the convex conjugate of $\psi^\Sp$ with $\cE$ evaluated at the law of the path.

\begin{lemma}[The spherical conjugate]
\label{l.spherical-duality}
For every $a \in [0,1)$ and every $p \in \cQ_\infty$ with $0 \le p \le a$, we have
\begin{equation}
\label{e.spherical-entropy-duality}
\psi^{\Sp,*}(p) = \cE(\boldsymbol{\nu}_p).
\end{equation}
For every $a \in [0,1)$, the function $\cE$ is continuous on $\cP([0,a])$ for the topology of weak convergence. For distinct $\boldsymbol{\zeta}, \boldsymbol{\zeta}' \in \cP([0,a])$ and every $\theta \in (0,1)$, we have
\begin{equation}
\label{e.E-strict-convexity}
\cE(\theta \boldsymbol{\zeta} + (1-\theta) \boldsymbol{\zeta}') < \theta \cE(\boldsymbol{\zeta}) + (1-\theta) \cE(\boldsymbol{\zeta}').
\end{equation}
Moreover, $\cE \ge 0$, with equality only at $\boldsymbol{\delta}_0$.
\end{lemma}
\begin{proof}
\emph{Step 1: construction of $\rho$ such that $p_\rho = p$.}
We fix $a \in [0,1)$ and $p \in \cQ_\infty$ with $0 \le p \le a$. In order to evaluate the supremum in \eqref{e.def.psiS-conj} at $p$, it suffices to construct $\rho \in \cQ_\infty$ for which $p_\rho = p$, since by \eqref{e.spherical-fenchel} such a $\rho$ attains the supremum. This construction occupies the rest of this step. We may replace $a$ by $p(1)$ and we recall $D_{\boldsymbol{\nu}_p}$ from \eqref{e.def.entropy}. 

If $a = 0$, then $p = 0$ and $\boldsymbol{\nu}_p = \boldsymbol{\delta}_0$. At $\rho = 0$, the minimizer in \eqref{e.def.psiS} is $b = 1$ by \eqref{e.multiplier}, so $p_0 = 0$ by \eqref{e.spherical-derivative}, and \eqref{e.spherical-fenchel} gives $\psi^{\Sp,*}(p_0) = -\psi^\Sp(0) = 0 = \cE(\boldsymbol{\delta}_0)$. This proves \eqref{e.spherical-entropy-duality} in this case. We assume from now on that $a > 0$. Since $p \le a$, we have $\boldsymbol{\nu}_p([0,s]) = 1$ for $s \in [a,1]$. Thus, for $x \in [0,a]$, we have
\begin{equation}
D_{\boldsymbol{\nu}_p}(x) \stackrel{\eqref{e.def.entropy}}{=} 1-a + \int_x^a \boldsymbol{\nu}_p([0,s]) \, \d s \ge 1-a > 0.
\label{e.D-path-lower-bound}
\end{equation}
In particular, $D_{\boldsymbol{\nu}_p}(a) = 1-a$. Moreover, $D_{\boldsymbol{\nu}_p}$ is $1$-Lipschitz, with
\begin{equation}
\label{e.D-path-derivative}
D_{\boldsymbol{\nu}_p}'(x) = -\boldsymbol{\nu}_p([0,x]) \qquad \text{for almost every } x \in [0,a].
\end{equation}
We define $G : [0,a] \to \R_+$ by
\begin{equation}
\label{e.spherical-G}
G(x) := \frac 1 2 \int_0^x D_{\boldsymbol{\nu}_p}(z)^{-2} \, \d z.
\end{equation}
The continuity of $D_{\boldsymbol{\nu}_p}$ and \eqref{e.D-path-lower-bound} show that $G$ is continuously differentiable and strictly increasing. We set
\begin{equation}
\rho := G \circ p \in \cQ_\infty \qquad \text{and} \qquad b := (1-a)^{-1}.
\label{e.spherical-inverse-path}
\end{equation}
The equality $\rho(1) = G(a)$ follows from $a = p(1)$ and the continuity of $G$.
For $x \in [0,a]$, the monotonicity of $G$ gives
\begin{equation*}
\rho(1) - \max\{\rho(u), G(x)\} = G(a) - G(\max\{p(u),x\}) = \int_x^a G'(z) \1_{\{z \ge p(u)\}} \, \d z.
\end{equation*}
Integrating this identity in $u$ and using Fubini's theorem, we obtain 
\begin{align}
2 K_\rho(G(x)) & \stackrel{\eqref{e.def.K}}{=} 2 \int_x^a \boldsymbol{\nu}_p([0,z]) G'(z) \, \d z \stackrel{\eqref{e.spherical-G}}{=} \int_x^a \frac{\boldsymbol{\nu}_p([0,z])}{D_{\boldsymbol{\nu}_p}(z)^2} \, \d z \notag \\
& \stackrel{\eqref{e.D-path-derivative}}{=} \Ll[\frac{1}{D_{\boldsymbol{\nu}_p}(z)}\Rr]_x^a = \frac{1}{1-a} - \frac{1}{D_{\boldsymbol{\nu}_p}(x)}.
\label{e.inverse-denominator}
\end{align}
Here, the last equality uses $D_{\boldsymbol{\nu}_p}(a) = 1-a$. Since $G(0) = 0$, \eqref{e.inverse-denominator} at $x = 0$ gives $b - 2K_\rho(0) = D_{\boldsymbol{\nu}_p}(0)^{-1} > 0$, so $b$ is admissible in \eqref{e.def.psiS}. To verify \eqref{e.multiplier}, we make the change of variables $s = G(x)$ and obtain
\begin{equation*}
\frac 1 b + 2 \int_0^{\rho(1)} \frac{\d s}{(b - 2 K_\rho(s))^2} \stackrel{\eqref{e.inverse-denominator}}{=} 1 - a + 2 \int_0^a D_{\boldsymbol{\nu}_p}(x)^2 G'(x) \, \d x \stackrel{\eqref{e.spherical-G}}{=} 1 - a + \int_0^a \d x = 1.
\end{equation*}
By \eqref{e.multiplier}, $b$ is the minimizer in \eqref{e.def.psiS} for $\rho$. For $u \in [0,1]$, the change of variables $s = G(x)$ has endpoints $G(0)=0$ and $G(p(u))=\rho(u)$ by \eqref{e.spherical-inverse-path}. It therefore gives
\begin{equation*}
p_\rho(u) \stackrel{\eqref{e.spherical-derivative}}{=} 2 \int_0^{\rho(u)} \frac{\d s}{(b-2K_\rho(s))^2} \stackrel{\eqref{e.inverse-denominator}}{=} 2 \int_0^{p(u)} D_{\boldsymbol{\nu}_p}(x)^2 G'(x) \, \d x \stackrel{\eqref{e.spherical-G}}{=} \int_0^{p(u)} \d x = p(u).
\end{equation*}
Thus $p_\rho = p$.

\emph{Step 2: evaluation of the conjugate.}
We keep $p$ and $a$ from Step~1 and the path $\rho$ and multiplier $b$ of \eqref{e.spherical-inverse-path}. For $x \in [0,a]$, Fubini's theorem gives 
\begin{align*}
D_{\boldsymbol{\nu}_p}(x) & \stackrel{\eqref{e.def.entropy}}{=} \int_x^1 \int_{[0,a]} \1_{\{z \le s\}} \, \d \boldsymbol{\nu}_p(z) \, \d s \\
& \stackrel{\text{Fubini}}{=} \int_{[0,a]} \int_x^1 \1_{\{z \le s\}} \, \d s \, \d \boldsymbol{\nu}_p(z) = \int_{[0,a]} (1-\max\{x,z\}) \, \d \boldsymbol{\nu}_p(z).
\end{align*}
Splitting this integral over $[0,x]$ and $(x,a]$, we obtain
\begin{equation}
D_{\boldsymbol{\nu}_p}(x) = (1-x)\boldsymbol{\nu}_p([0,x]) + \int_{(x,a]} (1-z) \, \d \boldsymbol{\nu}_p(z) = 1 - x \boldsymbol{\nu}_p([0,x]) - \int_{(x,a]} z \, \d \boldsymbol{\nu}_p(z),
\label{e.D-first-moment}
\end{equation}
where we used $\boldsymbol{\nu}_p([0,x])+\boldsymbol{\nu}_p((x,a]) = 1$.

Since $\rho = G \circ p$ by \eqref{e.spherical-inverse-path} and $G(z) = \int_0^z G'(x) \, \d x$ by \eqref{e.spherical-G}, we have 
\begin{equation}
\int_0^1 \rho(u) p(u) \, \d u \stackrel{\eqref{e.def.path-law}}{=} \int_{[0,a]} G(z) z \, \d \boldsymbol{\nu}_p(z) \stackrel{\text{Fubini}}{=} \int_0^a G'(x) \int_{(x,a]} z \, \d \boldsymbol{\nu}_p(z) \, \d x.
\label{e.spherical-path-pairing}
\end{equation}
The function $x \mapsto x/D_{\boldsymbol{\nu}_p}(x)$ is absolutely continuous, since its denominator is Lipschitz and bounded below by $1-a$ by \eqref{e.D-path-lower-bound}.
Its derivative satisfies
\begin{equation}
\frac{\d}{\d x} \frac{x}{D_{\boldsymbol{\nu}_p}(x)} \stackrel{\eqref{e.D-path-derivative}}{=} \frac{1}{D_{\boldsymbol{\nu}_p}(x)} + \frac{x\boldsymbol{\nu}_p([0,x])}{D_{\boldsymbol{\nu}_p}(x)^2} \qquad \text{for almost every } x \in [0,a].
\label{e.D-quotient-derivative}
\end{equation}
Since $\rho(1) = G(a) = \int_0^a G'(x) \, \d x$ by \eqref{e.spherical-inverse-path} and \eqref{e.spherical-G}, we obtain 
\begin{align}
\rho(1) - \int_0^1 \rho(u) p(u) \, \d u
& \stackrel{\eqref{e.spherical-path-pairing}}{=} \int_0^a G'(x) \Ll( 1 - \int_{(x,a]} z \, \d \boldsymbol{\nu}_p(z) \Rr) \d x \notag \\
& \stackrel{\eqref{e.D-first-moment}}{=} \int_0^a G'(x)\bigl(D_{\boldsymbol{\nu}_p}(x)+x\boldsymbol{\nu}_p([0,x])\bigr) \, \d x \notag \\
& \stackrel{\eqref{e.spherical-G}}{=} \frac 1 2 \int_0^a \Ll( \frac{1}{D_{\boldsymbol{\nu}_p}(x)} + \frac{x \boldsymbol{\nu}_p([0,x])}{D_{\boldsymbol{\nu}_p}(x)^2} \Rr) \d x \notag \\
& \stackrel{\eqref{e.D-quotient-derivative}}{=} \frac 1 2 \Ll[\frac{x}{D_{\boldsymbol{\nu}_p}(x)}\Rr]_0^a = \frac{a}{2(1-a)} = \frac{b-1}{2},
\label{e.spherical-pairing}
\end{align}
where we used $D_{\boldsymbol{\nu}_p}(a) = 1-a$ from \eqref{e.D-path-lower-bound} and $b = (1-a)^{-1}$ from \eqref{e.spherical-inverse-path}. The change of variables $s = G(x)$ also gives
\begin{equation*}
\int_0^{\rho(1)} \frac{\d s}{b-2K_\rho(s)} \stackrel{\eqref{e.inverse-denominator}}{=} \int_0^a D_{\boldsymbol{\nu}_p}(x) G'(x) \, \d x \stackrel{\eqref{e.spherical-G}}{=} \frac 1 2 \int_0^a \frac{\d x}{D_{\boldsymbol{\nu}_p}(x)}.
\end{equation*}
Since $b$ is the minimizer in \eqref{e.def.psiS} by Step~1, we obtain
\begin{equation}
\label{e.spherical-initial-value}
\psi^\Sp(\rho) = \rho(1) - \frac 1 2 \Ll( \int_0^a \frac{\d x}{D_{\boldsymbol{\nu}_p}(x)} + b - 1 - \log b \Rr).
\end{equation}
Since $p = p_\rho$, \eqref{e.spherical-fenchel} gives $\psi^{\Sp,*}(p) = \int_0^1 \rho p - \psi^\Sp(\rho)$. The terms $\rho(1)$ and $(b-1)/2$ cancel, and we obtain
\begin{equation*}
\psi^{\Sp,*}(p) \stackrel{\eqref{e.spherical-pairing},\,\eqref{e.spherical-initial-value}}{=} \frac 1 2 \Ll( \int_0^a \frac{\d x}{D_{\boldsymbol{\nu}_p}(x)} - \log b \Rr) \stackrel{\eqref{e.spherical-inverse-path}}{=} \frac 1 2 \Ll( \int_0^a \frac{\d x}{D_{\boldsymbol{\nu}_p}(x)} + \log(1-a) \Rr) \stackrel{\eqref{e.entropy-compact}}{=} \cE(\boldsymbol{\nu}_p).
\end{equation*}

\emph{Step 3: continuity and strict convexity of $\cE$.} We fix $a \in [0,1)$. If $a = 0$, then $\cP([0,a]) = \{\boldsymbol{\delta}_0\}$ and $\cE(\boldsymbol{\delta}_0) = 0$, so all the assertions hold. Suppose now that $a > 0$. For every $\boldsymbol{\zeta} \in \cP([0,a])$, the definition \eqref{e.def.entropy} gives
\begin{equation}
D_{\boldsymbol{\zeta}}(x) = 1-a + \int_x^a \boldsymbol{\zeta}([0,s]) \, \d s \ge 1-a \quad (x \in [0,a]) \quad \text{and} \quad D_{\boldsymbol{\zeta}}'(x) = -\boldsymbol{\zeta}([0,x]) \quad (\text{a.e. }x \in [0,a]).
\label{e.D-law-properties}
\end{equation}
For $\boldsymbol{\zeta},\boldsymbol{\zeta}' \in \cP([0,a])$, this implies
\begin{equation*}
\sup_{x \in [0,a]} |D_{\boldsymbol{\zeta}}(x)-D_{\boldsymbol{\zeta}'}(x)| \le \int_0^a |\boldsymbol{\zeta}([0,s])-\boldsymbol{\zeta}'([0,s])| \, \d s = W_1(\boldsymbol{\zeta},\boldsymbol{\zeta}').
\end{equation*}
Consequently, we have
\begin{equation*}
|\cE(\boldsymbol{\zeta})-\cE(\boldsymbol{\zeta}')| \stackrel{\eqref{e.entropy-compact}}{\le} \frac 1 2 \int_0^a \frac{|D_{\boldsymbol{\zeta}}(x)-D_{\boldsymbol{\zeta}'}(x)|}{D_{\boldsymbol{\zeta}}(x)D_{\boldsymbol{\zeta}'}(x)} \, \d x \stackrel{\eqref{e.D-law-properties}}{\le} \frac{a}{2(1-a)^2} \, W_1(\boldsymbol{\zeta},\boldsymbol{\zeta}').
\end{equation*}
Since weak convergence on $\cP([0,a])$ is equivalent to convergence for $W_1$, this proves the continuity of $\cE$.

We next prove that $\boldsymbol{\zeta} \mapsto D_{\boldsymbol{\zeta}}$ is injective on $\cP([0,a])$. If $D_{\boldsymbol{\zeta}} = D_{\boldsymbol{\zeta}'}$ on $[0,a]$, their almost-everywhere derivatives agree, and \eqref{e.D-law-properties} gives $\boldsymbol{\zeta}([0,x]) = \boldsymbol{\zeta}'([0,x])$ for almost every $x \in [0,a]$, implying $\boldsymbol{\zeta} = \boldsymbol{\zeta}'$.

Now let us fix distinct $\boldsymbol{\zeta},\boldsymbol{\zeta}' \in \cP([0,a])$ and $\theta \in (0,1)$. By the injectivity just proved and the continuity of $D_{\boldsymbol{\zeta}}$ and $D_{\boldsymbol{\zeta}'}$, these functions differ on a set of positive Lebesgue measure. The definition \eqref{e.def.entropy} is affine in the law, so for every $x \in [0,a]$,
\begin{equation*}
\frac{1}{D_{\theta\boldsymbol{\zeta}+(1-\theta)\boldsymbol{\zeta}'}(x)} = \frac{1}{\theta D_{\boldsymbol{\zeta}}(x)+(1-\theta)D_{\boldsymbol{\zeta}'}(x)} \le \frac{\theta}{D_{\boldsymbol{\zeta}}(x)} + \frac{1-\theta}{D_{\boldsymbol{\zeta}'}(x)}.
\end{equation*}
The inequality is strict whenever $D_{\boldsymbol{\zeta}}(x) \ne D_{\boldsymbol{\zeta}'}(x)$, by the strict convexity of $r \mapsto r^{-1}$ on $(0,\infty)$. Integrating over $[0,a]$, multiplying by $1/2$ and adding $\log(1-a)/2$, we obtain \eqref{e.E-strict-convexity} by \eqref{e.entropy-compact}.

Finally, $D_{\boldsymbol{\zeta}}(x) \le 1-x$ by \eqref{e.def.entropy}, and hence
\begin{equation*}
\cE(\boldsymbol{\zeta}) \stackrel{\eqref{e.entropy-compact}}{=} \frac 1 2 \int_0^a \Ll(\frac{1}{D_{\boldsymbol{\zeta}}(x)}-\frac{1}{1-x}\Rr) \, \d x \ge 0.
\end{equation*}
If $\cE(\boldsymbol{\zeta}) = 0$, the nonnegative continuous integrand vanishes throughout $[0,a]$. Thus $D_{\boldsymbol{\zeta}}(x) = 1-x = D_{\boldsymbol{\delta}_0}(x)$ on $[0,a]$, and injectivity gives $\boldsymbol{\zeta} = \boldsymbol{\delta}_0$. The converse follows from \eqref{e.def.entropy}.
\end{proof}

\subsection{Properties  of $\msc P_a$}
\label{s.field}
Recall that $\phi_{a,\boldsymbol{\zeta}}:[0,1]\times\R\to\R$ solves the Parisi equation \eqref{e.scalar-PDE}, and that  by definition $\msc P_a(\boldsymbol{\zeta},h) = \phi_{a,\boldsymbol{\zeta}}(0,h)$, see \eqref{e.def.Parisi-radius}. In this subsection, we study the functional $\msc P_a$. 

\begin{lemma}
\label{l.squared-field}
For every $a \ge 0$ and $\boldsymbol{\zeta} \in \cP_{<1}$, the function $h \mapsto \msc P_a(\boldsymbol{\zeta},h)$ is even and $1$-Lipschitz, the function $z \mapsto \msc P_a(\boldsymbol{\zeta},\sqrt z)$ is concave on $[0,\infty)$, and for every $h \in \R$,
\begin{equation}
\label{e.field-increment}
0 \le \msc P_a(\boldsymbol{\zeta},h) - \msc P_a(\boldsymbol{\zeta},0) \le \log \cosh h.
\end{equation}
Consequently, for every square-integrable real random variable $H$,
\begin{equation}
\label{e.Jensen-field}
\E \, \msc P_a(\boldsymbol{\zeta},H) \le \msc P_a \bigl( \boldsymbol{\zeta}, \sqrt{\E H^2} \bigr).
\end{equation}
\end{lemma}
\begin{proof}
We fix $a \ge 0$ and $\boldsymbol{\zeta} \in \cP_{<1}$. If $a = 0$, then $\msc P_0(\boldsymbol{\zeta},h) = \log\cosh h$, so evenness, the Lipschitz bound, and \eqref{e.field-increment} follow directly. For $x \ge 0$, we have $\tanh x = \int_0^x \cosh^{-2} y \, \d y \ge x\cosh^{-2} x$, and hence
\begin{equation*}
\frac{\d^2}{\d z^2}\log\cosh\sqrt z = \frac{\sqrt z\,\cosh^{-2}\sqrt z-\tanh\sqrt z}{4z^{3/2}} \le 0 \qquad (z > 0).
\end{equation*}
Continuity at zero extends concavity to $[0,\infty)$, and \eqref{e.Jensen-field} follows as at the end of the proof. We assume from now on that $a > 0$.

\emph{Step 1: evenness, the Lipschitz bound, and \eqref{e.field-increment}.} 
In the notation of \cite[Section~6.1]{boursier2026free}, with terminal time $a$, the change of time $v=as$ gives $\phi_{a,\boldsymbol{\zeta}}(s,x)=\Phi_{\boldsymbol{\nu}_{aQ_{\boldsymbol{\zeta}}}}(as,x)$, since $\boldsymbol{\nu}_{aQ_{\boldsymbol{\zeta}}}([0,as])=\boldsymbol{\zeta}([0,s])$ by \eqref{e.def.path-law}. Thus, by \cite[Lemma~6.1]{boursier2026free}, the spatial derivatives of $\phi_{a,\boldsymbol{\zeta}}$ of every positive order are continuous and bounded on $[0,1]\times\R$, and
\begin{equation}
|\dr_x\phi_{a,\boldsymbol{\zeta}}| < 1 \qquad \text{and} \qquad 0 \le \dr_x^2\phi_{a,\boldsymbol{\zeta}} \le 1-(\dr_x\phi_{a,\boldsymbol{\zeta}})^2.
\label{e.scalar-derivative-bounds}
\end{equation}
The function $(s,x) \mapsto \phi_{a,\boldsymbol{\zeta}}(s,-x)$ solves \eqref{e.scalar-PDE} with the same terminal value. Uniqueness therefore gives evenness in $x$. Together with \eqref{e.def.Parisi-radius} and \eqref{e.scalar-derivative-bounds}, this proves that $\msc P_a(\boldsymbol{\zeta},\cdot)$ is even and $1$-Lipschitz.

We now fix $s \in [0,1]$ and define $u : [0,\infty) \to \R$ by $u(x) := \dr_x\phi_{a,\boldsymbol{\zeta}}(s,x)$. Evenness gives $u(0)=0$, and \eqref{e.scalar-derivative-bounds} gives $0\le u'\le1-u^2$. Since $u<1$ by \eqref{e.scalar-derivative-bounds}, we may integrate to obtain
\begin{equation*}
\operatorname{arctanh}u(x) = \int_0^x\frac{u'(y)}{1-u(y)^2}\,\d y \le \int_0^x\d y = x \qquad (x \ge 0).
\end{equation*}
Consequently, we have
\begin{equation}
0 \le \dr_x\phi_{a,\boldsymbol{\zeta}}(s,x) \le \tanh x \qquad ((s,x) \in [0,1]\times[0,\infty)).
\label{e.field-gradient-bound}
\end{equation}
Hence, the function $x\mapsto\phi_{a,\boldsymbol{\zeta}}(0,x)$ is nondecreasing on $[0,\infty)$. Its evenness therefore gives $\phi_{a,\boldsymbol{\zeta}}(0,h)\ge\phi_{a,\boldsymbol{\zeta}}(0,0)$ for every $h\in\R$. Using this and integrating the above display, we obtain
\begin{equation*}
0 \le \msc P_a(\boldsymbol{\zeta},h)-\msc P_a(\boldsymbol{\zeta},0) \stackrel{\eqref{e.def.Parisi-radius}}{=} \int_0^{|h|}\dr_x\phi_{a,\boldsymbol{\zeta}}(0,x)\,\d x \stackrel{\eqref{e.field-gradient-bound}}{\le} \int_0^{|h|}\tanh x\,\d x = \log\cosh h.
\end{equation*}
This proves \eqref{e.field-increment}.

\emph{Step 2: concavity for finitely supported laws.} Suppose that $\boldsymbol{\zeta}$ is finitely supported. We define $f : [0,a]\times\R \to \R$ and $\vartheta : [0,a] \to [0,1]$ by
\begin{equation}
f(\tau,x) := \phi_{a,\boldsymbol{\zeta}}(1-\tau/a,x) \qquad \text{and} \qquad \vartheta(\tau) := \boldsymbol{\zeta}([0,1-\tau/a]).
\label{e.field-time-reversal}
\end{equation}
The coefficient $\vartheta$ is a step function. On each open interval on which it is constant, the chain rule and \eqref{e.scalar-PDE} give
\begin{equation}
\dr_\tau f = -\frac1a\dr_s\phi_{a,\boldsymbol{\zeta}}(1-\tau/a,x) = \dr_x^2 f+\vartheta (\dr_x f)^2 \qquad \text{and} \qquad f(0,x) = \log\cosh x.
\label{e.forward-scalar}
\end{equation}
Moreover, $f(a,h)=\msc P_a(\boldsymbol{\zeta},h)$ by \eqref{e.def.Parisi-radius}. We will prove that
\begin{equation}
\dr_x^3 f(\tau,x) \le 0 \qquad ((\tau,x) \in [0,a]\times[0,\infty)).
\label{e.field-third-derivative-sign}
\end{equation}
For this purpose, let $w := \dr_x^3 f$. Since
\begin{equation*}
\dr_x^3\bigl((\dr_x f)^2\bigr) = \dr_x\bigl(2(\dr_x^2 f)^2+2\,\dr_x f\,\dr_x^3 f\bigr) = 6\,\dr_x^2 f\,\dr_x^3 f+2\,\dr_x f\,\dr_x^4 f = 6\,\dr_x^2 f\,w+2\,\dr_x f\,\dr_x w,
\end{equation*}
differentiating \eqref{e.forward-scalar} three times in $x$ yields
\begin{equation}
\dr_\tau w = \dr_x^2 w+2\vartheta\,\dr_x f\,\dr_x w+6\vartheta\,\dr_x^2 f\,w.
\label{e.third-derivative-PDE}
\end{equation}
The initial condition in \eqref{e.forward-scalar} and evenness in $x$ give
\begin{equation*}
w(0,x) = -2\cosh^{-2}x\tanh x \le 0 \quad (x \ge 0) \qquad \text{and} \qquad w(\tau,0)=0 \quad (\tau \in [0,a]).
\end{equation*}
All spatial derivatives of $f$ are continuous and bounded, as recalled in Step~1. In particular, $w$ is bounded, and its values at the endpoints of the intervals on which $\vartheta$ is constant agree. Also, $|\dr_x f|\le1$ and $0\le \dr_x^2 f\le1$ by \eqref{e.scalar-derivative-bounds} and \eqref{e.field-time-reversal}.

We apply the maximum principle on bounded spatial intervals and then remove their right endpoints. We define $v : [0,a]\times[0,\infty) \to \R$ by $v(\tau,x):=e^{-6\tau}w(\tau,x)$. Equation \eqref{e.third-derivative-PDE} gives
\begin{equation}
\dr_\tau v-\dr_x^2 v-2\vartheta\,\dr_x f\,\dr_x v+6(1-\vartheta\,\dr_x^2 f)v=0.
\label{e.field-damped-PDE}
\end{equation}
Here $|2\vartheta\,\dr_x f|\le2$ and $6(1-\vartheta\,\dr_x^2 f)\ge0$. The function $e^{5\tau}(1+x^2)$ is a strict supersolution of \eqref{e.field-damped-PDE} on $x\ge0$, since
\begin{align*}
&\bigl(\dr_\tau-\dr_x^2-2\vartheta (\dr_x f)\dr_x+6(1-\vartheta\,\dr_x^2 f)\bigr)\bigl(e^{5\tau}(1+x^2)\bigr) \\
&\qquad = e^{5\tau}\bigl(5(1+x^2)-2-4\vartheta x\,\dr_x f+6(1-\vartheta\,\dr_x^2 f)(1+x^2)\bigr) \\
&\qquad \ge e^{5\tau}(5x^2-4x+3) = e^{5\tau}\bigl(5(x-2/5)^2+11/5\bigr)>0.
\end{align*}
We fix $\ep>0$ and choose $R>0$ such that $\ep(1+R^2)\ge\sup_{[0,a]\times[0,\infty)}|v|$. Then $v-\ep e^{5\tau}(1+x^2)$ is nonpositive on $\{0\} \times [0,R]$ and on $[0,a]\times \{0,R\}$. Applying the maximum principle to \eqref{e.field-damped-PDE} successively on the finitely many intervals on which $\vartheta$ is constant, and using continuity at their endpoints, we obtain
\begin{equation*}
v(\tau,x) \le \ep e^{5\tau}(1+x^2) \qquad ((\tau,x) \in [0,a]\times[0,R]).
\end{equation*}
For every fixed $(\tau,x)\in[0,a]\times[0,\infty)$ and every $\ep>0$, we may choose such an $R>x$. Letting $\ep\downarrow0$ gives $v(\tau,x)\le0$, hence $w(\tau,x)\le0$. This proves \eqref{e.field-third-derivative-sign}.

Since $\dr_x f(\tau,0)=0$ by evenness, we have, for $x\ge0$,
\begin{equation*}
\dr_x f(\tau,x) = \int_0^x \dr_x^2 f(\tau,y)\,\d y \stackrel{\eqref{e.field-third-derivative-sign}}{\ge} \int_0^x \dr_x^2 f(\tau,x)\,\d y = x\,\dr_x^2 f(\tau,x).
\end{equation*}
The change of variables $x=\sqrt z$ therefore gives
\begin{equation*}
\frac{\d^2}{\d z^2}f(\tau,\sqrt z) = \frac{\d}{\d z}\frac{\dr_x f(\tau,\sqrt z)}{2\sqrt z} = \frac{\sqrt z\,\dr_x^2 f(\tau,\sqrt z)-\dr_x f(\tau,\sqrt z)}{4z^{3/2}} \le 0 \qquad (z>0).
\end{equation*}
Continuity at zero gives concavity on $[0,\infty)$. Taking $\tau=a$ proves the assertion for finitely supported $\boldsymbol{\zeta}$.

\emph{Step 3: general laws and \eqref{e.Jensen-field}.} For arbitrary $\boldsymbol{\zeta},\boldsymbol{\zeta}'\in\cP_{<1}$ and $h\in\R$, we have
\begin{equation}
|\msc P_a(\boldsymbol{\zeta},h)-\msc P_a(\boldsymbol{\zeta}',h)| \stackrel{\eqref{e.field-duality}}{=} |\psi^\Is(aQ_{\boldsymbol{\zeta}};h)-\psi^\Is(aQ_{\boldsymbol{\zeta}'};h)| \stackrel{\mathrm{L.\,\ref{l.Ising-data}}}{\le} a|Q_{\boldsymbol{\zeta}}-Q_{\boldsymbol{\zeta}'}|_{L^1} = aW_1(\boldsymbol{\zeta},\boldsymbol{\zeta}').
\label{e.P-law-continuity}
\end{equation}
We choose $b<1$ such that $\boldsymbol{\zeta}\in\cP([0,b])$. For $n\in\N$, let $\boldsymbol{\zeta}_n$ be the image of $\boldsymbol{\zeta}$ under $x\mapsto\lfloor nx\rfloor/n$. This law is finitely supported in $[0,b]$, and the coupling given by this map yields
\begin{equation*}
W_1(\boldsymbol{\zeta},\boldsymbol{\zeta}_n) \le \int_{[0,b]}\left|x-\frac{\lfloor nx\rfloor}{n}\right|\,\d\boldsymbol{\zeta}(x) \le \frac1n.
\end{equation*}
Consequently, we have
\begin{equation}
\sup_{h\in\R}|\msc P_a(\boldsymbol{\zeta}_n,h)-\msc P_a(\boldsymbol{\zeta},h)| \stackrel{\eqref{e.P-law-continuity}}{\le} \frac an.
\label{e.field-law-approximation}
\end{equation}
By \eqref{e.field-law-approximation}, the map $z\mapsto\msc P_a(\boldsymbol{\zeta}_n,\sqrt z)$ converges uniformly on $[0,\infty)$, as $n\to\infty$, to $z\mapsto\msc P_a(\boldsymbol{\zeta},\sqrt z)$. Since each function in the sequence is concave by Step~2, the limit is concave as well. 

Finally, for every $a\ge0$ and every real random variable $H$ with $\E H^2<\infty$, the Lipschitz bound gives
\begin{equation*}
\E|\msc P_a(\boldsymbol{\zeta},H)| \le |\msc P_a(\boldsymbol{\zeta},0)|+\E|H| \le |\msc P_a(\boldsymbol{\zeta},0)|+\sqrt{\E H^2}<\infty.
\end{equation*}
Evenness and the concavity just proved therefore allow us to apply Jensen's inequality to $H^2$:
\begin{equation*}
\E\,\msc P_a(\boldsymbol{\zeta},H) = \E\,\msc P_a(\boldsymbol{\zeta},\sqrt{H^2}) \stackrel{\text{Jensen}}{\le} \msc P_a\bigl(\boldsymbol{\zeta},\sqrt{\E H^2}\bigr).
\end{equation*}
This proves \eqref{e.Jensen-field}.
\end{proof}

%
%
%
%
%
%
\section{Multi-species model on a slice}
\label{s.mixed}

In this section, we identify the limit free energy of several Ising species with prescribed magnetizations coupled to one spherical species. Theorem~\ref{t.mixed-slice} gives the limit as a viscosity solution and a Hopf formula. For Ising slices, this is \cite[Theorem~1.2]{boursier2026free}; related constrained Hamilton--Jacobi arguments appear in \cite[Section~4]{ho2026convergence} and \cite[Propositions~3.1 and~5.6]{chen2026formula}.

\subsection{A multi-species spin glass model}
\label{s.mixed-setting}

Let $J$ be a finite set, and let $\cS := J \cup \{\Sp\}$ be the set of species, where $\Sp$ labels the spherical species. For each $n \in \N$, we partition the coordinates $\{1,\ldots,n\}$ into species of sizes $(n_s)_{s \in \cS}$, with $\sum_s n_s = n$. We fix $\lambda_s\in(0,1)$ with $\sum_{s\in\cS}\lambda_s=1$, and assume
\begin{equation}
\label{e.proportions}
\lambda_{n,s} := \frac{n_s}{n} \xrightarrow[n \to \infty]{} \lambda_s > 0, \qquad \text{for every } s \in \cS.
\end{equation}
For $a=(a_j)_{j\in J}$ with $a_j\in\{0,\ldots,n_j\}$, set $m_{n,a,j}:=2a_j/n_j-1$ and $m_{n,a}:=(m_{n,a,j})_{j\in J}$. We define
\begin{equation}
\label{e.def.slice-measure}
\Omega_{n_j,a_j} := \bigg\{ \si \in \{-1,1\}^{n_j} : \sum_{i=1}^{n_j} \1_{\{\si_i = 1\}} = a_j \bigg\}\qquad\text{and}\qquad\boldsymbol{\mu}_{n,a} := \bigotimes_{j \in J} \bigg( 2^{-n_j} \sum_{\si \in \Omega_{n_j,a_j}} \boldsymbol{\delta}_{\si} \bigg) \otimes \boldsymbol{\mu}_{n_\Sp},
\end{equation}
where $\boldsymbol{\mu}_{n_\Sp}$ is the uniform probability measure on $\Sph_{n_\Sp}$. We call the support of $\boldsymbol{\mu}_{n,a}$ a \textit{slice}, which is the product of the Ising sets $\Omega_{n_j,a_j}$ and the sphere $\Sph_{n_\Sp}$. 

Let $P_s : \R^n \to \R^n$ be the orthogonal projection onto the coordinates of species $s$. For $\sigma, \sigma' \in \R^n$, we define the species overlaps $R(\sigma,\sigma') \in \R^{\cS}$ by $R_s(\sigma,\sigma') := n^{-1} P_s \sigma \cdot P_s \sigma'$. On $\supp\boldsymbol{\mu}_{n,a}$, we have
\begin{equation}
\label{e.fixed-norms}
|P_s \sigma|^2 = n_s, \qquad |\sigma|^2 = n, \qquad \text{and} \qquad R(\sigma,\sigma) = \lambda_n := (\lambda_{n,s})_{s \in \cS}.
\end{equation}
Let $\xi : \R^{\cS} \to \R$ be a polynomial with nonnegative coefficients and $\xi(0)=0$. We write
\begin{equation}
\xi(r)=\sum_{\mathbf p\in(\N\cup\{0\})^{\cS}}c_{\mathbf p}\prod_{s\in\cS}r_s^{p_s}, \qquad c_{\mathbf p}\ge0, \qquad\text{and}\qquad |\mathbf p|:=\sum_{s\in\cS}p_s,
\label{e.covariance-expansion}
\end{equation}
where only finitely many coefficients $c_{\mathbf p}$ are nonzero. Let $H_n:\R^n\to\R$ be a centered Gaussian field with covariance
\begin{equation}
\E[H_n(\sigma)H_n(\sigma')]=n\xi(R(\sigma,\sigma')), \qquad \sigma,\sigma'\in\R^n.
\label{e.mixed-covariance}
\end{equation}

Let $\mathfrak R$ be a Ruelle probability cascade, independent of $H_n$, whose overlap is uniformly distributed on $[0,1]$; this is a random probability measure on the unit sphere of a separable Hilbert space, and we write $\ga \wedge \ga'$ for the scalar product of two points $\ga, \ga'$ in its support; see \cite[Section~4 and Proposition~4.1]{chen2025vector} for this construction and its Gaussian fields. Conditionally on $\mathfrak R$, for $q = (q_s)_{s \in \cS} \in \cQ_\infty^{\cS}$, let $(w_i^q(\ga))_{\ga \in \supp \mathfrak R}$, for $i \le n$, be independent centered Gaussian processes, independent of $H_n$, with conditional covariance $\E[w_i^q(\ga) w_i^q(\ga')\mid\mathfrak R] = q_s(\ga \wedge \ga')$ when $i$ belongs to species $s$. We set
\begin{align}
\label{e.def.enriched}
\bar F_{n,a}(t,q) := - \frac 1 n \E \log \iint \exp \Bigl( & \sqrt{2t} H_n(\sigma) + \sqrt 2 \sum_{i=1}^n w_i^q(\ga) \sigma_i \notag \\
& - n t \xi(\lambda_n) - n q(1) \cdot \lambda_n \Bigr) \, \d \boldsymbol{\mu}_{n,a}(\sigma) \, \d \mathfrak R(\ga).
\end{align}
Here, $\E$ averages over $H_n$, $\mathfrak R$, and the fields $w_i^q$. The quantities $nt\xi(\lambda_n)$ and $nq(1)\cdot\lambda_n$ are one half of the variances of the Gaussian fields (and do not depend on $\sigma$ here). By \cite[Proposition~2.1]{chen2026spherical}, which is the multi-species form of \cite[Proposition~5.1]{chen2025vector}, the function $\bar F_{n,a}$ extends uniquely to $[0,\infty)\times\cQ_1^{\cS}$ and satisfies, for every $t,t'\ge0$ and $q,q'\in\cQ_1^{\cS}$,
\begin{equation}
\label{e.mixed-Lipschitz}
|\bar F_{n,a}(t,q) - \bar F_{n,a}(t',q')| \le C_{0,\xi} |t - t'| + \sum_s \lambda_{n,s} |q_s - q_s'|_{L^1},
\end{equation}
where $C_{0,\xi}:=\sup_{[-1,1]^{\cS}}|\xi|$.
The same results give differentiability on $[0,\infty)\times\cQ_2^{\cS}$. The path derivative (see \eqref{e.free-energy-derivatives}) belongs to $\cQ_\infty^{\cS}$ and satisfies $0\le\dr_q\bar F_{n,a}(t,q)\le\lambda_n$.

\subsection{Hamilton--Jacobi remainder estimate}
\label{s.branch}

A path $q : [0,1] \to \R^{\cS}$ is called \textbf{regular} if it is continuous and piecewise affine, with $q_s(0) > 0$ and all one-sided slopes positive, for every $s \in \cS$. Our goal is to prove the following.
\begin{proposition}
\label{p.transfer}
For every $\kappa>0$, $t\ge0$, $n\geq \kappa^{-1}$, and regular path $q$ satisfying $\dot q_s\ge\kappa$ almost everywhere for every $s\in\cS$, we have
\begin{equation}
\label{e.mixed-residual}
\Ll| \dr_t \bar F_{n,a}(t,q) - \int_0^1 \xi \bigl( \dr_q \bar F_{n,a}(t,q)(u) \bigr) \, \d u \Rr| \le 15 \, C_{2,\xi} (\kappa n)^{-1/3},
\end{equation}
where $C_{2,\xi}:=\sup_{[-1,1]^{\cS}}\norm{\nabla^2\xi}_{\mathrm{op}}$. The bound is uniform over $a$. 
\end{proposition}

We first record the Gaussian representation, the symmetries of the slice, and the associated Parisi PDE.

\smallskip
\noindent\textbf{Representation of $H_n(\sigma)$.}
For $\sigma\in\R^n$, we write $\sigma_s\in\R^{n_s}$ for its species-$s$ block. In view of~\eqref{e.covariance-expansion}, we define $d\in\N\cup\{0\}$ and $A:\R^n\to\R^d$ by
\begin{equation}
d:=\sum_{\mathbf p:\,c_{\mathbf p}>0}\prod_{s\in\cS}n_s^{p_s} \qquad\text{and}\qquad A(\sigma):=\biggl(\sqrt{c_{\mathbf p}n^{1-|\mathbf p|}}\,\bigotimes_{s\in\cS}\sigma_s^{\otimes p_s}\biggr)_{\mathbf p:\,c_{\mathbf p}>0}.
\label{e.Gram-map}
\end{equation}
Taking tensor inner products in \eqref{e.Gram-map} and using \eqref{e.covariance-expansion}, we obtain
\begin{equation}
\label{e.Gram}
A(\sigma) \cdot A(\sigma') = n \xi(R(\sigma,\sigma'))
\qquad \text{and} \qquad
|A(\sigma)|^2 = n \xi(\lambda_n) \quad \text{on the support of } \boldsymbol{\mu}_{n,a}.
\end{equation}
Let $g$ be an independent standard Gaussian vector in $\R^d$. By \eqref{e.Gram} and \eqref{e.mixed-covariance}, the fields $H_n$ and $(g\cdot A(\sigma))_{\sigma\in\R^n}$ have the same law. We henceforth use this realization.

\smallskip
\noindent\textbf{Symmetry on the slice.}
The slice $\supp\boldsymbol{\mu}_{n,a}$ is preserved by the group 
\begin{equation}
\label{e.orbit-group}
\mathsf{Sym}:=\prod_{j \in J} \mfk S_{n_j} \times O(n_\Sp),
\end{equation}
acting blockwise by coordinate permutations in $\mfk S_{n_j}$ on the Ising species and orthogonal transformations in $O(n_\Sp)$ on the spherical species. For every $S=(S_s)_{s\in\cS}\in\mathsf{Sym}$, we have
\begin{equation}
SP_s=P_sS \qquad\text{and}\qquad \boldsymbol{\mu}_{n,a}\circ S^{-1}=\boldsymbol{\mu}_{n,a}.
\label{e.slice-symmetries}
\end{equation}

To describe the corresponding change of the Gaussian vector, we write $z=(z_{\mathbf p})_{\mathbf p:\,c_{\mathbf p}>0}\in\R^d$, using the same blocks as in \eqref{e.Gram-map}, and we define $O_S:\R^d\to\R^d$ by
\begin{equation}
O_Sz:=\biggl(\Bigl(\bigotimes_{s\in\cS}S_s^{\otimes p_s}\Bigr)z_{\mathbf p}\biggr)_{\mathbf p:\,c_{\mathbf p}>0}.
\label{e.disorder-transformation}
\end{equation}
Since each $S_s$ is orthogonal, \eqref{e.Gram-map} and \eqref{e.disorder-transformation} give
\begin{equation}
O_S^\top O_S=\Id, \qquad A(S\sigma)=O_SA(\sigma), \qquad\text{and}\qquad O_Sg\stackrel{\d}{=}g.
\label{e.disorder-symmetry}
\end{equation}

Finally, the action of $\mathsf{Sym}$ on $\supp\boldsymbol{\mu}_{n,a}$ is \textbf{transitive}, meaning that
\begin{equation}
\text{for every }\sigma,\sigma'\in\supp\boldsymbol{\mu}_{n,a},\quad \text{there exists }S\in\mathsf{Sym}\text{ such that }S\sigma=\sigma'.
\label{e.slice-transitivity}
\end{equation}
Indeed, coordinate permutations match any two Ising blocks with the same number $a_j$ of $+1$ entries, and an orthogonal matrix matches any two spherical vectors of norm $\sqrt{n_\Sp}$. 

\smallskip
\noindent\textbf{Parisi PDE representation.}
The cascade integral in \eqref{e.def.enriched} can be represented by the following Parisi PDE, as in \cite[Proposition~A.8]{boursier2026free}. For a fixed $t\ge0$ and a regular path $q$, we set $\sfD(u):=\sum_s\dot q_s(u)P_s$, and let $\Phi:[0,1]\times\R^n\times\R^d\to\R$ solve
\begin{equation}
\label{e.mixed-PDE}
\begin{cases}
\dr_u \Phi + \tr(\sfD(u) \nabla_x^2 \Phi) + u \nabla_x \Phi \cdot \sfD(u) \nabla_x \Phi = 0, \\[2pt]
\Phi(1,x;g) = \log \displaystyle \int e^{\sqrt{2t} g \cdot A(\sigma) + x \cdot \sigma} \, \d \boldsymbol{\mu}_{n,a}(\sigma).
\end{cases}
\end{equation}
The terminal datum has bounded derivatives of every positive order because $\boldsymbol{\mu}_{n,a}$ has compact support, so the proof of \cite[Lemma~2.1]{boursier2026free} gives a unique solution. The equation is solved separately on the intervals on which $\sfD$ is constant, with continuous matching at their endpoints. 

Let $X_0$ be an independent centered Gaussian vector with covariance $2 \sum_s q_s(0) P_s$ and let $B$ be an independent standard Brownian motion in $\R^n$. The diffusion $X$ is the solution of
\begin{equation}
\label{e.branch}
\d X_u = 2 u \sfD(u) \nabla_x \Phi(u,X_u;g) \, \d u + \sqrt{2 \sfD(u)} \, \d B_u.
\end{equation}
The representation \cite[(2.13) and Proposition~A.8]{boursier2026free} gives
\begin{equation}
\bar F_{n,a}(t,q)=t\xi(\lambda_n)+q(1)\cdot\lambda_n-\frac1n\E\Phi(0,X_0;g).
\label{e.free-energy-PDE}
\end{equation}
Here the expectation is over $(g,X_0)$. 

\begin{proof}[Proof of Proposition~\ref{p.transfer}]
The proof follows \cite[Theorem~5.1, Sections~5.1--5.4]{boursier2026free}, which proves \eqref{e.mixed-residual} for multi-species models with Ising spins only. We refer to that paper for the arguments that apply without change, and give the details at the places where the spherical species require a modification. In outline, Step~0 couples the diffusion \eqref{e.branch} with a configuration $\sigma$ on the slice, and uses \cite[Theorem~4.1]{boursier2026free} to express the derivatives $\bar F_{n,a}$ in terms of averages under certain Gibbs measures. The remainder in \eqref{e.mixed-residual} is then controlled by the variances of the overlap under these Gibbs measures. Steps~1 and~2 bound these variances, first conditionally on $\sigma$, by the Brascamp--Lieb inequality, and then unconditionally, using the symmetry of the slice to show that the conditional law does not depend on $\sigma$. Step~3 collects the estimates. The modifications required by our present setting are minimal: the a priori bounds used in Steps~1 and~3 rely on the fixed norms \eqref{e.fixed-norms}, which hold on the slice for both types of species, and the symmetry group used in Step~2 must include the orthogonal action on the spherical species.

\emph{Step 0: the conditional moments and free-energy derivatives.} For the diffusion \eqref{e.branch}, we extend the probability space carrying $(g,X_0,B)$ by a configuration $\sigma\in\supp\boldsymbol{\mu}_{n,a}$ whose conditional distribution is
\begin{equation}
\P(\sigma\in\mathrm d\tau\mid g,X_0,B)=\exp\bigl(\sqrt{2t}g\cdot A(\tau)+X_1\cdot\tau-\Phi(1,X_1;g)\bigr)\,\boldsymbol{\mu}_{n,a}(\mathrm d\tau).
\label{e.terminal-Gibbs}
\end{equation}
Here $\P$ is the joint law of $(g,X_0,B,\sigma)$. The terminal condition in \eqref{e.mixed-PDE} makes \eqref{e.terminal-Gibbs} a probability measure. 
Let $\E$, $\Cov$, $\Var$, and $\Law$ be taken under $\P$. We let $\cG_u$ denote the $\sigma$-algebra generated by $\bigl( g, X_0, (B_w)_{w \le u} \bigr)$, and set
\begin{equation*}
 \mathsf m_u := \E[\sigma \mid \cG_u], \qquad 
\sfC_u := \Cov(\sigma \mid \cG_u), \qquad\text{and}\qquad 
\mathsf H_u := \nabla_x^2 \Phi(u,X_u;g).
\end{equation*}
Thus $\Cov(\sigma\mid\cG_u)=\E[(\sigma-\mathsf m_u)(\sigma-\mathsf m_u)^\top\mid\cG_u]$. The conditional-moment identity \cite[Lemma~3.1]{boursier2026free} gives $\mathsf m_u=\nabla_x\Phi(u,X_u;g)$. Here, $X$ and $\mathsf m$ take values in $\R^n$, while $\sfC$ and $\mathsf H$ take values in $\R^{n\times n}$. We define the deterministic paths $p_{n,s} : [0,1] \to \R_+$ by
\begin{equation}
\label{e.def.pn}
p_{n,s}(u) := \frac 1 n \E |P_s \mathsf m_u|^2, \qquad \text{so that} \qquad \sum_s p_{n,s}(u) = \frac 1 n \E |\mathsf m_u|^2 \le 1.
\end{equation}
Differentiating \eqref{e.free-energy-PDE} as in \cite[Theorem~4.1]{boursier2026free} identifies $\dr_q\bar F_{n,a}(t,q)=p_n:=(p_{n,s})_{s\in\cS}$. Thus, for every continuous piecewise affine $\rho:[0,1]\to\R^{\cS}$, 
\begin{equation}
\lim_{\ep\to0}\frac{\bar F_{n,a}(t,q+\ep\rho)-\bar F_{n,a}(t,q)}{\ep}=\int_0^1p_n(u)\cdot\rho(u)\,\d u=\sum_{s\in\cS}\int_0^1p_{n,s}(u)\rho_s(u)\,\d u.
\label{e.free-energy-derivatives}
\end{equation}
We write $\sfC_{u,ss}:=P_s\sfC_uP_s$ and $\mathsf H_{u,ss}:=P_s\mathsf H_uP_s$ for the diagonal blocks.

\emph{Step 1: the conditional variance estimate.} 
Let $\boldsymbol{\gamma}$ be the joint Gaussian law of $g$ and $Z:=X_0+\int_0^\cdot\sqrt{2\sfD(w)}\,\d B_w$ on $\R^d\times C([0,1];\R^n)$. For $\tau\in\supp\boldsymbol{\mu}_{n,a}$, we define $f_\tau:\R^d\times C([0,1];\R^n)\to(0,\infty)$ by
\begin{equation}
\label{e.conditional-density}
f_\tau(g,x) := \exp \Ll( \sqrt{2t} \, g \cdot A(\tau) + x_1 \cdot \tau - \int_0^1 \Phi(u,x_u;g) \, \d u \Rr).
\end{equation}
The Girsanov calculation of \cite[Lemma~5.2]{boursier2026free}, combined with \eqref{e.terminal-Gibbs}, shows that under $\P$, the joint law of $(g,X,\sigma)$ has density $(g,x,\tau)\mapsto f_\tau(g,x)$ with respect to $\boldsymbol{\gamma}\otimes\boldsymbol{\mu}_{n,a}$. Consequently,
\begin{equation*}
\Law((g,X)\mid \sigma=\tau)=\frac{f_\tau\,\boldsymbol{\gamma}}{\int f_\tau\,\d\boldsymbol{\gamma}}, \qquad \tau\in\supp\boldsymbol{\mu}_{n,a}.
\end{equation*}

For $s\in\cS$, we define, as in \cite[(5.3)]{boursier2026free}, 
\begin{equation}
L_{u,s} := \lambda_{n,s} - \frac 1 n |P_s(\sigma - \mathsf m_u)|^2.
\label{e.def.L}
\end{equation}
For $\vartheta\in L^2([0,1])$, let $V:=\int_0^1\vartheta(u)L_{u,s}\,\d u$. The proof of \cite[Lemma~5.3]{boursier2026free} gives joint convexity of $\Phi$ in $(x,g)$, hence concavity of $\log f_\tau$. The finite-dimensional approximation and Brascamp--Lieb argument of \cite[Lemma~5.6]{boursier2026free} therefore give
\begin{equation}
\label{e.conditional-variance}
\Var(V\mid \sigma=\tau)\le\frac4{n^2}\int_0^1\vartheta(u)^2\E\bigl[P_s(\sigma-\mathsf m_u)\cdot\mathsf H_uP_s(\sigma-\mathsf m_u)\mid \sigma=\tau\bigr]\,\d u.
\end{equation}
The bounds needed in \cite[(5.14)--(5.17)]{boursier2026free} are $|P_s(\tau-\nabla_x\Phi)|\le2\sqrt{n_s}$ and $0\preceq\nabla_x^2\Phi\preceq n\Id$. They follow from \eqref{e.fixed-norms} and the conditional-moment identities for the diffusion started at any $(u,x)$, as in \cite[Lemma~3.1 and (3.13)]{boursier2026free}.

\emph{Step 2: removal of the conditioning.} For $S\in\mathsf{Sym}$ from \eqref{e.orbit-group}, the transformation $(g,x,\tau)\mapsto(O_Sg,Sx,S\tau)$ preserves $\boldsymbol{\gamma}\otimes\boldsymbol{\mu}_{n,a}$: $O_S$ is orthogonal, $S$ commutes with every $P_s$, and \eqref{e.slice-symmetries} preserves the reference measure. Uniqueness in \eqref{e.mixed-PDE} and \eqref{e.disorder-symmetry} give
\begin{equation}
\Phi(u,Sx;O_Sg)=\Phi(u,x;g) \qquad\text{and}\qquad \nabla_x\Phi(u,Sx;O_Sg)=S\nabla_x\Phi(u,x;g).
\label{e.mixed-branch-symmetry}
\end{equation}
By \eqref{e.mixed-branch-symmetry}, the density \eqref{e.conditional-density} satisfies $f_{S\tau}(O_Sg,Sx)=f_\tau(g,x)$. Since $(g,x)\mapsto(O_Sg,Sx)$ preserves~$\boldsymbol{\gamma}$, the formula for $\Law((g,X)\mid\sigma=\tau)$ in Step~1 shows that $\Law((g,X)\mid\sigma=S\tau)$ is the image of $\Law((g,X)\mid\sigma=\tau)$ under this map. The random variable $V$ is a function of $(g,X,\sigma)$ which, by \eqref{e.def.L} and \eqref{e.mixed-branch-symmetry}, is invariant under $(g,x,\tau)\mapsto(O_Sg,Sx,S\tau)$. Hence $\Law(V\mid\sigma=S\tau)=\Law(V\mid\sigma=\tau)$ for every $\tau$. By the transitivity \eqref{e.slice-transitivity}, this law does not depend on $\tau$, and hence $\Var(V)=\E\Var(V\mid\sigma)$. This replaces \cite[(5.22)--(5.23)]{boursier2026free}. Averaging the conditional estimate \eqref{e.conditional-variance} over $\tau$ and using
\begin{equation*}
\E\bigl[P_s(\sigma-\mathsf m_u)\cdot\mathsf H_uP_s(\sigma-\mathsf m_u)\mid\cG_u\bigr]=\tr(\mathsf H_{u,ss}\sfC_{u,ss}),
\end{equation*}
we obtain
\begin{equation}
\label{e.BL-gain}
\Var \Ll( \int_0^1 \vartheta(u) L_{u,s} \, \d u \Rr) \le \frac{4}{n^2} \int_0^1 \vartheta(u)^2 \, \E \tr(\mathsf H_{u,ss} \sfC_{u,ss}) \, \d u.
\end{equation}

\emph{Step 3: the Hamilton--Jacobi remainder.} For $u\in[0,1]$, let $\langle\cdot\rangle_u$ denote the expectation with respect to two replicas $\sigma^1,\sigma^2$ which, conditionally on $\cG_u$, are independent with the conditional law of $\sigma$, and let $R:=R(\sigma^1,\sigma^2)$, so that $\E\langle R\rangle_u=p_n(u)$ by \eqref{e.def.pn}. The time-derivative identity of \cite[Theorem~4.1]{boursier2026free}, with \eqref{e.Gram} and \eqref{e.free-energy-derivatives}, gives
\begin{equation}
\mathscr E := \dr_t\bar F_{n,a}(t,q)-\int_0^1\xi(p_n(u))\,\d u =\int_0^1\Bigl\{\E\langle\xi(R)\rangle_u-\xi\bigl(\E\langle R\rangle_u\bigr)\Bigr\}\,\d u.
\label{e.mixed-residual-representation}
\end{equation}
A second-order Taylor expansion of $\xi$ around $p_n(u)$, whose linear term has zero mean under $\E\langle\cdot\rangle_u$, and the conditional variance decomposition of $R$ given $\cG_u$, as in \cite[(4.9) and (5.6)]{boursier2026free}, give the first inequality in \eqref{e.mixed-residual-bound} below. The second one follows as in \cite[(5.25) and (5.28)]{boursier2026free} from \eqref{e.BL-gain} and from \cite[(3.9), (3.15) and (3.17)]{boursier2026free}, namely the response bound $\mathsf H_u\succeq u\sfC_u$, the martingale property of $\mathsf m$, and $\int_0^1\E\tr(\mathsf H_u^2)\,\d u\le n/(2\kappa)$. None of these depends on the type of the species, and the argument otherwise uses only \eqref{e.fixed-norms}, through $|P_s\mathsf m_u|\le\sqrt{n_s}$ and the resulting bounds
\begin{equation*}
|P_s(\sigma-\mathsf m_u)|\le2\sqrt{n_s},\qquad |P_s(2\sigma-\mathsf m_u-\mathsf m_w)|\le4\sqrt{n_s},\qquad\text{and}\qquad -3\lambda_{n,s}\le L_{u,s}\le\lambda_{n,s}.
\end{equation*}
Consequently, for every $\delta\in\{1,1/2,1/3,\ldots\}$, we obtain
\begin{equation}
|\mathscr E|\le\frac{C_{2,\xi}}2\int_0^1\sum_s\left(\frac{\E\tr(\sfC_{u,ss}^2)}{n^2}+\Var(L_{u,s})\right)\,\d u\le\frac{C_{2,\xi}}2\left(13\delta+\frac5{2\kappa n\delta^2}\right).
\label{e.mixed-residual-bound}
\end{equation}
Choosing $\delta:=1/\lfloor(\kappa n)^{1/3}\rfloor$ in \eqref{e.mixed-residual-bound} gives $|\mathscr E|\le15C_{2,\xi}(\kappa n)^{-1/3}$, which is \eqref{e.mixed-residual}.
\end{proof}

\subsection{The initial condition and the Hamilton--Jacobi comparison}
\label{s.HJ}

For $m \in [-1,1]^J$, we define $\psi_{\lambda,m} : \cQ_1^{\cS} \to \R$ by
\begin{equation}
\label{e.def.mixed-initial}
\psi_{\lambda,m}(q) := \sum_{j \in J} \lambda_j g_{m_j}(q_j) + \lambda_\Sp \psi^\Sp(q_\Sp),
\end{equation}
where $g_m$ and $\psi^\Sp$ are the functions of \eqref{e.def.gm} and \eqref{e.def.psiS}. We view $\cQ_2^{\cS}$ as a closed convex cone in the Hilbert space $L^2([0,1);\R^{\cS})$, whose inner product we denote by $\langle \cdot, \cdot \rangle_{L^2}$, and we write
\begin{equation}
\label{e.def.dual-cone}
(\cQ_2^{\cS})^* := \bigl\{ \kappa \in L^2([0,1);\R^{\cS}) \: : \: \langle \kappa, p \rangle_{L^2} \ge 0 \text{ for every } p \in \cQ_2^{\cS} \bigr\}
\end{equation}
for its dual cone.

\begin{lemma}[The initial condition]
\label{l.mixed-initial-limit}
For every $R > 0$, we have
\begin{align}\label{e.def.dn}
    \lim_{n \to \infty}\sup_a \ \sup_{q \in \cQ_2^{\cS} :\: |q|_{L^2} \le R} |\bar F_{n,a}(0,q) - \psi_{\lambda,m_{n,a}}(q)|= 0.
\end{align}
For every $m \in [-1,1]^J$, the function $\psi_{\lambda,m}$ is convex and $1$-Lipschitz for the $L^1$ norm on $\cQ_1^{\cS}$, and for every $q, q' \in \cQ_2^{\cS}$, we have
\begin{equation}
\label{e.initial-monotone}
q' - q \in (\cQ_2^{\cS})^* \quad \Longrightarrow \quad \psi_{\lambda,m}(q') \ge \psi_{\lambda,m}(q).
\end{equation}
\end{lemma}

\begin{proof}
Let $\bar F^\Is_{N,k}(0,\rho)$ and $\bar F^\Sp_N(0,\rho)$ denote the one-species versions of \eqref{e.def.enriched}, with reference measure $2^{-N}\sum_{\si\in\Omega_{N,k}}\boldsymbol{\delta}_\si$ and $\boldsymbol{\mu}_N$, respectively. The product decomposition \eqref{e.def.slice-measure} and \cite[(2.66) and Lemma~A.4(1)]{chen2026spherical} give
\begin{equation}
\bar F_{n,a}(0,q)=\sum_{j\in J}\lambda_{n,j}\bar F^\Is_{n_j,a_j}(0,q_j)+\lambda_{n,\Sp}\bar F^\Sp_{n_\Sp}(0,q_\Sp).
\label{e.mixed-initial-factorization}
\end{equation}
By \cite[Proposition~7.5]{boursier2026free}, for every $j \in J$ and $R > 0$, we have
\begin{equation*}
\lim_{n_j \to \infty} \ \sup_{0 \le a_j \le n_j} \ \sup_{\rho \in \cQ_1 :\: |\rho|_{L^1} \le R} \bigl| \bar F^\Is_{n_j,a_j}(0,\rho) - g_{m_{n,a,j}}(\rho) \bigr| = 0.
\end{equation*}
By \cite[Proposition~2.11]{chen2026spherical}, $\bar F^\Sp_{n_\Sp}(0,\rho)$ converges to $\psi^\Sp(\rho)$ as $n_\Sp \to \infty$, for every $\rho \in \cQ_1$. Since $L^2$-bounded subsets of $\cQ_2$ are relatively compact in $L^1$ \cite[Lemma~3.4]{chen2025vector}, and since $\bar F^\Sp_{n_\Sp}(0,\cdot)$ and $\psi^\Sp$ are $1$-Lipschitz for the $L^1$ norm \cite[Proposition~2.1 and Corollary~2.10]{chen2026spherical}, this convergence is uniform over $\{\rho \in \cQ_2 : |\rho|_{L^2} \le R\}$. By \eqref{e.mixed-initial-factorization} and \eqref{e.proportions}, it remains to compare $\psi_{\lambda,m}$ with the function $\psi_{\lambda_n,m}$ obtained by replacing $\lambda$ with $\lambda_n$ in \eqref{e.def.mixed-initial}. Since $g_m(0) = I(m) \le \log 2$, $\psi^\Sp(0) = 0$, and both functions are $1$-Lipschitz for the $L^1$ norm, we have, for every $m \in [-1,1]^J$ and $|q|_{L^2} \le R$,
\begin{equation*}
|\psi_{\lambda_n,m}(q)-\psi_{\lambda,m}(q)|\le\sum_s|\lambda_{n,s}-\lambda_s|(\log2+|q_s|_{L^1})\le(\log2+R)\sum_s|\lambda_{n,s}-\lambda_s|,
\end{equation*}
which tends to zero by \eqref{e.proportions}. This proves \eqref{e.def.dn}.

Convexity, the Lipschitz bound, and \eqref{e.initial-monotone} follow species by species from \cite[Corollary~6.11]{boursier2026free} for the Ising terms and \cite[Lemma~2.9, Corollary~2.10 and Lemma~5.1]{chen2026spherical} for the spherical term.
\end{proof}

For $m\in[-1,1]^J$, let $f_{\lambda,m}:[0,\infty)\times\cQ_2^{\cS}\to\R$ be given by the Hopf formula
\begin{equation}
\label{e.mixed-Hopf}
f_{\lambda,m}(t,q) = \sup_{\substack{p \in \cQ_\infty^{\cS} \\ 0 \le p \le \lambda}} \ \inf_{q' \in \cQ_\infty^{\cS}} \Ll\{ \psi_{\lambda,m}(q') + \langle q - q', p \rangle_{L^2} + t \int_0^1 \xi(p(u)) \, \d u \Rr\}.
\end{equation}
By Lemma~\ref{l.mixed-initial-limit} and \cite[Theorem~4.7]{chen2025HJ}, this is the unique Lipschitz viscosity solution of $\partial_t f - \int \xi(\partial_q f)=0 $ with $f(0,\cdot)=\psi_{\lambda,m}$.

\begin{theorem}
\label{t.mixed-slice}
Under \eqref{e.proportions}, for every $t \ge 0$ and $q \in \cQ_2^{\cS}$, we have
\begin{equation}
\label{e.mixed-slice-limit}
\lim_{n \to \infty} \sup_a |\bar F_{n,a}(t,q) - f_{\lambda,m_{n,a}}(t,q)| = 0.
\end{equation}
The convergence is uniform on compact subsets of $[0,\infty) \times \cQ_2^{\cS}$. For every $(t,q)$, the function $m \mapsto f_{\lambda,m}(t,q)$ is continuous on $[-1,1]^J$. 
\end{theorem}
\begin{proof}
We follow the proofs of \cite[Propositions~8.4 and~8.5]{boursier2026free}. Their inputs are the Lipschitz bound \eqref{e.mixed-Lipschitz}, the derivative representation \eqref{e.free-energy-derivatives} with $\dr_q\bar F_{n,a}\in\cQ_\infty^{\cS}$ and $0\le\dr_q\bar F_{n,a}\le\lambda_n$, the residual estimate \eqref{e.mixed-residual}, and the initial-condition properties in Lemma~\ref{l.mixed-initial-limit}. These have been verified for the mixed model, so the interpolation, comparison, and limit in those proofs apply without further changes. Let $L<\infty$, depending only on $\xi$, be the finite propagation constant in \cite[Lemma~8.2]{boursier2026free}. Writing the initial error in \cite[(8.16)]{boursier2026free} explicitly, their estimate \cite[(8.17)]{boursier2026free} gives, for every $\kappa\in[n^{-1},1]$ and $R>|q|_{L^2}+Lt$, 
\begin{equation}
\sup_a |\bar F_{n,a}(t,q) - f_{\lambda,m_{n,a}}(t,q)| \le \sup_{a,\,q'\in\cQ_2^{\cS}:\:|q'|_{L^2}\le R}|\bar F_{n,a}(0,q')-\psi_{\lambda,m_{n,a}}(q')| + 3 \kappa + 15 \, C_{2,\xi} \, t \, (\kappa n)^{-1/3}.
\label{e.mixed-comparison}
\end{equation}
Setting $\kappa:=n^{-1/4}\in[n^{-1},1]$ and using Lemma~\ref{l.mixed-initial-limit} proves \eqref{e.mixed-slice-limit}. On a compact subset of $[0,\infty)\times\cQ_2^{\cS}$, both $t$ and $|q|_{L^2}$ are bounded, so one choice of $R$ in \eqref{e.mixed-comparison} proves uniform convergence.

For continuity in $m$, put $I_\lambda(m):=\sum_{j\in J}\lambda_jI(m_j)$. The spherical term in \eqref{e.def.mixed-initial} is independent of $m$. Thus \eqref{e.gm-continuity} supplies the same initial-data estimate as \cite[(6.52)]{boursier2026free}, and the application of \cite[Lemma~8.2]{boursier2026free} leading to \cite[(8.21)]{boursier2026free} gives
\begin{equation}
\label{e.mixed-m-continuity}
|f_{\lambda,m}(t,q)-f_{\lambda,m'}(t,q)-I_\lambda(m)+I_\lambda(m')|\le2(|q|_{L^2}+Lt)|m-m'|_\infty.
\end{equation}
Since $I$ is continuous on $[-1,1]$, this proves the assertion.
\end{proof}

%
%
%
%
%
%
\section{Reduction of the Hopf formula to one overlap law}
\label{s.reduction}
We now specialize Theorem~\ref{t.mixed-slice} to the bipartite model that is relevant to the study of \eqref{e.def.hybrid}. As discussed in the introduction, our analysis strategy for this model involves the identification of $2^k + 1$ types of spins. In this section, we record how the Hopf formula for this model can then be considerably simplified, ultimately ending up as a variational problem that involves a single probability measure.

\subsection{The bilinear covariance function and the path conjugates}
\label{s.bilinear}

Throughout this section, we fix $k \in \N \cup \{0\}$ and $\al > 0$. Let $E_k := \{-1,1\}^k$, with $E_0$ a singleton, and let $\epsav$ denote the uniform average over $\ep \in E_k$. The Ising species will be indexed by $E_k$: in Section~\ref{s.gaussian}, the spins of species $\ep$ are those at which the $k$ retained patterns take the values $\ep$. We set
\begin{equation}
\label{e.def.species-k}
\cS_k := E_k \cup \{\Sp\}, \qquad \lambda_\ep := \frac{2^{-k}}{1+\al}, \qquad \lambda_\Sp := \frac{\al}{1+\al},\qquad\text{and}\qquad \xi_{\al,k}(a) := (1+\al) \, a_\Sp \sum_{\ep \in E_k} a_\ep,
\end{equation}
for $a \in \R^{\cS_k}$.
For $m = (m_\ep)_{\ep \in E_k} \in [-1,1]^{E_k}$, we define the initial condition $\psi_{\al,k,m} : \cQ_1^{\cS_k} \to \R$ by
\begin{equation}
\label{e.def.psi-k}
\psi_{\al,k,m}(q) := \frac{1}{1+\al} \epsav \, g_{m_\ep}(q_\ep) + \frac{\al}{1+\al} \psi^\Sp(q_\Sp),
\end{equation}
which is \eqref{e.def.mixed-initial} in this setting. With the proportions and the covariance function as in \eqref{e.def.species-k} and this initial condition, the Hopf formula \eqref{e.mixed-Hopf} at $q = 0$ defines $f_{\al,k} : [0,\infty) \times [-1,1]^{E_k} \to \R$ by
\begin{equation}
\label{e.def.f-alpha-k}
f_{\al,k}(t,m) := \sup_{\substack{p' \in \cQ_\infty^{\cS_k} \\ 0 \le p' \le \lambda}} \ \inf_{q' \in \cQ_\infty^{\cS_k}} \Ll\{ \psi_{\al,k,m}(q') - \langle q', p' \rangle_{L^2} + t \int_0^1 \xi_{\al,k}(p'(u)) \, \d u \Rr\}.
\end{equation}

For $m \in [-1,1]$, we define the convex conjugate $g_m^* : \cQ_\infty \to \R \cup \{+\infty\}$ of $g_m$ by
\begin{equation}
\label{e.def.gm-conj}
g_m^*(p) := \sup_{\rho \in \cQ_\infty} \Ll\{ \int_0^1 \rho(u) p(u) \, \d u - g_m(\rho) \Rr\}.
\end{equation}
By definition, $g_m(\rho) \ge \int_0^1 \rho p - g_m^*(p)$ for every $p \in \cQ_\infty$, and the path $z$ in \eqref{e.gm-supporting} gives equality at $p = z$, with $0 \le z \le 1$. Hence, for every $\rho \in \cQ_\infty$,
\begin{equation}
\label{e.gm-biconjugacy}
g_m(\rho) = \max_{\substack{p \in \cQ_\infty \\ 0 \le p \le 1}} \Ll\{ \int_0^1 \rho(u) p(u) \, \d u - g_m^*(p) \Rr\};
\end{equation}
this is also \cite[Lemma~9.1]{boursier2026free}. Finally, for $t \ge 0$, we set $a_t := 2t/(1+2t)$, and we recall the notation $Q_{\boldsymbol{\zeta}}$ for the quantile function defined just below \eqref{e.def.path-law}.

\begin{proposition}[The slice free energy as a maximum over one law]
\label{p.one-law-slice}
For every $t \ge 0$ and $m \in [-1,1]^{E_k}$, we have
\begin{equation}
\label{e.one-law-slice}
(1+\al) f_{\al,k}(t,m) = \sup_{\boldsymbol{\zeta} \in \cP([0,a_t])} \Ll\{ \epsav \, g_{m_\ep}(\al t Q_{\boldsymbol{\zeta}}) - \al \cE(\boldsymbol{\zeta}) \Rr\}.
\end{equation}
The supremum is attained, and its value is unchanged if the set $\cP([0,a_t])$ is replaced by $\cP_{<1}$.
\end{proposition}
\begin{proof}
\emph{Step 1: conjugate form of the Hopf formula.} In this step, we rewrite the Hopf formula in a form involving the convex conjugates of the initial condition, as in \cite[Lemma~9.2]{boursier2026free}. For $p' = (p'_s)_{s \in \cS_k}$ in \eqref{e.def.f-alpha-k}, we write $p'_\ep = \lambda_\ep p_\ep$ for $\ep \in E_k$ and $p'_\Sp = \lambda_\Sp p_\Sp$, with $p_\ep, p_\Sp \in \cQ_\infty$ and $0 \le p_\ep, p_\Sp \le 1$. For fixed $(p_\ep)$ and $p_\Sp$, the terms in \eqref{e.def.f-alpha-k} involving $q'$ give
\begin{align*}
& \inf_{q' \in \cQ_\infty^{\cS_k}} \Ll\{ \psi_{\al,k,m}(q') - \sum_\ep \lambda_\ep \langle q'_\ep, p_\ep \rangle_{L^2} - \lambda_\Sp \langle q'_\Sp, p_\Sp \rangle_{L^2} \Rr\} \\
& \qquad \stackrel{\eqref{e.def.psi-k}}{=} \sum_\ep \lambda_\ep \inf_{\rho \in \cQ_\infty} \Ll\{ g_{m_\ep}(\rho) - \langle \rho, p_\ep \rangle_{L^2} \Rr\} + \lambda_\Sp \inf_{\rho \in \cQ_\infty} \Ll\{ \psi^\Sp(\rho) - \langle \rho, p_\Sp \rangle_{L^2} \Rr\} \\
& \qquad \stackrel{\eqref{e.def.gm-conj},\,\eqref{e.def.psiS-conj}}{=} - \sum_\ep \lambda_\ep g_{m_\ep}^*(p_\ep) - \lambda_\Sp \psi^{\Sp,*}(p_\Sp).
\end{align*}
The interaction term satisfies
\begin{equation*}
(1+\al) t \int_0^1 \xi_{\al,k} \bigl( (\lambda_\ep p_\ep(u))_{\ep \in E_k}, \lambda_\Sp p_\Sp(u) \bigr) \, \d u \stackrel{\eqref{e.def.species-k}}{=} (1+\al)^2 t \lambda_\Sp \sum_\ep \lambda_\ep \langle p_\Sp, p_\ep \rangle_{L^2} \stackrel{\eqref{e.def.species-k}}{=} \al t \langle p_\Sp, \epsav p_\ep \rangle_{L^2}.
\end{equation*}
Multiplying \eqref{e.def.f-alpha-k} by $1+\al$ and substituting these two identities, we obtain
\begin{equation}
\label{e.bilinear-conjugates}
(1+\al) f_{\al,k}(t,m) = \sup_{\substack{p_\Sp, p_\ep \in \cQ_\infty \\ 0 \le p_\Sp, p_\ep \le 1}} \Ll\{ \al t \langle p_\Sp, \epsav p_\ep \rangle_{L^2} - \epsav \, g^*_{m_\ep}(p_\ep) - \al \psi^{\Sp,*}(p_\Sp) \Rr\}.
\end{equation}

\emph{Step 2: a uniform bound on the endpoint of $p_\Sp$.} Fix $(p_\ep)$ as in \eqref{e.bilinear-conjugates}, and set $\rho := t \epsav p_\ep \in \cQ_\infty$. The terms depending on $p_\Sp$ in \eqref{e.bilinear-conjugates} satisfy
\begin{align*}
\sup_{\substack{p_\Sp \in \cQ_\infty \\ 0 \le p_\Sp \le 1}} \Ll\{ \al t \langle p_\Sp, \epsav p_\ep \rangle_{L^2} - \al \psi^{\Sp,*}(p_\Sp) \Rr\}
& = \al \sup_{\substack{p_\Sp \in \cQ_\infty \\ 0 \le p_\Sp \le 1}} \Ll\{ \langle \rho, p_\Sp \rangle_{L^2} - \psi^{\Sp,*}(p_\Sp) \Rr\} \\
& \stackrel{\eqref{e.spherical-biconjugacy}}{=} \al \psi^\Sp(\rho) \stackrel{\eqref{e.spherical-fenchel}}{=} \al \Ll\{ \langle \rho, p_\rho \rangle_{L^2} - \psi^{\Sp,*}(p_\rho) \Rr\}.
\end{align*}
Thus the supremum is attained at $p_\Sp = p_\rho$. Since $\rho \le t$, we have, for every $u \in [0,1]$,
\begin{equation*}
p_\rho(u) \le p_\rho(1) \stackrel{\eqref{e.spherical-derivative}}{\le} \frac{2 \rho(1)}{1 + 2 \rho(1)} \le \frac{2t}{1+2t} = a_t.
\end{equation*}
Restricting \eqref{e.bilinear-conjugates} to $p_\Sp \le a_t$ therefore leaves its value unchanged.

\emph{Step 3: elimination of the Ising dual paths.} For every fixed $p_\Sp \in \cQ_\infty$ with $0 \le p_\Sp \le 1$, the maximization over the Ising dual paths separates as follows:
\begin{equation}
\label{e.Ising-path-elimination}
\sup_{\substack{p_\ep \in \cQ_\infty \\ 0 \le p_\ep \le 1}} \epsav \Ll\{ \al t \langle p_\Sp, p_\ep \rangle_{L^2} - g_{m_\ep}^*(p_\ep) \Rr\} = \epsav \max_{\substack{p \in \cQ_\infty \\ 0 \le p \le 1}} \Ll\{ \langle \al t p_\Sp, p \rangle_{L^2} - g_{m_\ep}^*(p) \Rr\} \stackrel{\eqref{e.gm-biconjugacy}}{=} \epsav \, g_{m_\ep}(\al t p_\Sp).
\end{equation}
Combining Step~2 with \eqref{e.bilinear-conjugates} and \eqref{e.Ising-path-elimination} gives
\begin{equation}
\label{e.spherical-path-form}
(1+\al) f_{\al,k}(t,m) = \sup_{\substack{p_\Sp \in \cQ_\infty \\ 0 \le p_\Sp \le a_t}} \Ll\{ \epsav \, g_{m_\ep}(\al t p_\Sp) - \al \psi^{\Sp,*}(p_\Sp) \Rr\}.
\end{equation}
Every such $p_\Sp$ is the quantile function of its law $\boldsymbol{\nu}_{p_\Sp} \in \cP([0,a_t])$. Conversely, every law in $\cP([0,a_t])$ has such a quantile path. For these paths, we have
\begin{equation*}
\epsav \, g_{m_\ep}(\al t p_\Sp) - \al \psi^{\Sp,*}(p_\Sp) \stackrel{\eqref{e.spherical-entropy-duality}}{=} \epsav \, g_{m_\ep}(\al t Q_{\boldsymbol{\nu}_{p_\Sp}}) - \al \cE(\boldsymbol{\nu}_{p_\Sp}).
\end{equation*}
This identifies the supremum in \eqref{e.spherical-path-form} with the one in \eqref{e.one-law-slice}. By part (4) of Lemma~\ref{l.Ising-data} and Lemma~\ref{l.spherical-duality}, the functional $\boldsymbol{\zeta} \mapsto \epsav \, g_{m_\ep}(\al t Q_{\boldsymbol{\zeta}}) - \al \cE(\boldsymbol{\zeta})$ is continuous on the compact set $\cP([0,a_t])$, so its supremum is attained.

To extend the domain to $\cP_{<1}$, we instead apply \eqref{e.Ising-path-elimination} directly to \eqref{e.bilinear-conjugates}. This gives
\begin{equation}
\label{e.unrestricted-spherical-path-form}
(1+\al) f_{\al,k}(t,m) = \sup_{\substack{p_\Sp \in \cQ_\infty \\ 0 \le p_\Sp \le 1}} \Ll\{ \epsav \, g_{m_\ep}(\al t p_\Sp) - \al \psi^{\Sp,*}(p_\Sp) \Rr\}.
\end{equation}
Since $a_t < 1$, the admissible sets satisfy
\begin{equation*}
\{ p_\Sp \in \cQ_\infty : p_\Sp \le a_t \} \subseteq \{ p_\Sp \in \cQ_\infty : p_\Sp(1) < 1 \} \subseteq \{ p_\Sp \in \cQ_\infty : p_\Sp \le 1 \}.
\end{equation*}
The suprema on the two outer sets agree by \eqref{e.spherical-path-form} and \eqref{e.unrestricted-spherical-path-form}, so the supremum on the middle set has the same value. Applying \eqref{e.spherical-entropy-duality} to each path in the middle set identifies this supremum with the one over $\cP_{<1}$.
\end{proof}

\subsection{Releasing the magnetization constraints}
\label{s.minimax}
Recall the definition of $\msc P_a(\boldsymbol{\zeta},h)$ in \eqref{e.def.Parisi-radius}. Using Proposition~\ref{p.one-law-slice} and a min-max inversion, one can compute the convex conjugate of $(1+\al) f_{\al,k}(t,m)$ with respect to $m$ and relate it to $\msc P_a(\boldsymbol{\zeta},h)$. More precisely, the following holds.

\begin{proposition}[Minimax over the magnetizations]
\label{p.law-minimax}
For every $t \ge 0$ and every family of fields $(h_\ep)_{\ep \in E_k} \in \R^{E_k}$, we have
\begin{equation}
\label{e.law-minimax}
\al t + \max_{m \in [-1,1]^{E_k}} \Ll\{ \epsav [h_\ep m_\ep] - (1+\al) f_{\al,k}(t,m) \Rr\} = \inf_{\boldsymbol{\zeta} \in \cP_{<1}} \Ll\{ \epsav \, \msc P_{\al t}(\boldsymbol{\zeta},h_\ep) + \al \cE(\boldsymbol{\zeta}) \Rr\}.
\end{equation}
The function minimized on the right-hand side is strictly convex on $\cP_{<1}$. The infimum is attained at a unique $\boldsymbol{\zeta}$, and this minimizer belongs to $\cP([0,a_t])$.
\end{proposition}
\begin{proof}
\emph{Step 1: minimax on the compact set.} We define $L : [-1,1]^{E_k} \times \cP_{<1} \to \R$ by
\begin{equation}
\label{e.def.minimax-L}
L(m,\boldsymbol{\zeta}) := \epsav \bigl[ g_{m_\ep}(\al t Q_{\boldsymbol{\zeta}}) - h_\ep m_\ep \bigr] - \al \cE(\boldsymbol{\zeta}).
\end{equation}
All minima and maxima indexed by $m$ in this proof are over $[-1,1]^{E_k}$. We first restrict the law to $\cP([0,a_t])$. For fixed $\rho \in \cQ_1$, we write $g_m(\rho) = [g_m(\rho) - I(m)] + I(m)$. Both summands are convex and continuous in $m$ by part (2) of Lemma~\ref{l.Ising-data} and \eqref{e.def.I}, respectively. Thus $L$ is convex and continuous in $m$. By part (4) of Lemma~\ref{l.Ising-data} and Lemma~\ref{l.spherical-duality}, it is concave and continuous in~$\boldsymbol{\zeta}$. Both domains are compact and convex, so Sion's minimax theorem \cite[Theorem~3.4]{sion1958general} gives
\begin{equation}
\label{e.Sion}
\min_m \max_{\boldsymbol{\zeta} \in \cP([0,a_t])} L(m,\boldsymbol{\zeta}) = \max_{\boldsymbol{\zeta} \in \cP([0,a_t])} \min_m L(m,\boldsymbol{\zeta}).
\end{equation}
This is the exchange performed in the proof of \cite[Theorem~10.2]{boursier2026free}, with the affine pairing against a potential replaced here by the concave function $-\al \cE$.

For each fixed $m$, we have
\begin{equation}
\label{e.minimax-fixed-m}
\max_{\boldsymbol{\zeta} \in \cP([0,a_t])} L(m,\boldsymbol{\zeta}) \stackrel{\eqref{e.def.minimax-L},\,\eqref{e.one-law-slice}}{=} (1+\al) f_{\al,k}(t,m) - \epsav [h_\ep m_\ep].
\end{equation}
For each fixed $\boldsymbol{\zeta} \in \cP_{<1}$, the minimization separates over the finitely many scalar magnetizations, giving
\begin{align}
\min_m L(m,\boldsymbol{\zeta}) &\stackrel{\eqref{e.def.minimax-L}}{=} \epsav \min_{x \in [-1,1]} \Ll\{ g_x(\al t Q_{\boldsymbol{\zeta}}) - h_\ep x \Rr\} - \al \cE(\boldsymbol{\zeta})
\notag
\\
\label{e.minimax-fixed-law}
&\stackrel{\eqref{e.field-duality}}{=} \epsav \, \psi^\Is(\al t Q_{\boldsymbol{\zeta}}; h_\ep) - \al \cE(\boldsymbol{\zeta}) \stackrel{\eqref{e.field-duality}}{=} \al t - \epsav \, \msc P_{\al t}(\boldsymbol{\zeta},h_\ep) - \al \cE(\boldsymbol{\zeta}).
\end{align}
Combining the previous displays, we obtain
\begin{align}
\notag
& \al t + \max_m \Ll\{ \epsav [h_\ep m_\ep] - (1+\al) f_{\al,k}(t,m) \Rr\} \stackrel{\eqref{e.minimax-fixed-m}}{=} \al t - \min_m \max_{\boldsymbol{\zeta} \in \cP([0,a_t])} L(m,\boldsymbol{\zeta}) \\
\label{e.law-minimax-compact}
& \qquad \stackrel{\eqref{e.Sion}}{=} \al t - \max_{\boldsymbol{\zeta} \in \cP([0,a_t])} \min_m L(m,\boldsymbol{\zeta}) \stackrel{\eqref{e.minimax-fixed-law}}{=} \min_{\boldsymbol{\zeta} \in \cP([0,a_t])} \Ll\{ \epsav \, \msc P_{\al t}(\boldsymbol{\zeta},h_\ep) + \al \cE(\boldsymbol{\zeta}) \Rr\}.
\end{align}

\emph{Step 2: the domain $\cP_{<1}$.} By Proposition~\ref{p.one-law-slice}, we have $\sup_{\boldsymbol{\zeta} \in \cP_{<1}} L(m,\boldsymbol{\zeta}) = \max_{\boldsymbol{\zeta} \in \cP([0,a_t])} L(m,\boldsymbol{\zeta})$ for every $m$. The elementary minimax inequality therefore gives
\begin{equation*}
\sup_{\boldsymbol{\zeta} \in \cP_{<1}} \inf_m L(m,\boldsymbol{\zeta}) \le \inf_m \sup_{\boldsymbol{\zeta} \in \cP_{<1}} L(m,\boldsymbol{\zeta}) \stackrel{\eqref{e.Sion}}{=} \max_{\boldsymbol{\zeta} \in \cP([0,a_t])} \min_m L(m,\boldsymbol{\zeta}),
\end{equation*}
and the reverse inequality holds since $\cP([0,a_t]) \subseteq \cP_{<1}$. The two domains thus give the same value. Subtracting this value from $\al t$ and using \eqref{e.minimax-fixed-law} and \eqref{e.law-minimax-compact} proves \eqref{e.law-minimax}.

\emph{Step 3: strict convexity and uniqueness.} The first identity in \eqref{e.field-duality} represents $\boldsymbol{\zeta} \mapsto \psi^\Is(\al t Q_{\boldsymbol{\zeta}};h)$ as an infimum of concave functions of $\boldsymbol{\zeta}$ by \eqref{e.gm-law-concavity}. Thus this function is concave, and the second identity in \eqref{e.field-duality} shows that $\msc P_{\al t}(\cdot,h)$ is convex. Since $\al > 0$ and $\cE$ is strictly convex by Lemma~\ref{l.spherical-duality}, the functional minimized in \eqref{e.law-minimax} is strictly convex. By \eqref{e.P-law-continuity} and Lemma~\ref{l.spherical-duality}, it is continuous on the compact set $\cP([0,a_t])$, so its minimum is attained there and, by Step~2, equals the infimum over $\cP_{<1}$. Any two laws in $\cP_{<1}$ belong to a common $\cP([0,b])$ with $b < 1$, so strict convexity gives uniqueness on $\cP_{<1}$ as well.
\end{proof}

%
%
%
%
%
%
\section{Gaussian patterns with finitely many Bernoulli fields}
\label{s.gaussian}

In this section, we prove Theorem~\ref{t.hybrid}. Throughout, $\al, \be > 0$ and $k \in \N \cup \{0\}$ are fixed. We recall the random field $\cH_y$ from \eqref{e.def.field-law}, and the function $f_{\al,k}$ from Section~\ref{s.bilinear}. We first identify the free energy of the model in which the auxiliary Gaussian variables are constrained to a sphere, then integrate over the radius, and finally compare the contributions of the retained fields.

\subsection{The free energy on a sphere}
\label{s.shell}

For $u \ge 0$ and $y \in \R^k$, we define
\begin{equation}
\label{e.def.shell}
\bar B^{(k)}_{N,M}(u,y) := \frac 1 N \E \log \iint \exp \Ll( \sqrt{\frac u N} \sum_{i=1}^N \sum_{\mu=1}^M g_i^\mu \si_i \tau_\mu + \sqrt \be \sum_{i=1}^N \sum_{j=1}^k y_j \eta_i^j \si_i \Rr) \, \d P_N(\si) \, \d \boldsymbol{\mu}_M(\tau),
\end{equation}
where the expectation is over the Gaussian matrix $G = (g_i^\mu)_{i \le N, \mu \le M}$ and the Bernoulli variables $(\eta_i^j)$, and $\boldsymbol{\mu}_M$ is the uniform probability measure on the sphere $\Sph_M$. 
We use the notation $v_k(m) := \epsav[\ep \, m_\ep] \in \R^k$ for the vector of magnetizations along the retained patterns. For $k = 0$, $v_0(m)$ is the empty vector. We also recall that the random variable $\cH_y$ is defined in \eqref{e.def.field-law}.

\begin{proposition}[Free energy on a sphere]
\label{p.shell}
Under \eqref{e.loading}, the functions $\bar B^{(k)}_{N,M_N}$ converge locally uniformly on $[0,\infty) \times \R^k$ to the continuous function
\begin{equation}
\label{e.shell-limit}
B^{(k)}_\al(u,y) := \inf_{\boldsymbol{\zeta} \in \cP_{<1}} \Ll\{ \E \, \msc P_{\al u/2}(\boldsymbol{\zeta},\cH_y) + \al \cE(\boldsymbol{\zeta}) \Rr\}.
\end{equation}
The infimum is attained at a unique $\boldsymbol{\zeta}$, which belongs to $\cP([0,u/(1+u)])$. Equivalently,
\begin{equation}
\label{e.shell-slices}
B^{(k)}_\al(u,y) = \frac{\al u}{2} + \max_{m \in [-1,1]^{E_k}} \Ll\{ \sqrt \be \, y \cdot v_k(m) - (1+\al) f_{\al,k}(u/2,m) \Rr\}.
\end{equation}
\end{proposition}
\begin{proof}
We fix $u \ge 0$ and $y \in \R^k$, and put $M := M_N$. 

\emph{Step 1: deterministic numbers of field types.} The field felt by $\si_i$ is $h_\ep := \sqrt\be\,y\cdot\ep$ when $(\eta_i^1,\ldots,\eta_i^k)=\ep\in E_k$. Let $N_\ep$ be the number of sites of type $\ep$. We choose deterministic integers $\bar N_\ep$ that sum to $N$ and differ from $N2^{-k}$ by at most one. Fix an ordering of $E_k$ and partition $\{1,\ldots,N\}$ into consecutive blocks $(I_\ep)_{\ep \in E_k}$ of sizes $\bar N_\ep$ in this order. We assign type $\ep$, and hence field $h_\ep$, to every site in $I_\ep$. Let $\widetilde B_{N,M}(u,y)$ be defined as in \eqref{e.def.shell}, but with the term $\sqrt \be \sum_{i,j} y_j \eta_i^j \si_i$ replaced by $\sum_\ep h_\ep \sum_{i \in I_\ep} \si_i$, so that the expectation is only over $G$.

The variance of $N_\ep$ is at most $N$, so $\E\sum_\ep|N_\ep-\bar N_\ep|\le 2^k(\sqrt N+1)$. Changing at most $\frac12\sum_\ep|N_\ep-\bar N_\ep|$ types gives exactly $\bar N_\ep$ sites of each type, and each change alters the logarithm of the partition function by at most $2\sqrt{\be k}|y|$. A permutation of the sites then puts the types in the blocks $I_\ep$. Since the law of $G$ is invariant under permutations of its rows, this last operation does not change the expectation over $G$ of the logarithm of the partition function. Dividing by $N$ and averaging over the types gives
\begin{equation}
\label{e.shell-type-replacement}
\bigl|\bar B^{(k)}_{N,M}(u,y)-\widetilde B_{N,M}(u,y)\bigr| \le C_{k,\be}|y|N^{-1/2}, \qquad u \ge 0.
\end{equation}
Thus it suffices to determine the limit of $\widetilde B_{N,M}(u,y)$.

\emph{Step 2: identification of the free energy of a slice.} Let $n := N + M$, and consider the model of Section~\ref{s.mixed} with the Ising species $\ep \in E_k$ on the blocks $I_\ep$, of sizes $\bar N_\ep$, and a spherical species of dimension $M$. The Hamiltonian
$H_n(\si,\tau) := \sqrt{\frac{1+\al}{n}} \sum_{i,\mu} g_i^\mu \si_i \tau_\mu$
has covariance $n \xi_{\al,k}(R)$, with $\xi_{\al,k}$ as in \eqref{e.def.species-k}. Let $t_n(u) := un/[2N(1+\al)]$. A direct substitution gives
\begin{equation}
\label{e.shell-normalization}
\sqrt{2 t_n(u)} H_n = \sqrt{\frac u N} \sum_{i,\mu} g_i^\mu \si_i \tau_\mu, \qquad 
n t_n(u) \xi_{\al,k}(\lambda_n) = \frac{M u}{2}, \qquad \text{and} \qquad 
\lim_{N \to \infty} t_n(u) = \frac u 2.
\end{equation}
We take $N$ large enough that every $\bar N_\ep$ is positive. For each $a=(a_\ep)_{\ep\in E_k}$ with $a_\ep\in\{0,\ldots,\bar N_\ep\}$, recall the magnetizations introduced before \eqref{e.def.slice-measure}, which here read $m_{n,a,\ep}=\frac{2a_\ep}{\bar N_\ep}-1$ and $m_{n,a}=(m_{n,a,\ep})_{\ep\in E_k}$, and the slice measure $\boldsymbol{\mu}_{n,a}$ of \eqref{e.def.slice-measure}, which here is a measure on $\{-1,1\}^N \times \Sph_M$.
Let~$B_a$ be the random variable obtained from the definition of $\widetilde B_{N,M}(u,y)$ by replacing $P_N \otimes \boldsymbol{\mu}_M$ with~$\boldsymbol{\mu}_{n,a}$ and omitting the expectation. On the slice, the field term is the constant $\sum_\ep\bar N_\ep h_\ep m_{n,a,\ep}$. Equations \eqref{e.shell-normalization} and \eqref{e.def.enriched} at $q=0$ give
\begin{equation}
\label{e.shell-slice-mean}
\E B_a = \frac{Mu}{2N} + \sum_\ep \frac{\bar N_\ep}{N}h_\ep m_{n,a,\ep} - \frac n N\bar F_{n,a}(t_n(u),0).
\end{equation}
As $N \to \infty$, we have $\bar N_\ep/N \to 2^{-k}$ and $n/N \to 1+\al$, so that \eqref{e.proportions} holds with the proportions of \eqref{e.def.species-k}. By \eqref{e.mixed-Lipschitz}, we may replace $t_n(u)$ by $u/2$ in \eqref{e.shell-slice-mean} at the price of an error that vanishes uniformly in $a$. Theorem~\ref{t.mixed-slice} at $t = u/2$ and $q = 0$, together with the identity $\epsav [h_\ep m_\ep] = \sqrt\be\,y\cdot v_k(m)$, then yields the uniform slice limit 
\begin{equation}
\label{e.shell-slice-limit}
\lim_{N \to \infty}\sup_a \biggl|\E B_a - \frac{\al u}{2} - \sqrt\be\,y\cdot v_k(m_{n,a}) + (1+\al)f_{\al,k}(u/2,m_{n,a})\biggr| = 0.
\end{equation}

\emph{Step 3: sum over the slices.} Since $\sum_a \boldsymbol{\mu}_{n,a} = P_N \otimes \boldsymbol{\mu}_M$, we have $\widetilde B_{N,M}(u,y) = N^{-1} \E \log \sum_a e^{N B_a}$. We now follow the gluing argument of \cite[Lemma~7.1]{boursier2026free}, in the case $q=0$, with only $G$ random after Step~1. Since $|\si|^2=N$ and $|\tau|^2=M$, each $B_a$ is $\sqrt{uM}/N$-Lipschitz in $G$. For $u>0$ and $s>0$, Gaussian concentration and a union bound over at most $(N+1)^{2^k}$ slices give $\P(\max_a|B_a-\E B_a|>s)\le 2(N+1)^{2^k}\exp(-N^2s^2/(2uM))$. Integrating this tail bound gives $\E\max_a|B_a-\E B_a| \le C \sqrt{uM \, 2^k \log(N+1)}/N$, which also holds at $u=0$, when all $B_a$ are deterministic. Since $\max_a \E B_a \le \E \max_a B_a \le \max_a \E B_a + \E\max_a|B_a-\E B_a|$ and $\max_a B_a \le N^{-1} \log \sum_a e^{NB_a} \le \max_a B_a+2^k\log(N+1)/N$, we obtain
\begin{equation}
0 \le \widetilde B_{N,M}(u,y)-\max_a\E B_a \le \frac{2^k\log(N+1)}{N} + C\sqrt{\frac{uM}{N^2}\,2^k\log(N+1)}.
\label{e.shell-slice-sum}
\end{equation}
Since every $\bar N_\ep$ tends to infinity, the grids of magnetizations become dense in $[-1,1]^{E_k}$. Combining \eqref{e.shell-type-replacement}, \eqref{e.shell-slice-limit} and \eqref{e.shell-slice-sum} with the continuity in $m$ from \eqref{e.mixed-m-continuity}, we obtain
\begin{equation*}
\lim_{N \to \infty}\bar B^{(k)}_{N,M_N}(u,y) = \frac{\al u}{2} + \max_{m \in [-1,1]^{E_k}}\Ll\{\sqrt\be\,y\cdot v_k(m)-(1+\al)f_{\al,k}(u/2,m)\Rr\}.
\end{equation*}
This is the expression in \eqref{e.shell-slices}, as in the proof of \cite[Theorem~10.1]{boursier2026free}. Since $\cH_y$ has the law of $h_\ep$ with $\ep$ uniform on $E_k$, Proposition~\ref{p.law-minimax} with $t=u/2$ and $h_\ep=\sqrt\be\,y\cdot\ep$ identifies this expression with \eqref{e.shell-limit}, and proves the uniqueness of the minimizer and the bound $a_{u/2}=u/(1+u)$ on its support.

\emph{Step 4: local uniformity.} By Gaussian integration by parts, the derivative of $\bar B^{(k)}_{N,M}(u,y)$ in $u$ equals $\frac{1}{2N^2} \E \bigl\langle |\si|^2 |\tau|^2 - (\si \cdot \si')(\tau \cdot \tau') \bigr\rangle$, where $\langle \cdot \rangle$ is the Gibbs measure in \eqref{e.def.shell} and $(\si',\tau')$ is an independent copy of $(\si,\tau)$ under it. Since $|\si|^2 = N$ and $|\tau|^2 = M$, this derivative lies in $[0,M/N]$, and thus
\begin{equation*}
|\bar B^{(k)}_{N,M}(u,y) - \bar B^{(k)}_{N,M}(u',y)| \le \frac M N |u - u'|.
\end{equation*}
For each realization of the $\eta_i^j$, changing $y$ into $y'$ changes the logarithm of the partition function by at most $\sqrt \be \sum_i |(y-y') \cdot \eta_i|$, where $\eta_i := (\eta_i^1,\ldots,\eta_i^k)$. Taking expectations, dividing by $N$, and using $\epsav [(y-y') \cdot \ep]^2 = |y - y'|^2$, we obtain
\begin{equation}
\label{e.hybrid-field-Lipschitz}
|\bar B^{(k)}_{N,M}(u,y) - \bar B^{(k)}_{N,M}(u,y')| \le \sqrt \be \, \epsav |(y - y') \cdot \ep| \le \sqrt \be |y - y'|.
\end{equation}
The same estimate holds for $\bar F^{(k)}_{N,M}$. Since $M/N \to \al$, these Lipschitz bounds pass to the limit, so $B^{(k)}_\al$ is continuous, and a finite-net argument proves the local uniform convergence.
\end{proof}

\subsection{Integration over the radius}
\label{s.radial}

Recall that $\boldsymbol{\gamma}_M$ is the standard Gaussian probability measure on $\R^M$. The Hubbard--Stratonovich identity gives 
\begin{equation}
\label{e.HS}
\int e^{\be |G^\top \si|^2/(2N) + \sum_i h_i \si_i} \, \d P_N(\si) = \iint e^{\sqrt{\be/N} \, \si^\top G z + \sum_i h_i \si_i} \, \d P_N(\si) \, \d \boldsymbol{\gamma}_M(z).
\end{equation}
We write $z = r \tau$ with $r \ge 0$ and $\tau \in \Sph_M$, so that $|z| = r \sqrt M$. Under $\boldsymbol{\gamma}_M$, the variable $\tau$ is uniform on $\Sph_M$ conditionally on $r$, and the density of $r$ is
\begin{equation}
\label{e.radial-density}
\rho_M(r) := \frac{2^{1-M/2} M^{M/2}}{\Gamma(M/2)} \, r^{M-1} e^{-M r^2/2}, \qquad r > 0.
\end{equation}
We recall $\bar F^{(k)}_{N,M}(y)$ from \eqref{e.def.hybrid}, and let $B^{(k)}_{N,M}$ denote the random free energy inside the expectation in \eqref{e.def.shell}. Since $\sqrt{\be/N} \, \si^\top G z = \sqrt{\be r^2/N} \sum_{i,\mu} g_i^\mu \si_i \tau_\mu$, the identities \eqref{e.HS} and \eqref{e.radial-density} give
\begin{equation}
\label{e.radial-integral}
\bar F^{(k)}_{N,M}(y) = \E \Ll[ \frac 1 N \log \int_0^\infty \exp \bigl( N B^{(k)}_{N,M}(\be r^2, y) \bigr) \rho_M(r) \, \d r \Rr].
\end{equation}

\begin{proof}[Proof of Theorem~\ref{t.hybrid}]
We fix $y\in\R^k$ and put $M:=M_N$. For every Borel set $I\subset(0,\infty)$, we write
\begin{equation}
Z_N(I,y):=\int_I\exp\bigl(NB^{(k)}_{N,M}(\be r^2,y)\bigr)\rho_M(r)\,\d r.
\label{e.radial-restricted-Z}
\end{equation}
Thus \eqref{e.radial-integral} reads $\bar F^{(k)}_{N,M}(y)=N^{-1}\E\log Z_N((0,\infty),y)$.

\emph{Step 1: compact intervals of radii.} 
We fix $0<r_0<1<R<\infty$, and recall $\cR(r)=\log r-(r^2-1)/2$ from \eqref{e.def.Parisi-radius}. Stirling's formula in \eqref{e.radial-density} gives
\begin{equation}
\lim_{N\to\infty}\sup_{r\in[r_0,R]}\biggl|\frac1N\log\rho_M(r)-\al\cR(r)\biggr|=0.
\label{e.radial-density-limit}
\end{equation}
For each fixed $r$, the random variable $B^{(k)}_{N,M}(\be r^2,y)$ is $r\sqrt{\be M}/N$-Lipschitz in $G$, and the balancing argument in Step~1 of the proof of Proposition~\ref{p.shell} shows that its conditional expectation given the $\eta_i^j$ differs from $\bar B^{(k)}_{N,M}(\be r^2,y)$ by a quantity that tends to zero in probability. By Gaussian concentration, $B^{(k)}_{N,M}(\be r^2,y)-\bar B^{(k)}_{N,M}(\be r^2,y)$ therefore tends to zero in probability as $N\to\infty$. To pass from a finite net to the whole interval, we use
\begin{equation*}
|B^{(k)}_{N,M}(\be r^2,y)-B^{(k)}_{N,M}(\be r'^2,y)|\le\frac{\sqrt{\be M}}{N}\norm G_{\mathrm{op}}|r-r'|.
\end{equation*}
The prefactor is bounded in probability because $\norm G_{\mathrm{op}}/\sqrt N$ has bounded moments. Proposition~\ref{p.shell} and a finite-net argument therefore give
\begin{equation}
\lim_{N\to\infty}\sup_{r\in[r_0,R]}\bigl|B^{(k)}_{N,M}(\be r^2,y)-B^{(k)}_\al(\be r^2,y)\bigr|=0 \qquad\text{in probability}.
\label{e.radial-shell-uniform}
\end{equation}
Applying the Laplace principle to \eqref{e.radial-restricted-Z}, using \eqref{e.radial-density-limit} and \eqref{e.radial-shell-uniform}, yields
\begin{equation}
\lim_{N\to\infty}\frac1N\log Z_N([r_0,R],y)=\max_{r\in[r_0,R]}\Ll\{\al\cR(r)+B^{(k)}_\al(\be r^2,y)\Rr\} \qquad\text{in probability}.
\label{e.radial-compact-limit}
\end{equation}

\emph{Step 2: radii near zero and infinity.} We prove that removing the cutoffs has vanishing cost,
\begin{equation}
\lim_{r_0\downarrow0}\lim_{R\to\infty}\limsup_{N\to\infty}\frac1N\E\log\frac{Z_N((0,\infty),y)}{Z_N([r_0,R],y)}=0.
\label{e.radial-cutoff-error}
\end{equation}
Summing over $\si$ produces a product of hyperbolic cosines, so for every realization,
\begin{equation}
\label{e.shell-norm-bound}
0 \le B^{(k)}_{N,M}(\be r^2,y) \le \frac{r \sqrt{\be M}}{N} \norm G_{\mathrm{op}} + \sqrt{\be k} |y|.
\end{equation}
Fix $0<a<b<\infty$ such that $a\le \frac{M}{N}\le b$ for large $N$, let $\Lambda<\infty$, and first suppose that $\norm G_{\mathrm{op}}\le\Lambda\sqrt N$. By Stirling's formula, the prefactor in \eqref{e.radial-density} is at most $e^{CN}$. For $r\ge1$, the logarithm of the integrand in \eqref{e.radial-restricted-Z} is thus at most $N(C+\Lambda\sqrt{\be b}\,r+\sqrt{\be k}|y|+b\log r-ar^2/2)$ by \eqref{e.shell-norm-bound}. For $r\le r_0$, we use \eqref{e.shell-norm-bound} and $\int_0^{r_0}\rho_M(r)\,\d r\le e^{CN}r_0^M$. So, for $N$ large and $R$ sufficiently large, we have
\begin{align}
Z_N([R,\infty),y)\le e^{-cNR^2} \quad\text{and}\quad
Z_N((0,r_0],y)\le\exp\Ll(N\bigl[C+\Lambda\sqrt{\be b}\,r_0+\sqrt{\be k}|y|+a\log r_0\bigr]\Rr),
\label{e.radial-tail-bounds}
\end{align}
where $C,c>0$. Since $B^{(k)}_{N,M}\ge0$ and $\lim_{N\to\infty}\int_{r_0}^R\rho_M(r)\,\d r=1$, we also have $Z_N([r_0,R],y)\ge1/2$ for all large $N$. We choose $r_0$ sufficiently small and $R$ sufficiently large that both upper bounds in \eqref{e.radial-tail-bounds} are at most $e^{-N}$. Splitting the full radial integral into the three intervals and using $\log(1+x)\le x$, we obtain
\begin{equation*}
0\le\frac1N\log\frac{Z_N((0,\infty),y)}{Z_N([r_0,R],y)}\le\frac{4e^{-N}}{N}\qquad \text{on $\{\norm G_{\mathrm{op}}\le\Lambda\sqrt N\}$}.
\end{equation*}

For $\Lambda$ fixed and large enough, the event $\{\norm G_{\mathrm{op}}>\Lambda\sqrt N\}$ has exponentially small probability. The full partition function satisfies
\begin{equation}
0\le\frac1N\log Z_N((0,\infty),y)\le\frac{\be}{2N}\norm G_{\mathrm{op}}^2+\sqrt{\be k}|y|.
\label{e.radial-free-energy-bound}
\end{equation}
Since $\norm G_{\mathrm{op}}/\sqrt N$ has a fourth moment bounded uniformly in $N$, the right-hand side is bounded in $L^2$, so the free energies $N^{-1}\log Z_N((0,\infty),y)$ are uniformly integrable. For fixed $r_0$ and $R$, the truncated free energies $N^{-1}\log Z_N([r_0,R],y)$ lie between $-\log2/N$ and the full ones for all large $N$, so they are uniformly integrable as well, and the quantity $N^{-1}\log[Z_N((0,\infty),y)/Z_N([r_0,R],y)]$ is at most the right-hand side of \eqref{e.radial-free-energy-bound} plus $\log2/N$. Its expectation on $\{\norm G_{\mathrm{op}}>\Lambda\sqrt N\}$ therefore tends to zero, which together with the previous display proves \eqref{e.radial-cutoff-error}. The uniform integrability also permits taking expectations in \eqref{e.radial-compact-limit}.

\emph{Step 3: identification and attainment.} By \eqref{e.radial-compact-limit}, the uniform integrability from Step~2, and \eqref{e.radial-cutoff-error}, we obtain
\begin{equation*}
\lim_{N \to \infty} \bar F^{(k)}_{N,M_N}(y) = \sup_{r > 0} \Ll\{ \al \cR(r) + B^{(k)}_\al(\be r^2,y) \Rr\}.
\end{equation*}
Inserting \eqref{e.shell-limit} gives \eqref{e.hybrid-formula}. Taking expectations in \eqref{e.shell-norm-bound} and then letting $N\to\infty$ gives $0\le B^{(k)}_\al(\be r^2,y)\le C_{\al,\be}r+\sqrt{\be k}|y|$. Thus
\begin{equation*}
\al\cR(r)+B^{(k)}_\al(\be r^2,y)\le\al\Ll(\log r-\frac{r^2-1}{2}\Rr)+C_{\al,\be}r+\sqrt{\be k}|y|.
\end{equation*}
The right-hand side tends to $-\infty$ as $r\downarrow0$ or $r\to\infty$, uniformly for $y$ in a bounded set, whereas the value of the variational objective at $r=1$ is nonnegative. Continuity therefore implies that the supremum over $r$ is attained and that all maximizing radii lie in a common compact interval for $y$ in a bounded set. Proposition~\ref{p.shell} gives the uniqueness and the support of the inner minimizer. Finally, \eqref{e.hybrid-field-Lipschitz} for $\bar F^{(k)}_{N,M}$ upgrades the convergence to local uniformity in $y$.
\end{proof}

\subsection{Comparison of the retained fields}
\label{s.field-comparison}

We finally show that the free energy with $k$ retained fields $y \in \R^k$ is at most that with a single field of strength $|y|$. We recall that $F^{(k)}(y)$ and $A_k(\al,\be)$ are defined in \eqref{e.hybrid-formula} and \eqref{e.def.Ak} respectively.

\begin{proposition}[Several fields are dominated by one field]
\label{p.field-comparison}
For every $k \in \N$, $u \ge 0$ and $y \in \R^k$,
\begin{equation}
\label{e.field-comparison}
B^{(k)}_\al(u,y) \le B^{(1)}_\al(u,|y|) 
\qquad \text{and} \qquad 
F^{(k)}(y) \le F^{(1)}(|y|).
\end{equation}
In particular, for every $k \in \N$, we have $A_k(\al,\be) = A_1(\al,\be)$.
\end{proposition}
\begin{proof}
The random field $\cH_y$ in \eqref{e.def.field-law} satisfies $\E \cH_y^2 = \be |y|^2$. For every $a \ge 0$ and $\boldsymbol{\zeta} \in \cP_{<1}$, the inequality \eqref{e.Jensen-field} of Lemma~\ref{l.squared-field} gives 
\begin{equation*}
\E \, \msc P_a(\boldsymbol{\zeta},\cH_y) \le \msc P_a(\boldsymbol{\zeta},\sqrt \be |y|) = \E \, \msc P_a(\boldsymbol{\zeta},\cH_{|y|}),    
\end{equation*}
where the last equality uses the evenness of $\msc P_a(\boldsymbol{\zeta},\cdot)$ and $\cH_{|y|} = \sqrt \be |y| \ep_1$. We add $\al \cE(\boldsymbol{\zeta})$ to both sides and take the infimum over $\boldsymbol{\zeta}$ in \eqref{e.shell-limit}, then the supremum over $r$ in \eqref{e.hybrid-formula}. This proves \eqref{e.field-comparison}. Since the quadratic cost $|y|^2/2$ in \eqref{e.def.Ak} is unchanged by the replacement of $y$ by $|y|$, we deduce that $A_k(\al,\be) \le A_1(\al,\be)$. Conversely, for every $s \in \R$, we have $\cH_{(s,0,\ldots,0)} = \cH_s$, so $F^{(k)}((s,0,\ldots,0)) = F^{(1)}(s)$, and $A_k(\al,\be) \ge A_1(\al,\be)$.
\end{proof}

%
%
%
%
%
%
\section{Comparison of Bernoulli and Gaussian disorder under cutoffs}
\label{s.comparison}

In this section, we fix $\be > 0$, and we compare the free energies of the bipartite system with Bernoulli and with Gaussian disorder, when the auxiliary Gaussian variables are restricted to a bounded domain. We also show that the cutoffs defining this domain can be removed at a small cost for Gaussian disorder, and we estimate the effect of adding Gaussian columns to the disorder matrix. For an $N \times M$ matrix $A$, a weight $w : \{-1,1\}^N \to [0,\infty)$ with $\int w \, \d P_N > 0$, and a deterministic Borel set $D \subset \R^M$ of positive Gaussian measure, we define
\begin{equation}
\label{e.def.restricted}
\cF_{N,M}(A,w;D) := \frac 1 N \log \iint_{\{-1,1\}^N \times D} w(\si) \, e^{\sqrt{\be/N} \, \si^\top A z} \, \d P_N(\si) \, \d \boldsymbol{\gamma}_M(z).
\end{equation}
We omit $D$ from the notation when $D = \R^M$. For random weights, we assume that $\E |\log \int w \, \d P_N| < \infty$. For $L, R, \de > 0$, we define the ball, the cube, and the domain
\begin{equation}
\label{e.def.cutoffs}
\Ball_{N,M}(L) := \{ z \in \R^M : |z| \le L \sqrt N \}, \quad 
\Cube_M(R) := [-R,R]^M, \quad 
\Dom_{N,M}(\de,L) := \Ball_{N,M}(L) \cap \Cube_M(\de \sqrt N),
\end{equation}
and we also allow $R = \infty$ in \eqref{e.def.cutoffs}, with $\Cube_M(\infty) := \R^M$.

We will use two standard estimates. For an $N \times M$ matrix $A$ with independent entries that are either uniform on $\{-1,1\}$ or standard Gaussian, the bound on the norm of subgaussian matrices \cite[Theorem~5.39]{vershynin2012introduction} gives
\begin{equation}
\label{e.matrix-norm}
\P \bigl( \norm A_{\mathrm{op}} > C (\sqrt N + \sqrt M) + t \bigr) \le 2 e^{-c t^2}, \qquad t \ge 0,
\end{equation}
and integrating this tail bound gives $\E \norm A_{\mathrm{op}}^r \le C_r (\sqrt N + \sqrt M)^r$ for every $r \ge 1$. For a convex function $U$ of independent random variables uniform on $\{-1,1\}$, which is globally $\ell$-Lipschitz for the Euclidean norm, Talagrand's concentration inequality for convex functions gives
\begin{equation}
\label{e.convex-concentration}
\P(|U - \E U| > t) \le C e^{-c t^2/\ell^2}, \qquad t \ge 0.
\end{equation}
This is the consequence of the convex distance inequality given in \cite[Section~4]{talagrand1995concentration}; replacing the median by the mean only changes the numerical constants. Integrating a union bound, we obtain, for every finite family $(U_j)_{j \le J}$ of such functions with a common Lipschitz constant $\ell$, and without any independence assumption between them,
\begin{equation}
\label{e.maximum-concentration}
\E \max_{j \le J} U_j \le \max_{j \le J} \E U_j + C \ell \sqrt{\log(2J)}.
\end{equation}

Our first lemma compares the two kinds of disorder on the domain $\Dom_{N,M}(\de,L)$. 
\begin{lemma}[Comparison of the two kinds of disorder]
\label{l.restricted-comparison}
Let $X$ and $G$ be $N \times M$ matrices with independent entries, uniform on $\{-1,1\}$ and standard Gaussian respectively, and let $w$ be a random weight independent of $X$ and $G$. For every $\de \in (0,1]$, $L > 0$, and every deterministic Borel set $D \subset \Dom_{N,M}(\de,L)$ of positive Gaussian measure,
\begin{equation}
\label{e.restricted-comparison}
|\E \cF_{N,M}(X,w;D) - \E \cF_{N,M}(G,w;D)| \le C_\be L^2 \de^2.
\end{equation}
\end{lemma}
\begin{proof}
\emph{Step 1: fourth derivatives.}
The estimates below are uniform in $w$, so we fix a realization of~$w$ and write $f(A) := \cF_{N,M}(A,w;D)$. Until the final averaging over $w$, the expectations are with respect to $X$ and $G$. Let $\langle \cdot \rangle_A$ denote the probability measure on $\{-1,1\}^N \times D$ proportional to the integrand in \eqref{e.def.restricted}. The derivatives of $f$ with respect to the entry $A_{i\mu}$ are the cumulants of $\sqrt{\be/N} \, \si_i z_\mu$ under $\langle \cdot \rangle_A$, divided by $N$. In particular, with $\kappa_{4,A}$ the fourth cumulant,
\begin{equation}
\label{e.fourth-derivative}
\dr^4_{i\mu} f(A) = \frac{\be^2}{N^3} \kappa_{4,A}(\si_i z_\mu)
\qquad \text{and} \qquad
|\dr^4_{i\mu} f(A)| \le \frac{C \be^2}{N^3} \langle z_\mu^4 \rangle_A,
\end{equation}
where the inequality follows by writing the fourth cumulant in terms of moments, and by bounding every product of moments by the fourth absolute moment using H\"older's inequality.

\emph{Step 2: interpolation with three matching moments.} We take $X$ and $G$ independent. Let $A_\theta := \sqrt \theta \, G + \sqrt{1-\theta} \, X$ and $h(\theta) := \E f(A_\theta)$ for $\theta \in [0,1]$. For $\theta \in (0,1)$, the Gaussian integration by parts gives
\begin{equation}
\label{e.interpolation-derivative}
h'(\theta) = \frac 1 2 \sum_{i,\mu} \E \, \dr^2_{i\mu} f(A_\theta) - \frac{1}{2\sqrt{1-\theta}} \sum_{i,\mu} \E \bigl[ X_{i\mu} \, \dr_{i\mu} f(A_\theta) \bigr].
\end{equation}
For a random variable $\ep$ uniform on $\{-1,1\}$ and a smooth function $v$, an integration gives
\begin{equation}
\label{e.Bernoulli-IBP}
\E_\ep [\ep \, v(b + r \ep)] - r \E_\ep [v'(b + r\ep)] = - \frac 1 4 \int_{-r}^r (r^2 - u^2) \, v'''(b+u) \, \d u.
\end{equation}
Fix $G$ and all entries of $X$ except $X_{i\mu}$. We apply \eqref{e.Bernoulli-IBP} with $b = \sqrt\theta \, G_{i\mu}$, $\ep = X_{i\mu}$, $r=\sqrt{1-\theta}$ and $v(x)=\dr_{i\mu} f(A_\theta)\big|_{(A_\theta)_{i\mu}=x}$. The terms containing $\dr^2_{i\mu} f$ cancel in \eqref{e.interpolation-derivative}. The remaining fourth derivatives are evaluated at $x = \sqrt\theta \, G_{i\mu} + u$, with $|u| \le r$. 

\emph{Step 3: uniform control of the intermediate Gibbs averages.} At each matrix appearing in the integral in Step 2, the exponent in \eqref{e.def.restricted} differs from that at $A_\theta$ by
\begin{equation*}
\sqrt{\be/N} \, \si_i z_\mu (u - r X_{i\mu}),
\qquad
\bigl|\sqrt{\be/N} \, \si_i z_\mu (u - r X_{i\mu})\bigr| \le 2r \sqrt\be \, \de \le 2\sqrt\be,
\end{equation*}
since $|u| \le r \le 1$ and $|z_\mu| \le \de\sqrt N$ on $D$. 
The Gibbs averages of nonnegative functions at the two matrices therefore differ by a factor at most $e^{4 \sqrt \be}$. Using \eqref{e.fourth-derivative} and $\int_{-r}^r (r^2 - u^2) \, \d u = 4 r^3/3$, we obtain
\begin{equation*}
|h'(\theta)| \le \frac{C_\be (1-\theta)}{N^2} \, \E \Ll\langle \sum_{\mu=1}^M z_\mu^4 \Rr\rangle_{A_\theta}.
\end{equation*}
On $D$, we have $\sum_\mu z_\mu^4 \le \max_\mu z_\mu^2 \, |z|^2 \le L^2 \de^2 N^2$. Integrating over $\theta \in [0,1]$ and then averaging over~$w$ proves \eqref{e.restricted-comparison}. 
\end{proof}

We next show that the restriction to the ball $\Ball_{N,M}(L)$ changes the free energy by an exponentially small amount, for either kind of disorder, provided that $L$ is large enough.

\begin{lemma}[Truncation of the norm]
\label{l.norm-truncation}
Let $b < \infty$, and let $A$ be an $N \times M$ matrix with $M \le bN$ and independent entries, either uniform on $\{-1,1\}$ or standard Gaussian. There exist $L_0 = L_0(b,\be) \ge 1$ and $\Lambda = \Lambda(b) < \infty$ satisfying $\P(\norm A_{\mathrm{op}} > \Lambda \sqrt N) \le C e^{-cN}$ and such that, on the event $\{\norm A_{\mathrm{op}} \le \Lambda \sqrt N\}$, we have for every $L \ge L_0$, $R \in (0,\infty]$, and every weight $w$ that
\begin{equation}
\label{e.norm-truncation}
0 \le \cF_{N,M}(A,w;\Cube_M(R)) - \cF_{N,M}(A,w;\Cube_M(R) \cap \Ball_{N,M}(L)) \le C e^{-cN}.
\end{equation}
Moreover, the difference in \eqref{e.norm-truncation} is always nonnegative, and its expectation is at most $C_{b,\be} \, e^{-cN}$ for every $L \ge L_0$, $R \in (0,\infty]$, and every random weight $w$, which may depend on $A$.
\end{lemma}
\begin{proof}
We fix $A$ and $w$, and let $\langle\cdot\rangle$ denote the Gibbs measure on $\{-1,1\}^N\times\R^M$ proportional to the integrand in \eqref{e.def.restricted} with $D=\R^M$. Let $Y\in\R^M$ be a standard Gaussian vector independent of the Gibbs sample. Completing the square shows that, conditionally on $\si$, the variable $z$ has the law of $Y+u(\si)$, where
\begin{equation}
\label{e.auxiliary-mean}
u(\si) := \sqrt{\be/N} \, A^\top \si 
\qquad \text{satisfies} \qquad 
|u(\si)| \le \sqrt \be \norm A_{\mathrm{op}}.
\end{equation}
The difference in \eqref{e.norm-truncation} can therefore be described as
\begin{equation}
\text{\ diff.\ in \eqref{e.norm-truncation}\ }=\frac{1}{N}\log\frac{\langle\1_{\{z\in\Cube_M(R)\}}\rangle}{\langle\1_{\{z\in\Cube_M(R)\cap\Ball_{N,M}(L)\}}\rangle}=\frac{1}{N}\log\frac{\langle\P_Y(Y+u(\si)\in\Cube_M(R))\rangle}{\langle\P_Y(Y+u(\si)\in\Cube_M(R)\cap\Ball_{N,M}(L))\rangle}.
\label{e.diff=fracP}
\end{equation}
In particular, the difference is nonnegative. 

For fixed $u\in\R^M$, conditioning on $Y+u\in\Cube_M(R)$ preserves independence of the coordinates and conditions each $|Y_\mu+u_\mu|$ to be at most $R$. For every $s\ge0$, we have
\begin{equation*}
\P_Y\bigl(|Y_\mu+u_\mu|\le s \ \big| \  |Y_\mu+u_\mu|\le R\bigr)=\frac{\P_Y(|Y_\mu+u_\mu|\le\min\{s,R\})}{\P_Y(|Y_\mu+u_\mu|\le R)}\ge\P_Y(|Y_\mu+u_\mu|\le s).
\end{equation*}
Applying this inequality successively to the coordinates of the event $\{\sum_\mu|Y_\mu+u_\mu|^2\le L^2N\}$ gives
\begin{equation}
\P_Y(Y+u\in\Cube_M(R)\cap\Ball_{N,M}(L))\ge\P_Y(Y+u\in\Cube_M(R))\,\P_Y(Y+u\in\Ball_{N,M}(L)).
\label{e.cube-ball-probability}
\end{equation}

Let $\Lambda:=C(1+\sqrt b)+1$, where $C$ is the constant in \eqref{e.matrix-norm}, and let $L_0:=\sqrt\be\Lambda+\sqrt b+2$. Since $M\le bN$, the bound \eqref{e.matrix-norm} with $t=\sqrt N$ gives $\P(\norm A_{\mathrm{op}}>\Lambda\sqrt N)\le2e^{-cN}$. 
Since $y\mapsto|y|$ is $1$-Lipschitz and $\E_Y|Y|\le\sqrt M\le\sqrt{bN}$, Gaussian concentration gives $\P_Y(|Y|>(\sqrt b+2)\sqrt N)\le e^{-2N}$. On $\{\norm A_{\mathrm{op}}\le\Lambda\sqrt N\}$, for every $L\ge L_0$, the bound \eqref{e.auxiliary-mean} therefore implies
\begin{equation}
\label{e.norm-probability}
\inf_\si \P_Y \bigl( |Y + u(\si)| \le L \sqrt N \bigr) \ge 1 - e^{-2N}.
\end{equation}
By \eqref{e.cube-ball-probability} and \eqref{e.norm-probability}, the denominator in \eqref{e.diff=fracP} is at least $1-e^{-2N}$ times its numerator. Hence the difference in \eqref{e.norm-truncation} is at most $-N^{-1}\log(1-e^{-2N})\le2e^{-2N}$. This proves \eqref{e.norm-truncation} on $\{\norm A_{\mathrm{op}}\le\Lambda\sqrt N\}$.

For arbitrary $A$, the symmetry of the ball and Jensen's inequality give, for every $u\in\R^M$,
\begin{equation*}
\P_Y(Y + u \in \Ball_{N,M}(L)) = e^{-|u|^2/2} \int_{\Ball_{N,M}(L)} e^{u \cdot z} \, \d \boldsymbol{\gamma}_M(z) \ge e^{-|u|^2/2} \boldsymbol{\gamma}_M(\Ball_{N,M}(L)),
\end{equation*}
and the last Gaussian mass is at least $1-e^{-2N}\ge1/2$ for $L\ge L_0$. Using \eqref{e.auxiliary-mean} in \eqref{e.cube-ball-probability}, we obtain
\begin{equation*}
\P_Y(Y+u(\si)\in\Cube_M(R)\cap\Ball_{N,M}(L))\ge\frac12 e^{-\be\norm A_{\mathrm{op}}^2/2}\,\P_Y(Y+u(\si)\in\Cube_M(R)).
\end{equation*}
Averaging under $\langle\cdot\rangle$ and substituting into \eqref{e.diff=fracP} yields the deterministic bound
\begin{equation}
0\le\text{\ diff.\ in \eqref{e.norm-truncation}\ }\le\frac{\be}{2N}\norm A_{\mathrm{op}}^2+\frac{\log2}{N}.
\label{e.norm-truncation-deterministic}
\end{equation}
Using \eqref{e.norm-truncation} on the event $\{\norm A_{\mathrm{op}}\le\Lambda\sqrt N\}$ and \eqref{e.norm-truncation-deterministic} on its complement, followed by an application of the Cauchy--Schwarz inequality and the moment bound $\E\norm A_{\mathrm{op}}^4\le C_bN^2$ from \eqref{e.matrix-norm}, we obtain
\begin{align*}
&\E\bigl[\text{\ diff.\ in \eqref{e.norm-truncation}\ }\bigr]\le2e^{-2N}+\E\left[\left(\frac{\be}{2N}\norm A_{\mathrm{op}}^2+\frac{\log2}{N}\right)\1_{\{\norm A_{\mathrm{op}}>\Lambda\sqrt N\}}\right]\le C_{b,\be}\,e^{-cN}.
\end{align*}
All bounds hold for every $w$, so no independence between $A$ and $w$ is needed.
\end{proof}

As emphasized in the introduction, for Bernoulli patterns we will show absence of universality in the disorder for certain regimes of $\al$ and $\be$; in other words, for sufficiently small $\de > 0$, the asymptotic behaviors of $\cF_{N,M}(A, 1)$ and $\cF_{N,M}(A,1; \Cube_M(\delta \sqrt{N}))$ will differ in these cases. We show below that this does not occur in the case of Gaussian patterns. 

\begin{lemma}[Truncation of the coordinates for Gaussian disorder]
\label{l.Gaussian-coordinate}
For every $0 < a \le b < \infty$, $aN \le M \le bN$, and $R \ge 2$, uniformly over the weights $w$ independent of $G$,
\begin{equation}
\label{e.Gaussian-coordinate}
0 \le \E \cF_{N,M}(G,w) - \E \cF_{N,M}(G,w;\Cube_M(R)) \le \frac{C_{a,b,\be}}{R^2}.
\end{equation}
\end{lemma}
\begin{proof}
We condition on $w$, and let $\langle \cdot \rangle$ denote the Gibbs measure on $\{-1,1\}^N \times \R^M$ proportional to the integrand in \eqref{e.def.restricted}, with $A = G$ and $D = \R^M$. The ratio of the partition function restricted to the cube to the full one is $\langle \1_{\{z \in \Cube_M(R)\}} \rangle$. As in the proof of Lemma~\ref{l.norm-truncation}, given $\si$, the variable $z$ has the law $\mcl N(u, \Id_M)$ with $u = \sqrt{\be/N} \, G^\top \si$, so this ratio equals $\langle \prod_\mu \P(|Y + u_\mu| \le R) \rangle$, where $Y$ is a scalar standard Gaussian. Let $\ell_R(v) := - \log \P(|Y + v| \le R)$. Jensen's inequality gives
\begin{equation}
\label{e.coordinate-Jensen}
\cF_{N,M}(G,w) - \cF_{N,M}(G,w;\Cube_M(R)) \le \frac 1 N \sum_{\mu=1}^M \langle \ell_R(u_\mu) \rangle.
\end{equation}
For $|v| \le R/2$ and $R \ge 2$, the Gaussian tail gives $\ell_R(v) \le C e^{-R^2/8}$, after increasing the constant to cover a bounded range of $R$. For arbitrary $v$, integrating the shifted Gaussian density over $[-1,1]$ gives $\ell_R(v) \le C(1 + v^2)$. Hence
\begin{equation}
\label{e.ell-tail}
\ell_R(v) \le C e^{-R^2/8} + C (1 + v^2) \1_{\{|v| > R/2\}}.
\end{equation}

For every orthogonal matrix $O \in O(M)$, the matrix $GO$ has the same law as $G$. Since $\boldsymbol{\gamma}_M$ is rotationally invariant, the change of variables $z \mapsto O^\top z$ shows that the law of $\si$ under $\langle \cdot \rangle$ is unchanged when $G$ is replaced by $GO$, while the auxiliary mean becomes $O^\top u$. The law of $u$ averaged over $G$ and over the Gibbs measure is therefore rotationally invariant. The fourth moment of a coordinate of a uniformly distributed point on a sphere gives
\begin{equation}
\label{e.rotational-fourth}
\E \langle u_\mu^4 \rangle = \frac{3}{M(M+2)} \E \langle |u|^4 \rangle \le \frac{3 \be^2}{M(M+2)} \E \norm G_{\mathrm{op}}^4 \le C_{a,b,\be},
\end{equation}
where we used \eqref{e.auxiliary-mean} and \eqref{e.matrix-norm}; the identity for the sphere also holds when $M = 1$. On $\{|v| > R/2\}$ with $R \ge 2$, we have $1 + v^2 \le 8 v^4/R^2$. The estimates \eqref{e.ell-tail} and \eqref{e.rotational-fourth} therefore bound the expected Gibbs average in \eqref{e.coordinate-Jensen} by $C/R^2$ for each $\mu$. Summing over $\mu$ and using $M/N \le b$ proves \eqref{e.Gaussian-coordinate}.
\end{proof}

Combining Lemmas~\ref{l.norm-truncation} and~\ref{l.Gaussian-coordinate}, we obtain, for $aN \le M \le bN$, $\de \sqrt N \ge 2$, and $L \ge L_0(b,\be)$,
\begin{equation}
\label{e.Gaussian-both-cutoffs}
0 \le \E \cF_{N,M}(G,w) - \E \cF_{N,M}(G,w;\Dom_{N,M}(\de,L)) \le \frac{C_{a,b,\be}}{\de^2 N} + C_{b,\be} e^{-cN},
\end{equation}
for every weight $w$ independent of $G$. 

In the next section, we replace all but a bounded number of the Bernoulli patterns by Gaussian ones. The resulting models have slightly fewer than $M_N$ Gaussian columns, and the following lemma allows us to adjust this number: the expected free energy $\bar F^{(k)}_{N,M}$ in \eqref{e.def.hybrid} is nondecreasing in the number $M$ of Gaussian columns, and each additional column increases it by at most a constant divided by~$N$.

\begin{lemma}[Adding Gaussian columns]
\label{l.Gaussian-columns}
For every $0 < a \le b < \infty$ and $0 \le M \le M'$ with $aN \le M' \le bN$, uniformly in $k$ and $y$,
\begin{equation}
\label{e.Gaussian-columns}
0 \le \bar F^{(k)}_{N,M'}(y) - \bar F^{(k)}_{N,M}(y) \le C_{a,b,\be} \frac{M' - M}{N}.
\end{equation}
\end{lemma}
\begin{proof}
We multiply the last $M' - M$ nonnegative quadratic terms in the exponent of \eqref{e.def.hybrid} by $s \in [0,1]$. The expected free energy is nondecreasing and convex in $s$, so its increment over $[0,1]$ is bounded by its derivative at $s = 1$. At this endpoint, the Gaussian columns are exchangeable and independent of the fields, so the derivative is at most
\begin{equation*}
\frac{\be (M' - M)}{2 N^2 M'} \E \langle |G^\top \si|^2 \rangle \le \frac{\be (M' - M)}{2 N M'} \E \norm G_{\mathrm{op}}^2 \le C_{a,b,\be} \frac{M' - M}{N},
\end{equation*}
by \eqref{e.matrix-norm}. This proves \eqref{e.Gaussian-columns}.
\end{proof}

%
%
%
%
%
%
\section{Retained patterns and the Bernoulli limit}
\label{s.finite-pattern}

In this section, we prove that the Bernoulli free energy is close, uniformly in $N$, to the free energy of the model with Gaussian patterns and finitely many retained Bernoulli patterns, and we deduce Theorem~\ref{t.sup-Ak}, Proposition~\ref{p.structure}, as well as \eqref{e.def.FB}, \eqref{e.main-formula} and \eqref{e.one-pattern}. The remaining assertions of Theorem~\ref{t.hopfield} and Corollary~\ref{c.one-pattern} are proved in Section~\ref{s.consequences}. We recall the free energy $\bar F^{(k)}_{N,M}$ from \eqref{e.def.hybrid} and the limiting optimization $A_k$ from \eqref{e.def.Ak}. For $L > 0$ and $k \in \N \cup \{0\}$, we define
\begin{equation}
\label{e.def.ANk}
A_{N,k}(L) := \sup_{|y| \le L} \Ll\{ \bar F^{(k)}_{N,M_N}(y) - \frac{|y|^2}{2} \Rr\}.
\end{equation}

\subsection{A uniform finite-size estimate}
\label{s.finite-size}

The main step is the following estimate, which compares $\bar F_N(\be)$ (defined in \eqref{e.def.free-energy}) with $A_{N,K}(L)$ for a number $K$ of retained patterns that can remain bounded. 

\begin{proposition}[Approximation by finitely many retained patterns]
\label{p.finite-pattern}
Let $0 < a \le b < \infty$ and $\be > 0$. There exists $L_0 < \infty$ such that, for every $L \ge L_0$, $N \ge 2$, and $\de \in (0,1]$, with $K := \lceil L^2/\de^2 \rceil$, the assumptions
\begin{equation}
\label{e.finite-pattern-assumptions}
aN \le M_N \le bN, \qquad K \le M_N/2, \qquad \text{and} \qquad \de \sqrt N \ge 2
\end{equation}
imply
\begin{equation}
\label{e.finite-pattern-bound}
|\bar F_N(\be) - A_{N,K}(L)| \le C_\be L^2 \de^2 + C_{a,b,\be,L} \Ll( \sqrt{\frac{K \log N}{N}} + \frac{K \log N}{N} \Rr).
\end{equation}
\end{proposition}
\begin{proof}
\emph{Step 1: truncation of the norm.} Let $L_0 = L_0(b,\be) \ge 1$ be given by Lemma~\ref{l.norm-truncation}, and let $L \ge L_0$. By \eqref{e.HS} and \eqref{e.def.restricted}, we have $\bar F_N(\be) = \E \cF_{N,M_N}(X_N,1)$. Let $Z_N$ and $Z_N^L$ denote the partition functions in \eqref{e.def.restricted} for $A = X_N$, $w = 1$, and $D = \R^{M_N}$ and $D = \Ball_{N,M_N}(L)$ respectively, so that $\bar F_N(\be) = \frac 1 N \E \log Z_N$. By Lemma~\ref{l.norm-truncation}, applied with $R=\infty$, we have
\begin{equation}
\label{e.Bernoulli-norm-loss}
0 \le \bar F_N(\be) - \frac 1 N \E \log Z_N^L \le C_{b,\be} e^{-cN}.
\end{equation}

\emph{Step 2: free energy with retained patterns.} We fix a deterministic set $S = \{\mu_1 < \cdots < \mu_s\} \subset \{1,\ldots,M_N\}$ with $s \le K$, and write $S^{\complement} := \{1,\ldots,M_N\}\setminus S$. Let $X_S$ and $X_{S^{\complement}}$ denote the matrices of the columns of $X_N$ indexed by these sets. For $y \in \R^s$, we define the weight and the restricted free energy, using the domain $\Dom_{N,M_N-s}(\de,L)$ from \eqref{e.def.cutoffs}, by 
\begin{equation}
\label{e.def.retained}
w_{S,y}(\si) := \exp \Ll( \sqrt \be \sum_{i=1}^N \sum_{j=1}^s y_j x_i^{\mu_j} \si_i \Rr)
\quad \text{and} \quad
\mcl F_S(y) := \cF_{N,M_N - s}(X_{S^{\complement}}, w_{S,y}; \Dom_{N,M_N-s}(\de,L)).
\end{equation}
These definitions arise by fixing $z_S = \sqrt N y$ in the partition function $Z_N$ of $\cF_{N,M_N}(X_N,1)$. Indeed, the integrand factors as $e^{\sqrt{\be/N} \, \si^\top X_N z}\big|_{z_S=\sqrt N y} = w_{S,y}(\si) \, e^{\sqrt{\be/N} \, \si^\top X_{S^{\complement}}z_{S^{\complement}}}$.
Under the same change of variables, the Gaussian measure becomes
\begin{equation}
\d\boldsymbol{\gamma}_{M_N}(z) = \Ll(\frac{N}{2\pi}\Rr)^{s/2} e^{-N|y|^2/2} \, \d y \, \d\boldsymbol{\gamma}_{M_N-s}(z_{S^{\complement}}), \qquad z_S=\sqrt N y.
\label{e.retained-Gaussian-measure}
\end{equation}
The weight $w_{S,y}$ depends only on $X_S$, and is therefore independent of $X_{S^{\complement}}$. Let $G$ be an $N \times (M_N - s)$ matrix with independent standard Gaussian entries, independent of $X_N$. Lemma~\ref{l.restricted-comparison} then gives 
\begin{equation}
\label{e.retained-comparison}
|\E \mcl F_S(y) - \E \cF_{N,M_N-s}(G, w_{S,y}; \Dom_{N,M_N-s}(\de,L))| \le C_\be L^2 \de^2.
\end{equation}
Without cutoffs, the Hubbard--Stratonovich identity \eqref{e.HS} identifies this Gaussian model with the free energy defined in \eqref{e.def.hybrid}, namely 
\begin{equation}
\label{e.retained-hybrid}
\E \cF_{N,M_N-s}(G, w_{S,y}) = \bar F^{(s)}_{N,M_N-s}(y).
\end{equation}

For fixed $S$ and $y$, the function $\mcl F_S(y)$ is convex in the entries of $X_N$, being the logarithm of an integral of exponentials of affine functions of these entries. Let $\langle\cdot\rangle$ denote the Gibbs measure associated with \eqref{e.def.retained}. The derivative with respect to $x_i^\mu$ is $\sqrt{\be/N}\,\langle\si_i z_\mu\rangle/N$ for $\mu\in S^{\complement}$, and $\sqrt\be\,y_j\langle\si_i\rangle/N$ for $\mu=\mu_j\in S$. Since $|z_{S^{\complement}}|\le L\sqrt N$ on $\Dom_{N,M_N-s}(\de,L)$, we obtain
\begin{equation}
|\nabla_{X_N}\mcl F_S(y)|^2 \le \frac{\be}{N^2}\langle|z_{S^{\complement}}|^2\rangle + \frac{\be|y|^2}{N} \le \frac{\be(L^2+|y|^2)}{N}.
\label{e.retained-disorder-Lipschitz}
\end{equation}
Thus $\mcl F_S(y)$ is globally $\sqrt{2\be}\,L/\sqrt N$-Lipschitz in $X_N$ when $|y|\le L$. For Bernoulli $X_N$, each derivative in $y$ satisfies
\begin{equation}
\left|\partial_{y_j}\mcl F_S(y)\right| = \left|\frac{\sqrt\be}{N}\sum_{i=1}^N x_i^{\mu_j}\langle\si_i\rangle\right| \le \sqrt\be, \qquad |\nabla_y\mcl F_S(y)|\le\sqrt{\be s}.
\label{e.retained-field-Lipschitz}
\end{equation}

\emph{Step 3: upper bound.} We partition the ball $\Ball_{N,M_N}(L)$ from \eqref{e.def.cutoffs} according to $S(z) := \{\mu : |z_\mu| > \de\sqrt N\}$. Since $K\de^2\ge L^2$, every $z\in\Ball_{N,M_N}(L)$ satisfies 
\begin{equation*}
|S(z)|\de^2 N \le \sum_{\mu\in S(z)}z_\mu^2 \le |z|^2 \le L^2N \le K\de^2N.
\end{equation*}
On the piece where $S(z)=S$, we have $|z_S|\le L\sqrt N$ and $z_{S^{\complement}}\in\Dom_{N,M_N-s}(\de,L)$. Enlarging each piece to this product of sets and then setting $z_S=\sqrt N y$, we obtain 
\begin{align}
Z_N^L &= \sum_{s=0}^K\ \sum_{|S|=s}\int_{\{z\in\Ball_{N,M_N}(L):\,S(z)=S\}}\int e^{\sqrt{\be/N}\,\si^\top X_Nz}\,\d P_N(\si)\,\d\boldsymbol{\gamma}_{M_N}(z)\notag\\
&\le \sum_{s=0}^K\ \sum_{|S|=s}\int_{\{|z_S|\le L\sqrt N, \  z_{S^{\complement}}\in\Dom_{N,M_N-s}(\de,L)\}}\int e^{\sqrt{\be/N}\,\si^\top X_Nz}\,\d P_N(\si)\,\d\boldsymbol{\gamma}_{M_N}(z)\notag\\
&\stackrel{\eqref{e.retained-Gaussian-measure},\,\eqref{e.def.retained}}{=} \sum_{s=0}^K \ \sum_{|S| = s} \Ll( \frac{N}{2\pi} \Rr)^{s/2} \int_{|y| \le L} \exp \Ll( N \Ll[ \mcl F_S(y) - \frac{|y|^2}{2} \Rr] \Rr) \d y,
\label{e.pattern-decomposition}
\end{align}
where for $s=0$ the integral is the evaluation at the empty vector. The assumption $M_N\le bN$ in \eqref{e.finite-pattern-assumptions} bounds the number of subsets by 
\begin{equation*}
\sum_{s=0}^K\binom{M_N}{s} \le (K+1)M_N^K \le (K+1)(bN)^K.
\end{equation*}
The ball of radius $L$ in $\R^s$ has volume at most $(2L)^s$. Bounding the integrand in \eqref{e.pattern-decomposition} by its supremum, and taking the logarithm of the number of subsets, the volumes, and the Gaussian prefactors, we deduce that
\begin{equation}
\label{e.entropy-subsets}
\frac 1 N \log Z_N^L \le \sup_{\substack{|S| \le K \\ |y| \le L}} \Ll( \mcl F_S(y) - \frac{|y|^2}{2} \Rr) + C_{b,L} \frac{K \log N}{N}.
\end{equation}
To bound the expectation of this supremum, we choose, for every $S$ with $|S|=s$, an $N^{-1}$-net of the ball of radius $L$ in $\R^s$ with at most $(3LN)^s$ points. The total number $J$ of pairs $(S,y)$ in these nets satisfies $\log(2J)\le C_{b,L}K\log N$. For fixed $(S,y)$, the function $\mcl F_S(y)-|y|^2/2$ is convex and, by \eqref{e.retained-disorder-Lipschitz}, is $\sqrt{2\be}\,L/\sqrt N$-Lipschitz in the Bernoulli entries of $X_N$. Hence \eqref{e.maximum-concentration} gives
\begin{equation*}
\E \max_{(S,y)} \Ll( \mcl F_S(y) - \frac{|y|^2}{2} \Rr) \le \max_{(S,y)} \Ll( \E \mcl F_S(y) - \frac{|y|^2}{2} \Rr) + C_{b,\be,L} \sqrt{\frac{K \log N}{N}},
\end{equation*}
where both maxima are over the $J$ pairs. By \eqref{e.retained-field-Lipschitz}, the function $\mcl F_S(y)-|y|^2/2$ is $(\sqrt{\be s}+L)$-Lipschitz on $\{|y|\le L\}$. Passing from the net to the full ball on the left therefore costs at most $(\sqrt{\be K}+L)/N\le C_{\be,L}\sqrt{K\log N/N}$, while enlarging the maximum on the right can only increase it. We obtain
\begin{equation}
\label{e.deterministic-subsets}
\E \sup_{\substack{|S| \le K \\ |y| \le L}} \Ll( \mcl F_S(y) - \frac{|y|^2}{2} \Rr) \le \sup_{\substack{|S| \le K \\ |y| \le L}} \Ll( \E \mcl F_S(y) - \frac{|y|^2}{2} \Rr) + C_{b,\be,L} \sqrt{\frac{K \log N}{N}}.
\end{equation}
For deterministic $S$ and $y$, we use \eqref{e.retained-comparison}, enlarge the domain of integration of the Gaussian model to $\R^{M_N-s}$, and then use \eqref{e.retained-hybrid}. Adding $s$ Gaussian columns can only increase the expected free energy by Lemma~\ref{l.Gaussian-columns}, while padding $y$ with $K-s$ zero coordinates leaves it unchanged. Consequently, we have
\begin{equation}
\label{e.upper-pattern-comparison}
\E \mcl F_S(y) \le \bar F^{(s)}_{N,M_N - s}(y) + C_\be L^2 \de^2 \le \bar F^{(K)}_{N,M_N}((y,0)) + C_\be L^2 \de^2.
\end{equation}
Since $|(y,0)|=|y|\le L$, the definition \eqref{e.def.ANk} bounds the right-hand side by $A_{N,K}(L)+|y|^2/2+C_\be L^2\de^2$. 
Combining \eqref{e.Bernoulli-norm-loss}, \eqref{e.entropy-subsets}, \eqref{e.deterministic-subsets} and \eqref{e.upper-pattern-comparison}, we obtain
\begin{equation*}
\bar F_N(\be) \le A_{N,K}(L) + C_\be L^2 \de^2 + C_{b,\be,L} \sqrt{\frac{K \log N}{N}} + C_{b,L} \frac{K \log N}{N} + C_{b,\be} e^{-cN}.
\end{equation*}
Since $e^{-cN}\le C/N\le CK\log N/N$, this proves the upper bound in \eqref{e.finite-pattern-bound}.

\emph{Step 4: lower bound.} Let $S:=\{1,\ldots,K\}$, and fix $y\in\R^K$ with $|y|\le L$. We restrict the integral defining $Z_N$ to $z_{S^{\complement}}\in\Dom_{N,M_N-K}(\de,L)$ and $|z_S/\sqrt N-y|\le N^{-1}$. This restriction is applied to $Z_N$, so there is no constraint on the full norm $|z|$. Positivity of the integrand and the change of variables $z_S=\sqrt N v$ give 
\begin{align}
Z_N &\ge \int_{\{|z_S/\sqrt N-y|\le N^{-1}, \ z_{S^{\complement}}\in\Dom_{N,M_N-K}(\de,L)\}}\int e^{\sqrt{\be/N}\,\si^\top X_Nz}\,\d P_N(\si)\,\d\boldsymbol{\gamma}_{M_N}(z)\notag\\
&\stackrel{\eqref{e.retained-Gaussian-measure},\,\eqref{e.def.retained}}{=} \Ll( \frac{N}{2\pi} \Rr)^{K/2} \int_{|v - y| \le N^{-1}} \exp \Ll( N \Ll[ \mcl F_S(v) - \frac{|v|^2}{2} \Rr] \Rr) \d v.
\label{e.lower-ball-integral}
\end{align}
By \eqref{e.retained-field-Lipschitz} and $|v|^2\le |y|^2+2L/N+N^{-2}$, every $v$ in this ball satisfies
\begin{equation*}
\mcl F_S(v) - \frac{|v|^2}{2} \ge \mcl F_S(y) - \frac{|y|^2}{2} - \frac{\sqrt{\be K} + L}{N} - \frac{1}{2N^2}.
\end{equation*}
The unit ball of $\R^K$ contains the cube $[-K^{-1/2},K^{-1/2}]^K$, so the ball of integration in \eqref{e.lower-ball-integral} has volume at least $(2/\sqrt K)^K N^{-K}$. Since $K\le M_N\le bN$ by \eqref{e.finite-pattern-assumptions}, the logarithm of $(N/(2\pi))^{K/2}$ times this volume is at least $-C_bK\log N$. Taking logarithms and expectations in \eqref{e.lower-ball-integral}, we obtain
\begin{equation}
\label{e.lower-ball}
\bar F_N(\be) \ge \E \mcl F_S(y) - \frac{|y|^2}{2} - C_{b,\be,L} \frac{K \log N}{N}.
\end{equation}
We now bound $\E\mcl F_S(y)$ from below. By \eqref{e.finite-pattern-assumptions}, we have $aN/2\le M_N-K\le bN$ and $\de\sqrt N\ge2$. Since $L\ge L_0$ and $w_{S,y}$ is independent of $G$, the estimate \eqref{e.Gaussian-both-cutoffs}, with $a/2$ in place of $a$, applies to the Gaussian model with $M_N-K$ columns. Together with \eqref{e.retained-comparison} and \eqref{e.retained-hybrid}, it yields
\begin{equation*}
\E \mcl F_S(y) \ge \bar F^{(K)}_{N,M_N - K}(y) - C_\be L^2 \de^2 - \frac{C_{a,b,\be}}{\de^2 N} - C_{b,\be} e^{-cN}.
\end{equation*}
Lemma~\ref{l.Gaussian-columns} gives $\bar F^{(K)}_{N,M_N-K}(y)\ge\bar F^{(K)}_{N,M_N}(y)-C_{a,b,\be}K/N$. Since $K\de^2\ge L^2\ge1$, we have $1/(\de^2N)\le K/N$, and $e^{-cN}\le C/N$. Substituting into \eqref{e.lower-ball}, we conclude that
\begin{equation*}
\bar F_N(\be) \ge \bar F^{(K)}_{N,M_N}(y)-\frac{|y|^2}{2}-C_\be L^2\de^2-C_{a,b,\be,L}\frac{K\log N}{N}.
\end{equation*}
Taking the supremum over $|y|\le L$ and using \eqref{e.def.ANk} proves the lower bound in \eqref{e.finite-pattern-bound}.
\end{proof}

\subsection{Identification of the limit}
\label{s.identification}

\begin{proof}[Proof of Theorem~\ref{t.sup-Ak}]
We choose two numbers $0 < a < b < \infty$ such that $M_N/N \in [a,b]$ for all large $N$. We also fix $L \ge L_0$, where $L_0$ is given by Proposition~\ref{p.finite-pattern}. For every fixed $k$, we have
\begin{equation}
\label{e.ANk-limit}
\lim_{N \to \infty} A_{N,k}(L) = A_k(L) := \sup_{|y| \le L} \Ll\{ F^{(k)}(y) - \frac{|y|^2}{2} \Rr\}.
\end{equation}
Indeed, the definitions give
\begin{equation*}
|A_{N,k}(L) - A_k(L)| \le \sup_{|y| \le L} \Ll|\bar F^{(k)}_{N,M_N}(y) - F^{(k)}(y)\Rr|,
\end{equation*}
and the right-hand side tends to zero as $N \to \infty$ by the local uniform convergence in Theorem~\ref{t.hybrid}.
For a fixed $\de \in (0,1]$, let $K := \lceil L^2/\de^2 \rceil$. The assumptions \eqref{e.finite-pattern-assumptions} hold for $N$ large enough, and \eqref{e.finite-pattern-bound} and \eqref{e.ANk-limit} give
\begin{equation}
\label{e.limsup-AK}
\limsup_{N \to \infty} |\bar F_N(\be) - A_K(L)| \le C_\be L^2 \de^2.
\end{equation}
Since $A_K(L) \le A_K(\al,\be) \le \sup_{k} A_k(\al,\be)$ (see \eqref{e.def.Ak}), we deduce that 
\begin{equation*}
\limsup_{N \to \infty} \bar F_N(\be) \le \sup_{k \in \N} A_k(\al,\be) + C_\be L^2 \de^2.
\end{equation*}
Conversely, let $j \in \N$. Since $\cH_{(y,0)} = \cH_y$, evaluating $F^{(K)}$ at $y \in \R^j$ padded with $K - j$ zero coordinates gives $F^{(j)}(y)$, so $A_K(L) \ge A_j(L)$ as soon as $K \ge j$, which holds for $\de$ small enough. By \eqref{e.limsup-AK},
\begin{equation*}
\liminf_{N \to \infty} \bar F_N(\be) \ge A_j(L) - C_\be L^2 \de^2.
\end{equation*}
Sending $\de \downarrow 0$ with $L$ and $j$ fixed in the last two displays, we obtain
\begin{equation*}
\sup_{k \in \N} A_k(\al,\be) \ge \limsup_{N \to \infty} \bar F_N(\be) \ge \liminf_{N \to \infty} \bar F_N(\be) \ge A_j(L).
\end{equation*}
As $L \to \infty$, the quantity $A_j(L)$ increases to $A_j(\al,\be)$. Taking then the supremum over $j$, we obtain
\begin{equation}
\label{e.Bernoulli-supAk}
\lim_{N \to \infty} \bar F_N(\be) = \sup_{k \in \N} A_k(\al,\be).
\end{equation}
By \eqref{e.matrix-norm}, we have $0 \le \bar F_N(\be) \le \frac{\be}{2N} \E \norm{X_N}_{\mathrm{op}}^2 \le C_{a,b} \be$, so the limit is finite. This proves the existence of the limit \eqref{e.def.FB} and the identity \eqref{e.sup-Ak}. 

Finally, the inequality $F^{(k)}(y) \le F^{(1)}(|y|)$ and the identity $A_k(\al,\be) = A_1(\al,\be)$ are given by Proposition~\ref{p.field-comparison}.
\end{proof}

\begin{proof}[Proof of \eqref{e.def.FB}--\eqref{e.main-formula} and of \eqref{e.one-pattern}]
The existence of the limit \eqref{e.def.FB} and the first identity in \eqref{e.one-pattern} follow from Theorem~\ref{t.sup-Ak}. It remains to prove the second identity in \eqref{e.one-pattern} and the formula \eqref{e.main-formula}. The function $F^{(1)}$ is even (by $\si \mapsto -\si$ in \eqref{e.def.hybrid}) and $\sqrt \be$-Lipschitz (by \eqref{e.hybrid-field-Lipschitz}, which also holds for $\bar F^{(1)}_{N,M}$, and Theorem~\ref{t.hybrid}). For $s \ge \sqrt \be$, we have 
\begin{equation*}
F^{(1)}(s) - \frac{s^2}{2} \le F^{(1)}(\sqrt \be) - \frac \be 2 - \frac{(s - \sqrt \be)^2}{2},
\end{equation*}
so $A_1(\al,\be) = \max_{0 \le s \le \sqrt \be} \{F^{(1)}(s) - s^2/2\}$. Substituting $h = \sqrt \be s$ and using the identity \eqref{e.gaussian-formula}, we obtain \eqref{e.one-pattern}. The formula \eqref{e.main-formula} follows from \eqref{e.one-pattern} by inserting \eqref{e.gaussian-formula} and exchanging the two outer suprema.
\end{proof}

\begin{proof}[Proof of Proposition~\ref{p.structure}]
Let $r > 0$, $h \in \R$ and $a := \al \be r^2/2$. For every $b \in [0,1)$, the continuity of $\boldsymbol{\zeta} \mapsto \msc P_a(\boldsymbol{\zeta},h)$ on $\cP([0,b])$ follows directly from \eqref{e.P-law-continuity} (weak convergence on $[0,b]$ is equivalent to convergence for $W_1$). Together with Lemma~\ref{l.spherical-duality}, this proves the asserted continuity of $\msc P_a(\cdot,h) + \al \cE$.
The strict convexity on $\cP_{<1}$, the uniqueness of the minimizer, and the restriction of its support to $[0,\be r^2/(1+\be r^2)]$ are given by Proposition~\ref{p.law-minimax}, applied with $k = 1$, $t = \be r^2/2$, $h_{1} = h_{-1} = h$, and $a_t = \be r^2/(1+\be r^2)$. This proves parts (1) and (2) of Proposition~\ref{p.structure}. For part (3), the replacement of the constraint $0 \le h \le \be$ by $h \in \R$ was shown in the previous proof. By Theorem~\ref{t.hybrid}, the maximizing radii for $F^{(1)}(s)$, $0 \le s \le \sqrt \be$, lie in a common compact interval, and the function $(r,s) \mapsto \al \cR(r) + B^{(1)}_\al(\be r^2, s)$ is continuous by Proposition~\ref{p.shell}. The supremum in \eqref{e.main-formula} is therefore attained.
\end{proof}

%
%
%
%
%
%
\section{Convergence, and comparison of the two models}
\label{s.consequences}

In this section, we complete the proofs of the results stated in the introduction. We first show that $F_N(\be)$ concentrates around its expectation. This gives the almost sure and $L^1$ convergence in Theorem~\ref{t.hopfield}, which we also show to be uniform over bounded intervals of $\be$. We then prove the assertions of Corollary~\ref{c.one-pattern} concerning $\be \le 1$, and finally the strict separation of Proposition~\ref{p.separation}. We set $F(\al,0) := 0$, consistently with $F_N(0) = 0$.

\begin{proposition}[Concentration and uniform convergence]
\label{p.random-convergence}
Under \eqref{e.loading}, for every $\be > 0$, there exists $C < \infty$ such that
\begin{equation}
\label{e.random-concentration}
\E |F_N(\be) - \bar F_N(\be)| \le C N^{-1/2}.
\end{equation}
Moreover, for every $T < \infty$,
\begin{equation}
\label{e.uniform-beta}
\lim_{N \to \infty} \sup_{\be \in [0,T]} |F_N(\be) - F(\al,\be)| = 0 \qquad \text{almost surely and in } L^1.
\end{equation}
\end{proposition}
\begin{proof}
\emph{Step 1: concentration at fixed $\be$.} We choose $L := L_0(b,\be)$ as in Lemma~\ref{l.norm-truncation}, with $b$ a bound on $M_N/N$, and let $F_N^L := \cF_{N,M_N}(X_N,1;\Ball_{N,M_N}(L))$. This is a convex function of the entries of $X_N$, and its Lipschitz constant for the Euclidean norm is at most $\sqrt \be L/\sqrt N$. Indeed, its gradient is $\sqrt{\be/N} \, \langle \si z^\top \rangle/N$, with $|\si|^2 = N$ and $|z|^2 \le L^2 N$. The inequality \eqref{e.convex-concentration} gives
\begin{equation*}
\P \bigl( |F_N^L - \E F_N^L| > t \bigr) \le C e^{-c N t^2/(\be L^2)}
\qquad \text{and} \qquad
\E |F_N^L - \E F_N^L| \le C_{\be,L} N^{-1/2}.
\end{equation*}
By Lemma~\ref{l.norm-truncation}, the expectation of $F_N - F_N^L$ is exponentially small, and this difference is itself exponentially small outside the event $\{\norm{X_N}_{\mathrm{op}} > \Lambda \sqrt N\}$, whose probability is exponentially small. The triangle inequality proves \eqref{e.random-concentration}. A union bound also shows that $\P(|F_N(\be) - \bar F_N(\be)| > \ep)$ is summable in $N$ for every $\ep > 0$, so by the Borel--Cantelli lemma, $F_N(\be) - \bar F_N(\be) \to 0$ almost surely. Together with \eqref{e.Bernoulli-supAk}, this proves the almost sure and $L^1$ convergence in Theorem~\ref{t.hopfield}.

\emph{Step 2: uniformity in $\be$.} For every realization and every $\be, \be' \ge 0$, we have
\begin{equation}
\label{e.beta-Lipschitz}
|F_N(\be) - F_N(\be')| \le \frac{\norm{X_N}_{\mathrm{op}}^2}{2N} |\be - \be'|.
\end{equation}
By \eqref{e.matrix-norm} and the Borel--Cantelli lemma, the prefactor is bounded almost surely for all $N$ large enough, and its first and second moments are bounded uniformly in $N$. On a single event of probability one, the convergence proved in Step 1 holds at every nonnegative rational $\be$. A finite rational net in $[0,T]$ and \eqref{e.beta-Lipschitz} give the almost sure uniform convergence, and show that the limit is Lipschitz in $\be$. Finally, $\sup_{\be \le T} F_N(\be) \le T \norm{X_N}_{\mathrm{op}}^2/(2N)$ and the uniform bound on the second moment give uniform integrability, which proves the $L^1$ assertion in \eqref{e.uniform-beta}.
\end{proof}

We now turn to the interval $0 < \be \le 1$ on which the two limit free energies coincide. 

\begin{proof}[Proof of \eqref{e.beta-at-most-one} and of the last assertion of Corollary~\ref{c.one-pattern}]
By \eqref{e.field-increment}, for every $r > 0$, $\boldsymbol{\zeta} \in \cP_{<1}$ and $h \in \R$, we have $0 \le \msc P_{\al \be r^2/2}(\boldsymbol{\zeta},h) - \msc P_{\al \be r^2/2}(\boldsymbol{\zeta},0) \le \log \cosh h$.
Adding $\al \cE(\boldsymbol{\zeta})$, taking the infimum over $\boldsymbol{\zeta}$ and then the supremum over $r$ in \eqref{e.gaussian-formula}, we obtain
\begin{equation}
\label{e.Gaussian-field-bound}
0 \le F^\Ga(\al,\be;h) - F^\Ga(\al,\be) \le \log \cosh h.
\end{equation}
For $0 < \be \le 1$ and $h \ne 0$, we have $\log \cosh h < h^2/2 \le h^2/(2\be)$, so by \eqref{e.Gaussian-field-bound}, every $h \ne 0$ gives a strictly smaller value than $h = 0$ in \eqref{e.one-pattern}. This proves \eqref{e.beta-at-most-one} and the assertion on the maximizer.
\end{proof}

Finally, we prove that the two limit free energies differ when $\al$ is small and $\be$ is large. 

\begin{proof}[Proof of Proposition~\ref{p.separation}]
\emph{Step 1: the lower bound by retrieval of one pattern.} For $N$ large enough, we have $M_N \ge 1$, and we restrict the integral in \eqref{e.def.free-energy} to the single configuration $\si = x^1$, which has probability $2^{-N}$ under $P_N$. The first squared overlap $(\sum_i x_i^1 x_i^1)^2$ is $N^2$, and for $\mu \ge 2$, the expectation of $(\sum_i x_i^\mu x_i^1)^2$ is $N$. We thus have
\begin{equation*}
\bar F_N(\be) \ge \frac{\be}{2N^2} \E |X_N^\top x^1|^2 - \log 2 = \frac \be 2 \Ll( 1 + \frac{M_N - 1}{N} \Rr) - \log 2,
\end{equation*}
and Theorem~\ref{t.hopfield} gives the lower bound in \eqref{e.two-bounds}.

\emph{Step 2: the upper bound by a Gaussian comparison.} Let $G$ be an $N \times M_N$ matrix with independent standard Gaussian entries, let $h \sim \mathcal N(0,I_N)$ and $z \sim \mathcal N(0,I_{M_N})$ be independent of each other and of $G$, and let $T_N := N^{-1} \max_{\si \in \{-1,1\}^N} |G^\top \si|$. We compare the centered Gaussian processes $(\si,u) \mapsto \si^\top G u$ and $(\si,u) \mapsto h \cdot \si + \sqrt N z \cdot u$, indexed by $\si \in \{-1,1\}^N$ and $u$ in the unit sphere of $\R^{M_N}$. For two indices $(\si,u)$ and $(\si',u')$, independence, we have
\begin{align*}
\Var\bigl(h \cdot (\si-\si') + \sqrt N z \cdot (u-u')\bigr) - \Var\bigl(\si^\top G u - (\si')^\top G u'\bigr)   = 2 (N - \si \cdot \si') (1 - u \cdot u') \ge 0.
\end{align*}
The Sudakov--Fernique inequality \cite[Theorem~2.2.3]{adler2007random} therefore gives 
\begin{equation}
\label{e.Gaussian-max}
\E T_N \le \frac 1 N \E \sum_{i=1}^N |h_i| + \frac{1}{\sqrt N} \E |z| = \sqrt{2/\pi} + \frac{\E|z|}{\sqrt N} \le \sqrt{2/\pi} + \sqrt{M_N/N}.
\end{equation}

The function $T_N$ is $N^{-1/2}$-Lipschitz in $G$ for the Frobenius norm, since $|\si^\top (G - G') u|/N \le \norm{G - G'}_{\mathrm F}/\sqrt N$, so the Gaussian Poincar\'e inequality gives $\Var(T_N) \le 1/N$. The normalized partition function is at most its largest integrand, so we have (recall $\bar F^{(0)}_{N,M_N}$ from~\eqref{e.def.hybrid})
\begin{equation*}
\bar F^{(0)}_{N,M_N} \le \frac \be 2 \E T_N^2 \le \frac \be 2 \Ll[ \Ll( \sqrt{2/\pi} + \sqrt{M_N/N} \Rr)^2 + \frac 1 N \Rr].
\end{equation*}
Letting $N \to \infty$ proves the upper bound in \eqref{e.two-bounds}.

\emph{Step 3: conclusion.} The two bounds in \eqref{e.two-bounds} give $F(\al,\be) > F^\Ga(\al,\be)$ as soon as $\frac \be 2 (1 + \al) - \log 2 > \frac \be 2 (\sqrt{2/\pi} + \sqrt \al)^2$, that is, as soon as the denominator in \eqref{e.strict-range} is positive and $\be$ satisfies the stated bound. The positivity of this denominator is equivalent to $1 - 2/\pi > 2 \sqrt{2 \al/\pi}$, which is the stated bound on $\al$.
\end{proof}

\medskip

\small
\bibliographystyle{plain}
\bibliography{ref}

\begin{thebibliography}{10}

\bibitem{adler2007random}
Robert~J. Adler and Jonathan~E. Taylor.
\newblock {\em Random fields and geometry}.
\newblock Springer Monographs in Mathematics. Springer, New York, 2007.

\bibitem{agliari2026infinite}
Elena Agliari, Adriano Barra, Pierluigi Bianco, Alberto Fachechi, and Diego
  Pallara.
\newblock The infinite volume limit of the {Hopfield} neural networks in the
  replica-symmetric regime.
\newblock {\em Ann. Henri Poincar\'{e}}, 2026.
\newblock Online first.

\bibitem{agliari2017neural}
Elena Agliari, Adriano Barra, Chiara Longo, and Daniele Tantari.
\newblock Neural networks retrieving {Boolean} patterns in a sea of {Gaussian}
  ones.
\newblock {\em J. Stat. Phys.}, 168:1085--1104, 2017.

\bibitem{amit1985spin}
Daniel~J. Amit, Hanoch Gutfreund, and Haim Sompolinsky.
\newblock Spin-glass models of neural networks.
\newblock {\em Phys. Rev. A}, 32(2):1007--1018, 1985.

\bibitem{amit1985storing}
Daniel~J. Amit, Hanoch Gutfreund, and Haim Sompolinsky.
\newblock Storing infinite numbers of patterns in a spin-glass model of neural
  networks.
\newblock {\em Phys. Rev. Lett.}, 55(14):1530--1533, 1985.

\bibitem{amit1987statistical}
Daniel~J. Amit, Hanoch Gutfreund, and Haim Sompolinsky.
\newblock Statistical mechanics of neural networks near saturation.
\newblock {\em Ann. Physics}, 173(1):30--67, 1987.

\bibitem{aronow2026quantitative}
Peter~M. Aronow and Patrick Lopatto.
\newblock Quantitative {Parisi} formulas and fluctuations in the
  {Sherrington--Kirkpatrick} model.
\newblock {\em Preprint, arXiv:2609.09103}, 2026.

\bibitem{auffinger2015unique}
Antonio Auffinger and Wei-Kuo Chen.
\newblock The {Parisi} formula has a unique minimizer.
\newblock {\em Comm. Math. Phys.}, 335(3):1429--1444, 2015.

\bibitem{barra2012equivalence}
Adriano Barra, Alberto Bernacchia, Enrica Santucci, and Pierluigi Contucci.
\newblock On the equivalence of {Hopfield} networks and {Boltzmann} machines.
\newblock {\em Neural Netw.}, 34:1--9, 2012.

\bibitem{barra2015multispecies}
Adriano Barra, Pierluigi Contucci, Emanuele Mingione, and Daniele Tantari.
\newblock Multi-species mean field spin glasses. {R}igorous results.
\newblock {\em Ann. Henri Poincar\'{e}}, 16(3):691--708, 2015.

\bibitem{barra2013mean}
Adriano Barra, Gino Del~Ferraro, and Daniele Tantari.
\newblock Mean field spin glasses treated with {PDE} techniques.
\newblock {\em Eur. Phys. J. B}, 86(7):Art. 332, 10, 2013.

\bibitem{barra2010breaking}
Adriano Barra, Aldo Di~Biasio, and Francesco Guerra.
\newblock Replica symmetry breaking in mean-field spin glasses through the
  {Hamilton--Jacobi} technique.
\newblock {\em J. Stat. Mech. Theory Exp.}, 2010(9):P09006, 22, 2010.

\bibitem{barra2026existence}
Adriano Barra and Fabrizio Durante.
\newblock Existence of the thermodynamic limit for the {Hopfield} model.
\newblock {\em Preprint, arXiv:2610.41938}, 2026.

\bibitem{barra2010replica}
Adriano Barra, Giuseppe Genovese, and Francesco Guerra.
\newblock The replica symmetric approximation of the analogical neural network.
\newblock {\em J. Stat. Phys.}, 140(4):784--796, 2010.

\bibitem{barra2012glassy}
Adriano Barra, Giuseppe Genovese, Francesco Guerra, and Daniele Tantari.
\newblock How glassy are neural networks?
\newblock {\em J. Stat. Mech. Theory Exp.}, 2012(7):P07009, 2012.

\bibitem{barra2018phase}
Adriano Barra, Giuseppe Genovese, Peter Sollich, and Daniele Tantari.
\newblock Phase diagram of restricted {Boltzmann} machines and generalized
  {Hopfield} networks with arbitrary priors.
\newblock {\em Phys. Rev. E}, 97(2):022310, 2018.

\bibitem{barra2008ergodic}
Adriano Barra and Francesco Guerra.
\newblock About the ergodic regime in the analogical {Hopfield} neural
  networks: moments of the partition function.
\newblock {\em J. Math. Phys.}, 49(12):125217, 2008.

\bibitem{boursier2026free}
Jeanne Boursier and Jean-Christophe Mourrat.
\newblock Free energy of non-convex multi-species {Ising} spin glasses with
  external fields.
\newblock {\em Preprint, arXiv:2610.18350}, 2026.

\bibitem{bovier1997retrieval}
Anton Bovier and V\'{e}ronique Gayrard.
\newblock The retrieval phase of the {Hopfield} model: a rigorous analysis of
  the overlap distribution.
\newblock {\em Probab. Theory Related Fields}, 107(1):61--98, 1997.

\bibitem{bovier1998review}
Anton Bovier and V\'{e}ronique Gayrard.
\newblock {Hopfield} models as generalized random mean field models.
\newblock In Anton Bovier and Pierre Picco, editors, {\em Mathematical aspects
  of spin glasses and neural networks}, volume~41 of {\em Progr. Probab.},
  pages 3--89. Birkh\"{a}user Boston, Boston, MA, 1998.

\bibitem{bovier1995gibbs}
Anton Bovier, V\'{e}ronique Gayrard, and Pierre Picco.
\newblock Gibbs states of the {Hopfield} model with extensively many patterns.
\newblock {\em J. Statist. Phys.}, 79(1--2):395--414, 1995.

\bibitem{camilli2026solution}
Francesco Camilli, Pierluigi Contucci, and Emanuele Mingione.
\newblock The solution of the {Hopfield} model with {Gaussian} patterns.
\newblock {\em Preprint, arXiv:2610.32618}, 2026.

\bibitem{carmona2006universality}
Philippe Carmona and Yueyun Hu.
\newblock Universality in {Sherrington--Kirkpatrick}'s spin glass model.
\newblock {\em Ann. Inst. H. Poincar\'e Probab. Statist.}, 42(2):215--222,
  2006.

\bibitem{chatterjee2005simple}
Sourav Chatterjee.
\newblock A simple invariance theorem.
\newblock {\em Preprint, arXiv:math/0508213}, 2005.

\bibitem{chen2026species}
Hong-Bin Chen.
\newblock On free energy of non-convex multi-species spin glasses.
\newblock {\em ALEA Lat. Am. J. Probab. Math. Stat.}, 23:429--473, 2026.

\bibitem{chen2026formula}
Hong-Bin Chen, Fu-Hsuan Ho, and Victor Issa.
\newblock A formula for the free energy of non-convex vector spin glasses.
\newblock {\em Preprint, arXiv:2610.18744}, 2026.

\bibitem{chen2026ising}
Hong-Bin Chen, Victor Issa, and Jean-Christophe Mourrat.
\newblock Free energy of non-convex multi-species spin glasses with centered
  {Ising} spins.
\newblock {\em Preprint, arXiv:2606.16636}, 2026.

\bibitem{chen2025vector}
Hong-Bin Chen and Jean-Christophe Mourrat.
\newblock On the free energy of vector spin glasses with non-convex
  interactions.
\newblock {\em Probab. Math. Phys.}, 6:1--80, 2025.

\bibitem{chen2026spherical}
Hong-Bin Chen and Jean-Christophe Mourrat.
\newblock Free energy of non-convex multi-species spherical spin glasses.
\newblock {\em Preprint, arXiv:2609.32710}, 2026.

\bibitem{chen2025HJ}
Hong-Bin Chen and Jiaming Xia.
\newblock {Hamilton--Jacobi} equations from mean-field spin glasses.
\newblock {\em Probab. Theory Related Fields}, 192(3--4):803--873, 2025.

\bibitem{crisanti1986saturation}
Andrea Crisanti, Daniel~J. Amit, and Hanoch Gutfreund.
\newblock Saturation level of the {Hopfield} model for neural network.
\newblock {\em Europhys. Lett.}, 2(4):337--341, 1986.

\bibitem{guerra2001sum}
Francesco Guerra.
\newblock Sum rules for the free energy in the mean field spin glass model.
\newblock In {\em Mathematical physics in mathematics and physics ({Siena},
  2000)}, volume~30 of {\em Fields Inst. Commun.}, pages 161--170. Amer. Math.
  Soc., Providence, RI, 2001.

\bibitem{hebb1949organization}
Donald~O. Hebb.
\newblock {\em The organization of behavior: a neuropsychological theory}.
\newblock Wiley, New York, 1949.

\bibitem{ho2026convergence}
Fu-Hsuan Ho.
\newblock Convergence of the free energy of vector spin glasses with non-convex
  interactions.
\newblock {\em Preprint, arXiv:2609.39307}, 2026.

\bibitem{hopfield1982neural}
John~J. Hopfield.
\newblock Neural networks and physical systems with emergent collective
  computational abilities.
\newblock {\em Proc. Natl. Acad. Sci. U.S.A.}, 79(8):2554--2558, 1982.

\bibitem{jagannath2016dynamic}
Aukosh Jagannath and Ian Tobasco.
\newblock A dynamic programming approach to the {Parisi} functional.
\newblock {\em Proc. Amer. Math. Soc.}, 144(7):3135--3150, 2016.

\bibitem{kim2026sharp}
Taegyun Kim.
\newblock A sharp universality dichotomy for the free energy of spherical spin
  glasses.
\newblock {\em Preprint, arXiv:2601.08599}, 2026.

\bibitem{little1974existence}
William~A. Little.
\newblock The existence of persistent states in the brain.
\newblock {\em Math. Biosci.}, 19(1--2):101--120, 1974.

\bibitem{mourrat2021nonconvex}
Jean-Christophe Mourrat.
\newblock Nonconvex interactions in mean-field spin glasses.
\newblock {\em Probab. Math. Phys.}, 2(2):281--339, 2021.

\bibitem{mourrat2022parisi}
Jean-Christophe Mourrat.
\newblock The {Parisi} formula is a {Hamilton--Jacobi} equation in
  {Wasserstein} space.
\newblock {\em Canad. J. Math.}, 74(3):607--629, 2022.

\bibitem{mourrat2023upper}
Jean-Christophe Mourrat.
\newblock Free energy upper bound for mean-field vector spin glasses.
\newblock {\em Ann. Inst. Henri Poincar\'{e} Probab. Stat.}, 59(3):1143--1182,
  2023.

\bibitem{panchenko2015multispecies}
Dmitry Panchenko.
\newblock The free energy in a multi-species {Sherrington--Kirkpatrick} model.
\newblock {\em Ann. Probab.}, 43(6):3494--3513, 2015.

\bibitem{pastur1994replica}
Leonid Pastur, Mariya Shcherbina, and Brunello Tirozzi.
\newblock The replica-symmetric solution without replica trick for the
  {Hopfield} model.
\newblock {\em J. Statist. Phys.}, 74(5--6):1161--1183, 1994.

\bibitem{pastur1977rigorously}
Leonid~A. Pastur and Alexander~L. Figotin.
\newblock Rigorously solvable model of spin glass.
\newblock {\em JETP Lett.}, 25(8):323--326, 1977.

\bibitem{pastur1978theory}
Leonid~A. Pastur and Alexander~L. Figotin.
\newblock Theory of disordered spin systems.
\newblock {\em Theoret. and Math. Phys.}, 35(2):403--414, 1978.

\bibitem{sawhney2024free}
Mehtaab Sawhney and Mark Sellke.
\newblock Free energy universality of spherical spin glasses.
\newblock {\em Preprint, arXiv:2408.13701}, 2024.

\bibitem{shcherbina1993free}
Mariya Shcherbina and Brunello Tirozzi.
\newblock The free energy of a class of {Hopfield} models.
\newblock {\em J. Statist. Phys.}, 72(1--2):113--125, 1993.

\bibitem{sion1958general}
Maurice Sion.
\newblock On general minimax theorems.
\newblock {\em Pacific J. Math.}, 8(1):171--176, 1958.

\bibitem{steffan1994replica}
Helmut Steffan and Reimer K\"{u}hn.
\newblock Replica symmetry breaking in attractor neural network models.
\newblock {\em Z. Phys. B}, 95(2):249--260, 1994.

\bibitem{subag2026geometry}
Eliran Subag.
\newblock Geometry of spherical spin glasses.
\newblock In Susan Friedlander and Yuri Tschinkel, editors, {\em Proceedings of
  the International Congress of Mathematicians 2026, Volume 6: Invited Lectures
  (Sections 12--14)}, pages 102--122. SIAM, 2026.

\bibitem{talagrand1995concentration}
Michel Talagrand.
\newblock Concentration of measure and isoperimetric inequalities in product
  spaces.
\newblock {\em Inst. Hautes \'{E}tudes Sci. Publ. Math.}, 81:73--205, 1995.

\bibitem{talagrand1998rigorous}
Michel Talagrand.
\newblock Rigorous results for the {Hopfield} model with many patterns.
\newblock {\em Probab. Theory Related Fields}, 110(2):177--276, 1998.

\bibitem{talagrand2000exponential}
Michel Talagrand.
\newblock Exponential inequalities and convergence of moments in the
  replica-symmetric regime of the {Hopfield} model.
\newblock {\em Ann. Probab.}, 28(4):1393--1469, 2000.

\bibitem{talagrand2011mean1}
Michel Talagrand.
\newblock {\em Mean field models for spin glasses. {V}olume {I}: Basic
  examples}, volume~54 of {\em Ergebnisse der Mathematik und ihrer
  Grenzgebiete. 3. Folge}.
\newblock Springer-Verlag, Berlin, 2011.

\bibitem{talagrand2011mean2}
Michel Talagrand.
\newblock {\em Mean field models for spin glasses. {V}olume {II}: Advanced
  replica-symmetry and low temperature}, volume~55 of {\em Ergebnisse der
  Mathematik und ihrer Grenzgebiete. 3. Folge}.
\newblock Springer, Heidelberg, 2011.

\bibitem{vershynin2012introduction}
Roman Vershynin.
\newblock Introduction to the non-asymptotic analysis of random matrices.
\newblock In Yonina~C. Eldar and Gitta Kutyniok, editors, {\em Compressed
  sensing: theory and applications}, pages 210--268. Cambridge University
  Press, 2012.

\end{thebibliography}

\end{document}